\documentclass[12pt]{article}

\usepackage{lineno}
\usepackage[hidelinks]{hyperref}
\modulolinenumbers[5]
\usepackage{amssymb}
\usepackage{amsmath}
\usepackage[square,numbers]{natbib}

\newcommand{\BYDEF}{\,\shortstack{\textup{\tiny def} \\ = }\,}

\def\o{\over}

\def\G{\Gamma}
\def\g{\gamma}
\def\gm{\gamma(dy,du)}
\def\xm{\xi(dy,du)}
\def\k{\kappa}
\def\a{\alpha}
\def\d{\delta}
\def\e{\epsilon}

\def\A{{\mathcal A}}

\def\P{{\mathcal P}}
\def\D{{\mathcal D}}
\def\M{{\mathcal M}}
\def\Mp{{\mathcal M_+}}

\def\reals{I\!\!R}
\def\O{\Omega}
\def\ph{\varphi}

\def\bp{\bar \psi}

\def\hf{\hfill{$\Box$}}
\def\T{\Theta}

\allowdisplaybreaks

\newtheorem{Theorem}{Theorem}[section]
\newtheorem{Proposition}[Theorem]{Proposition}

\newtheorem{Lemma}[Theorem]{Lemma}
\newtheorem{Corollary}[Theorem]{Corollary}

\numberwithin{equation}{section}

\begin{document}

{\bf \Large{LP Based Upper and Lower Bounds for Ces\`aro and Abel Limits of the Optimal Values in Problems of Control of Stochastic Discrete Time Systems\footnote{This is the author version of the paper of the same name which appears in Journal of Mathematical Analysis and Applications, v. 512(1), 2022}}}
\bigskip
\bigskip
\bigskip

\setcounter{footnote}{0}
\renewcommand*{\thefootnote}{\fnsymbol{footnote}}
{\bf K. Avrachenkov$^a$, V. Gaitsgory$^b$ and L. Gamertsfelder$^b$}\\
$^a$ {\it\small{INRIA Sophia Antipolis, 2004 Route des Lucioles, 06902 Valbonne, France} }\\
$^b$ {\it \small{Department of Mathematics and Statistics, Macquarie University, Eastern Road, Macquarie Park, NSW 2113, Australia}}
\renewcommand*{\thefootnote}{\arabic{footnote}}
{\let\thefootnote\relax\footnote{{Email addresses: \href{mailto:k.avrachenkov@inria.fr}{k.avrachenkov@inria.fr} (Konstantin Avrachenkov), \href{mailto:vladimir.gaitsgory@mq.edu.au}{vladimir.gaitsgory@mq.edu.au} (Vladimir Gaitsgory), \href{mailto:lucas.gamertsfelder@gmail.com}{lucas.gamertsfelder@gmail.com} (Lucas Gamertsfelder)}}}

\bigskip
\bigskip
\bigskip
	
{\bf Abstract.} \small{
		In this paper, we study asymptotic properties of    problems of control of stochastic discrete time systems (also known as Markov decision processes) with time averaging and time discounting optimality criteria, and we establish that the Ces\`aro and Abel limits of the optimal values in such problems can be evaluated with the help of a certain infinite-dimensional linear programming  problem and its dual.}

%
%
%
%
%

\section{Introduction}

In this paper, we analyze optimal values in problems of control of stochastic discrete time systems
considered on long or infinite time horizons. We work with stochastic discrete time systems presented in
the form of controlled stochastic recursions. It is known (see e.g., \cite{BaRi,BF1992,Kifer})
that a controlled stochastic recursion can be represented as a Markov Decision Process (MDP) and vice versa.
It is not rare that a different point of view provides new insights into a well studied research area. Viewing MDPs
as controlled stochastic recursions allows us to gain a better understanding of MDP properties
in the challenging non-ergodic case.

There are many approaches to study  optimal control of stochastic discrete time systems considered on long
or infinite time horizons. Let us mention three such approaches that  we are using in the present work. Firstly,
it is the optimization  of the expected discounted cost over the  infinite time horizon.
The expected discounted cost with a fixed discount factor can, however, be very insensitive to
the long-term effects of a control, and
 one of common ways to study the long-term MDPs
is to consider the lower or upper limit of the optimal expected discounted cost when the discount factor tends to one \cite{Feinbergatal2012,Jean-H-2,Jean-H-3,Filar92}.
These are called the  lower or, respectively, upper Abel limits of the optimal expected discounted cost.

Secondly it is also quite common
in the MDP theory, to consider the lower or upper limits of the averages of the expected
costs, the so-called lower or upper Ces\`aro limits, and then optimize these quantities.
In the MDP literature (see e.g., \cite{BorkarSurvey1993,Jean-H-2}
and references therein), this approach is referred to as a long-run average cost problem.
However, as was noted in \cite{Flynn1980}, in contrast to the discounted criterion, the long-run
average criterion disregards system dynamics on any finite length time interval. Therefore,
it was proposed in \cite{Flynn1980} to consider a sequence of optimal average values over finite
horizons with increasing lengths. This is the third optimality criterion that we study in this
work. Note that the Ces\`aro limit of the sequence of optimal average values may not exist in the general
case, and one needs to investigate the lower and upper Ces\`aro limits of these sequences.

In fact, there are many more MDP optimality criteria such as average overtaking optimality,
selective optimality criteria, etc. Note that in general different long-term optimality criteria can give
different values or may not produce an optimal policy in some settings \cite{Feinberg1980,Khlopin,Sorin92,LM1994,Renault09}.
For an extensive account of various MDP optimality criteria and a historic perspective we refer
an interested reader to the books \cite{Jean-H-2,Jean-H-3,KallenbergLN,Piun-1997,Piun-2013,Puterman}.

Linear programming is one of the main tools for both theoretical and numerical analysis of MDPs.
Various optimality criteria have been extensively studied by linear programming and other tools
under conditions ensuing ergodicity or communicability of the underlying processes. For overview of many works
in the ergodic setting we refer to \cite{BorkarSurvey1993,Jean-H-2,Puterman}. The case with the presence of several non-communicating ergodic classes is more challenging and has received much less attention, particularly,
in infinite state spaces. The authors of \cite{DenardoFox1968,Denardo1970} have proposed to apply
nested linear programs together with combinatorial search procedures for finding average cost optimal
controls in finite MDPs. Then in \cite{HK-1,HK-2,Kallenberg1983}, also for finite MDPs, the authors
have constructed a single linear programming formulation with two layers of constraints.
In \cite{HL1994} the two-layer linear programming approach to average cost control problem
has been extended to countable state spaces with multi-chain ergodic structure. There the authors
indicated that there may exist a duality gap between the primal and dual linear programs.
In \cite{Gonzalez-Hernandes}, the two-layer linear programming approach has been extended to MDPs
with  general (non-compact) Borel state spaces, and sufficient conditions for the strong duality in the corresponding two-layer linear programs have been studied.
Note that, as has been demonstrated in the recent work \cite{Ilya-arch}, the duality gap can appear in such program even in the case when
the MDPs have compact state spaces.

In the first part of our work, we show that the convex closed hull of the discounted occupational measures set (obtained with fixed initial conditions) is characterized by  a set of linear constraints depending on the parameter $\e$, with $1-\e $ being the discount factor, and that the optimal expected discounted value is equal to
the optimal value of an Infinite Dimensional Linear Program (IDLP) considered on this set.
Also in the first part, we establish that, as the time horizon $T$ tends to infinity, the closed convex
hull of the union (over all initial conditions) of the occupational measure sets converges to the set of stationary measures defined by linear constraints, and the closed convex hull of the union of the sets of  discounted occupational measures  converges to this set as $\e $ tends to zero (see Theorem \ref{Th-discounting-1}). Let us note that already in \cite{Manne1960} it was conjectured and in \cite{Denardo1970}
rigorously proved  for finite models that only one-layer linear
programming formulation is sufficient if one seeks to optimize over both the set of policies
and the set of initial states.

Since, in the present work, we use the controlled stochastic recursion representation of MDP as a point
of departure, the constraints of our primal programs are defined in terms of continuous test functions,
which, at the first sight, looks very different from the form of constraints in \cite{Jean-H-2,Gonzalez-Hernandes}.
We note, however,  that the constraints can be reformulated in terms of bounded Borel test functions.
The latter are readily shown to be equivalent to the constraints in a more common form based on the use of the transition probability kernel. As  in \cite{Jean-H-2,Jean-H-3}, we also show that in the case of joint optimization over
the sets of admissible policies and initial states, there is no duality gap,
which in particular implies that both Abel limit of the expected discounted costs and
the Ces\`aro limit of the optimal values over increasing finite horizons exist and coincide.
One difference with respect to the results obtained in \cite{Jean-H-2,Jean-H-3} is that we consider
the convergence of the optimal values as the time horizon moves to infinity while in \cite{Jean-H-2,Jean-H-3}
the long-term average cost problem is considered.

The most important results of the paper are obtained in the second part of the paper, in which we consider the setting with a fixed initial state and view
the linear program for the expected discounted cost as a singularly perturbed linear program
\cite{Kostia-0,Per-Gai}, with a perturbation parameter $\e $.
Specifically, the IDLP problem considered on the set of stationary measures is obtained by formally taking $\e =0$ in the perturbed problem. However, its optimal value can be strictly less than the limit of the optimal value of the perturbed IDLP problem, thus allowing one to characterize the latter as singularly perturbed (SP) problem. A way of dealing with SP problems is by augmenting the set of constraints obtained by just taking $\e =0$ with some additional constraints and decision variables (see \cite{Kostia-0,Kostia-1,Kostia-2,GR-1,Per-Gai}).
The augmented linear problem and its dual provide upper and lower bounds for the lower and upper Ces\`aro and Abel limits of the optimal values (see Theorem \ref{Th-lower-upper-limit-main}). In addition, we establish that, if a point-wise limit of the optimal values as the horizon moves to infinity exists and is continuous with respect to the initial state, then this limit is equal to the optimal value of the augmented dual problem, and also that, if a point-wise limit of the expected discounted costs as the discount factor goes to one exists and is continuous, then this limit is equal to the optimal value of the augmented dual problem too.

We  introduced an optimization problem, the optimal value of which is shown to be equal to that of the dual problem. The  feasible domain of this problem contains the feasible domain of the primal IDLP problem, and the sufficient condition for the strong duality is that the closure of the latter is equal to the former (see Proposition \ref{Prop-cl-D-y-0-dual}
and Corollary \ref{Cor-the-only-conclusion}). We have  shown that the upper limit of the   set of the occupational measures generated by the state-control trajectories emanated from a given initial point $y_0$ is contained in the feasible set of this  problem (see Proposition \ref{Lemma-conv-to cl-D-y-0}), with the limit of the  convex closed hull of the former being equal to the latter under the additional assumption that the pointwise limits of the optimal value functions exist and are continuous (see Corollary \ref{Cor-conv-feasible}).
We have also
provided a condition for the absence of duality gap when the time average optimality can be achieved with the use of  periodic generating regimes  (see Proposition \ref{Prop-PGR} and Corollary \ref{Cor-strong-dual-1}).
The importance of the fact that the limit of the optimal values is (or can be) equal
to the optimal value  of the augmented dual problem is emphasized by sufficient and necessary conditions for the long run average optimality (see Propositions \ref{Prop-optim-cond-suf} and \ref{Prop-nec-opt-cond}).

The approach that we use in the present paper was proposed in  \cite{BG-2019}, where LP based upper and lower bounds for  Ces\`aro and Abel limits of
optimal values  were obtained for \lq\lq non-ergodic" deterministic systems evolving in continuous time. In \cite{BGS-2019}, similar results were obtained
for deterministic systems evolving in discrete time under less restrictive assumptions than in \cite{BG-2019}. For example, in contrast to \cite{BG-2019}, no assumptions about the existence of the uniform Ces\`aro and Abel limits and their Lipschitz continuity as functions of the initial values have been made in \cite{BGS-2019} to establish the upper bound for these limits. Many results obtained in this paper are stated and proved similarly to their deterministic counterparts obtained in \cite{BGS-2019}.
However,
the part of Theorem \ref{Th-lower-upper-limit-main}  that establishes the upper bound for the  Ces\`aro and Abel limits of the optimal values (this being one of the most important results of the paper) is stronger than
the corresponding statement  in \cite{BGS-2019} (see Theorem 3.1 in \cite{BGS-2019}) since, in contrast to the latter, it is not assumed that the optimal value over the finite horizon and the optimal expected discounted cost are continuous
with respect to the initial state. (Note that the line of research started in \cite{BG-2019} has been also continued in \cite{GS-2021-1} and  \cite{GS-2021-2}, where some results for deterministic  continuous time systems
related to  relaxations of  conditions   used  in \cite{BG-2019}  have been obtained.)

The paper consists of ten sections. In Section \ref{Model-Prelim} we first define the model and state the assumptions. Then, we recall definitions of occupational (and discounted occupational) measures generated by the state-control trajectories of system, and we reformulate optimal control problems as problems of optimization on these sets. All obtained results are stated in Sections~\ref{Sec-Time-discounting} - \ref{Sec-opt-cond}, and most of their proofs (as well as the results needed for these proofs) are given in Sections \ref{Sec-Proof-Ab-Ch-estimates} - \ref{Sub-Sec-Auxiliary-Res}.

\section{Model and Preliminaries}
\label{Model-Prelim}

We consider the discrete time stochastic control system in a form of controlled stochastic recursion
\begin{equation}\label{A1}
y(t+1) = f(y(t), u(t), s(t)), \quad t  = 0, 1, 2, \ldots ,
\end{equation}
and we assume  that the following conditions are satisfied everywhere in what follows:
\begin{itemize}
	\itemsep0em
	\item The function $f(y,u,s) : Y\times \hat{U}  \times S \rightarrow \mathbb{R}^m$ is bounded, continuous in $(y,u) $ on $Y\times \hat{U} $ and  Borel measurable  in $s$ on $S$, where $Y$ is a compact subset of $\mathbb{R}^n$, $\hat{U}$ is a compact metric space,  and $S$ is a Polish space.

	\item  $s(t) \in S, \ t=0,1,...,$ is a sequence of independent, identically distributed random elements defined on a common probability space.
	
\item The controls $u(t),\ t=0,1,..., $ are defined by a sequence of   functions $\pi\BYDEF\{\pi_t(y),\  t=0,1,...\} $ that are Borel measurable selections
 of a multivalued map $U(\cdot): Y \leadsto \hat{U}$
 so that

	\begin{equation}\label{e-feedback-control}
	u(t)= \pi_t(y(t))\in U(y(t)),\ \  t=0,1,...,
	\end{equation}
where $U(\cdot)$ is upper semicontinuous and compact-valued (that is, $U(y) $ is compact for any $y\in Y $).
	
	\item $f(y,u,s) \in Y$ for any $y \in Y$, any $u \in U(y)$,  and any $s \in S$ (that is, the set $Y$ is forward invariant with respect to system (\ref{A1})).
	
\end{itemize}
Let $\Pi $ stand for the set of sequences of measurable selections of $U(\cdot) $:
$$
\Pi \BYDEF \left\{\pi = \{\pi_t(\cdot),\ t = 0,1,...\}\ | \ \pi_t(y)\in U(y)\  \forall \ y\in Y, \ \ \pi_t(\cdot) \ {\rm are \ Borel\ measurable}\right\}.
$$
For any sequence $\pi\in \Pi$ (for convenience, such sequences  will be referred to as {\it control plans}) and any initial condition $\ y(0)=y_0\in Y $, let $(y^{\pi, y_0}(\cdot), u^{\pi,  y_0}(\cdot)) $ stand for the state-control trajectory obtained in accordance with (\ref{A1}) and (\ref{e-feedback-control}).

Consider the following optimal control problems
\begin{equation}\label{A112-1}
{1\o T} \min_{\pi\in \Pi}E\left[\sum_{t=0}^{T-1} k(y^{\pi, y_0}(t),u^{\pi, y_0}(t))\right]\BYDEF v_T(y_0),
\end{equation}
\begin{equation}\label{A112-2}
\e\min_{\pi\in\Pi}\sum_{t=0}^{\infty}(1-\e)^t k(y^{\pi , y_0}(t),u^{\pi , y_0}(t))\BYDEF h_{\e}(y_0),
\end{equation}
where $k(y,u):Y\times \hat{U}\rightarrow \mathbb{R} $ is a continuous function with
\begin{equation}\label{intergand-bounded}
	|k(y,u)|\leq M \ \ \ \forall (y,u)\in Y\times \hat{U}
, \ \ \ \ M = {\rm const},
	\end{equation}
and $\e \in (0,1)$ (that is, $(1-\e)$ is a discount factor).
Note that, as follows from the proposition stated below,
 the minima in (\ref{A112-1})  and (\ref{A112-2}) are achieved and the optimal value functions $v_T(\cdot)$ and $h_{\e}(\cdot) $ are
lower semicontinuous.

\begin{Proposition}\label{Prop-lsc}
Under the assumptions made above, the optimal value functions $v_T(\cdot)$ and $h_{\e}(\cdot) $ are
lower semicontinuous and satisfy the equations:
\begin{equation}\label{e-App-1}
T v_T(y)= \min_{u\in U(y)}\{k(y,u) + (T-1)E[v_{T-1}(f(y,u,s))]\} \ \ \forall \ y\in Y, \ \ \ T=1,2,...,
\end{equation}
\begin{equation}\label{e-App-2}
h_{\e}(y)= \min_{u\in U(y)}\{ \e k(y,u) + (1-\e) E[h_{\e}(f(y,u,s))]\} \ \ \forall \ y\in Y.
\end{equation}
Also, the minima in (\ref{A112-1})  and (\ref{A112-2}) are achieved.
\end{Proposition}
{\bf Proof.} The proof follows Theorems 2.4.6 and 7.2.1 in \cite{BaRi}. $\ \Box$

We will be interested in evaluating
$\lim_{T\rightarrow\infty} v_T(y_0)$ and $\lim_{\e\rightarrow 0} h_{\e}(y_0)$  (these limits are commonly referred to as the {\it Ces\`aro limit} of the sequence of optimal values and the {\it Abel limit} of the expected discounted costs). More specifically, we will establish that $\limsup_{T\to \infty} v_T(y_0)$ and $ \limsup_{\e\rightarrow 0} h_{\e}(y_0)$ are bounded from above by the optimal value of a certain IDLP problem, and that  $\liminf_{T\to \infty} v_T(y_0)$ and $ \liminf_{\e\rightarrow 0} h_{\e}(y_0)$ are bounded from below by the optimal value of the corresponding dual problem (see Theorem \ref{Th-lower-upper-limit-main}). An immediate consequence of this result is the statement that the Ces\`aro and Abel limits exist and are equal to each other if there is no duality gap (see Corollary \ref{Cor-Strong-dual}).

Let us  conclude this section with the introduction of some notations and definitions that will be used in the subsequent sections. Let $G$ stand for the graph of  $U(\cdot)$,
\begin{equation*}
G \BYDEF \text{graph}(U) = \left\{(y,u) : u \in U(y),\text{ }  y \in Y \right\},
\end{equation*}
and let $\mathcal{P}(G)$ stand for the set of probability measures defined on Borel subsets of $G$. (Note that, due to upper semicontinuity of  $U(\cdot)$, the graph $G$ is a compact subset of $Y\times \hat U $.) 
Given a control plan $\pi\in \Pi $ and an initial condition $y(0)=y_0\in Y $,  denote by $\g^{\pi , y_0 ,T}\in \mathcal{P}(G) $ and $\g_d^{\pi , y_0  , \e}\in \mathcal{P}(G) $ the probability measures defined as follows: for
any Borel $Q\subset G $,
\begin{equation}\label{E6-taime-ave}
\g^{\pi , y_0 ,T}(Q)={1\o T} E\left[\sum_{t=0}^{T-1} 1_Q(y^{\pi , y_0 }(t),u^{\pi , y_0 }(t))\right],
\end{equation}
\begin{equation}\label{E6}
\g_d^{\pi,  y_0 , \e}(Q)=\e E\left[\sum_{t=0}^{\infty} (1-\e)^t 1_Q(y^{\pi , y_0 }(t),u^{\pi , y_0 }(t))\right],
\end{equation}
where  $1_Q(\cdot)$ is the indicator function of $Q$. The measures defined by (\ref{E6-taime-ave}) and (\ref{E6}) will be referred to as {\it occupational measure} and, respectively, {\it discounted occupational measure} generated by the control plan $\pi$. Note that from (\ref{E6-taime-ave}) and (\ref{E6}) it follows that
\begin{equation}\label{G88}
\int_{G} q(y,u) \g^{\pi ,  y_0 ,T}(dy,du)={1\o T}E \left[\sum_{t=0}^{T-1} q(y^{\pi , y_0 }(t),u^{\pi , y_0 }(t))\right]
\end{equation}
and
\begin{equation}\label{G8}
\int_{G} q(y,u) \g_d^{\pi ,  y_0,\e}(dy,du)=\e E \left[\sum_{t=0}^{\infty} (1-\e)^t q(y^{\pi , y_0 }(t),u^{\pi , y_0 }(t))\right]
\end{equation}
for any bounded Borel measurable function $q$ on $G$.  In fact, the definitions (\ref{E6-taime-ave}) and (\ref{E6}) are equivalent to that the equality
(\ref{G88}) and (\ref{G8}) are valid if $q(\cdot)$ is an indicator function of $Q$. Therefore, these equalities are valid for linear combinations of indicator functions. The validity of (\ref{G88}) and (\ref{G8}) for any Borel function  follows from the fact that any such function can be presented as uniform limit of linear combinations of indicator functions.

Let us denote by  $\G_T(y_0) $ the set of occupational measures and by $\T_{\e}(y_0) $ the set of discounted occupational measures:
\begin{equation}\label{union-1}
\G_T(y_0)\BYDEF\bigcup_{\pi\in\Pi}\{\g^{\pi ,  y_0 ,T}\},\ \ \ \ \ \ \  \T_{\e}(y_0)\BYDEF\bigcup_{_{\pi\in\Pi}}\{\g_d^{\pi ,  y_0,\e}\}.
\end{equation}
Note that, due to \eqref{G88} and \eqref{G8}, problems \eqref{A112-1} and (\ref{A112-2}) can be rewritten in the form
\begin{equation}\label{A3-T}
\min_{\g\in \G_T(y_0)} \int_{G} k(y,u)\g(dy,du)=v_T(y_0)
\end{equation}
and
\begin{equation}\label{A3}
\min_{\g\in  \T_{\e}(y_0)} \int_{G} k(y,u)\g(dy,du)=h_{\e}(y_0),
\end{equation}
respectively.

To describe convergence properties of occupational measures,  the following metric on
$\P(G)$ will be used:
\begin{equation}\label{e-rho}
\rho(\g',\g''):=\sum_{j=1}^{\infty} {1\o 2^j}\left|\int_G q_j(y,u)\g'(dy,du)-\int_G q_j(y,u)\g''(dy,du)\right|
\end{equation}
for $\g',\g''\in \P(G)$, where $q_j(\cdot),\,j=1,2,\dots,$ is a sequence of Lipschitz continuous functions dense in the unit ball of the space of continuous functions $C(G)$ from $G$ to $\mathbb{R}$.
This metric is consistent with the weak$^*$ convergence topology on $\P(G)$, that is,
a sequence $\g^k\in \P(G)$ converges to $\g\in \P(G)$ in this metric if and only if
$$
\lim_{k\to \infty}\int_G q(y,u)\g^k(dy,du)=\int_G q(y,u)\g(dy,du)
$$
for any $q\in C(G)$.

REMARK. Note that $\P(G)$ is weak$^*$ compact (as implied by Banach-Alaoglu theorem; see, e.g., Theorem 3.5.16 in \cite{Ash}). Therefore, it is also compact
in metric $\rho$ defined in (\ref{e-rho}). Also, any weak$^*$ closed subset of $\P(G)$ is compact in this metric.

Using the metric $\rho$, we can define the ``distance" $\rho(\g,\Gamma)$ between $\g\in \P(G)$ and $\Gamma\subset \P(G)$
and the Hausdorff metric $\rho_H(\Gamma_1,\Gamma_2)$ between $\Gamma_1\subset \P(G)$ and $\Gamma_2\subset \P(G)$ as follows:
$$
\rho(\g,\Gamma)\BYDEF\inf_{\g'\in \Gamma}\rho(\g,\g'),\quad
\rho_H(\Gamma_1,\Gamma_2)\BYDEF\max\{\sup_{\g\in \Gamma_1}\rho(\g,\Gamma_2),\sup_{\g\in \Gamma_2}\rho(\g,\Gamma_1)\}.
$$
Note that, although, by some abuse of terminology,  we refer to
$\rho_H(\cdot,\cdot)$ as  a metric on the set of subsets of
${\mathcal P} (G)$, it is, in fact, a semi-metric on this set
(since $\rho_H(\Gamma_1, \Gamma_2)=0$ implies  $\Gamma_1
= \Gamma_2$ if  $\Gamma_1$ and $\Gamma_2$ are closed, but the equality may not be true if at least one of these sets is not closed).

\section{IDLP presentations of optimal control problems: Vanishing discounting as a singular perturbation}\label{Sec-Time-discounting}
Define the sets $ W(\e,y_0)\subset \P(G)$ and $ W\subset \P(G)$ by the equations
\begin{equation}\label{w-disc}
\begin{aligned}
W(\e,y_0)\BYDEF&\bigg\{\g\in \P(G) |\, \int_G \big((1-\e)(E\big[\ph(f(y,u,s))\big]-\ph(y))+\e(\ph(y_0)-\ph(y))\big)\gm=0\ \ \forall \ph\in C(Y)\bigg\}
\end{aligned}
\end{equation}
and
\begin{equation}\label{M17}
\begin{aligned}
W\BYDEF\left\{\g\in \P(G) |\, \int_{G}\big(E\big[\ph(f(y,u,s))\big]-\ph(y)\big)
\gm=0\ \ \forall \ph\in C(Y)\right\},
\end{aligned}
\end{equation}
where $s $  stands for a generic random variable that has the same distribution as $s(t)$. Note that, as can be readily seen,  the sets $ W(\e,y_0)$ and $W$ are convex. Also they are weak$^*$ closed (due to the continuity of the functions under the  integrals (\ref{w-disc}) and (\ref{M17})
 implied by the assumptions made). Therefore, these sets are  compact in metric  $\rho$ defined in (\ref{e-rho}) (see Remark after the definition of $\rho$).
 For convenience,  $W(\e,y_0)$ and $W$ will be referred to as the sets of {\it discounted stationary} and {\it stationary} measures (respectively).

 Consider the following two optimization problems
\begin{equation}\label{D1}
\min_{\g\in W(\e,y_0)}  \int_{G} k(y,u)\g(dy,du)\BYDEF k^*(\e,y_0),
\end{equation}
\begin{equation}\label{M22}
\min_{\g\in W} \int_{G} k(y,u)\g(dy,du)\BYDEF k^*.
\end{equation}
Note that these problems belong  to the class of IDLP problems since both  the objective function and the constraints defining $W(\e,y_0) $ and $W$ are linear in $\gamma$.
The problem dual to  (\ref{D1}) and (\ref{M22}) can be written as follows
(see \cite{And-1}, \cite{And-2} and Section \ref{Sec-LP1} below):
\begin{equation}\label{e-dis-dual}
\sup_{\psi\in C(Y)}\min_{(y,u)\in G}\big\{k(y,u)+(1-\e)\big(E\big[\psi(f(y,u,s))\big]-\psi(y)\big)+\e\big(\psi(y_0)-\psi(y)\big)\big\}\BYDEF\mu^*( \e, y_0)
\end{equation}
and, respectively,
\begin{equation}\label{M21}
\sup_{\psi\in C(Y)}\min_{(y,u)\in G}\big\{k(y,u)+E\big[\psi(f(y,u,s))\big]-\psi(y)\big\}\BYDEF \mu^* .
\end{equation}

\begin{Proposition}\label{SD-1}
The following \lq\lq strong duality" equalities are valid:
\begin{equation}\label{e-SE-1}
k^*(\e,y_0)= \mu^*( \e, y_0) \ \ \forall \e\in (0,1) ,
\end{equation}
and
\begin{equation}\label{e-SE-2}
k^*= \mu^*.
\end{equation}
\end{Proposition}
{\bf Proof.} The proof is given in Section \ref{Sec-LP1}. $ \ \Box$

Relationships between the occupational measures sets introduced in (\ref{union-1}) and the sets $W(\e,y_0)$, $W$ are established by the following theorem.
\begin{Theorem}\label{Th-discounting-1}
The following equalities are valid
\begin{equation}\label{convergence-to-W-dis}
\bar{\rm co} \T_{\e}(y_0)=W(\e,y_0) \ \ \ \forall \ \e\in (0,1),
\end{equation}
\begin{equation}\label{convergence-to-W-dis-1}
\lim_{\e\rightarrow 0}\rho_H (\bar{\rm co} \T_{\e},W)=0, \ \ \ {\rm where}\ \ \ \T_{\e}\BYDEF\bigcup_{y_0\in Y}\{\T_{\e}(y_0)\},
\end{equation}
\begin{equation}\label{convergence-to-W-dis-2}
\lim_{T\rightarrow \infty}\rho_H (\bar{\rm co} \G_{T},W)=0, \ \ \ {\rm where}\ \ \ \G_T\BYDEF\bigcup_{y_0\in Y}\{\G_T(y_0)\},
\end{equation}
where ${\rm \bar{co}} $ stands for the closed convex hulls of the corresponding sets.
\end{Theorem}
{\bf Proof.} The proof of the theorem is given in Section \ref{Proofs-Secondary}. Note that it is very
 similar to the proofs of the corresponding results in the deterministic setting; see Corollary 2 and Theorem 5.4 in \cite{GPS-17}.  $\ \Box$

\begin{Corollary}\label{e-three-equalities}
The optimal values of the problems (\ref{A112-1}), (\ref{A112-2}) are related to the optimal values of the IDLP problems (\ref{D1}), (\ref{M22}) by the equalities:
\begin{equation}\label{convergence-to-W-dis-3}
h_{\e}(y_0) = k^*(\e,y_0) \ \ \ \forall \ \e\in (0,1),
\end{equation}
\begin{equation}\label{convergence-to-W-dis-4}
\lim_{\e\rightarrow 0}\min_{y_0\in Y}h_{\e}(y_0) = k^*,
\end{equation}
\begin{equation}\label{convergence-to-W-dis-5}
\lim_{T\rightarrow \infty}\min_{y_0\in Y}v_{T}(y_0) = k^*.
\end{equation}
\end{Corollary}

{\bf Proof.} The fact that (\ref{convergence-to-W-dis-3}) is implied by (\ref{convergence-to-W-dis}) follows from  the validity of the equalities
\begin{equation}\label{auxiliary-1}
h_{\e}(y_0)=\min_{\g \in  \T_\e(y_0)}\int_G k(y,u)\g(dy,du) = \min_{\g \in \bar{\rm co} \T_\e(y_0)}\int_G k(y,u)\g(dy,du)
\end{equation}
(see (\ref{A3}) and (\ref{D1})).
The fact that  (\ref{convergence-to-W-dis-4}) is implied by  (\ref{convergence-to-W-dis-1}) follows from that
\begin{equation}\label{auxiliary-2}
\min_{y_0\in Y}h_{\e}(y_0) =\min_{\g \in  \T_\e}\int_G k(y,u)\g(dy,du) = \min_{\g \in \bar{\rm co} \T_\e}\int_G k(y,u)\g(dy,du),
\end{equation}
and the fact that (\ref{convergence-to-W-dis-5}) is implied by  (\ref{convergence-to-W-dis-2}) follows from   that
\begin{equation}\label{auxiliary-3}
\min_{y_0\in Y}v_{T}(y_0) =\min_{\g \in  \G_{T}}\int_G k(y,u)\g(dy,du) = \min_{\g \in \bar{\rm co} \G_{T}}\int_G k(y,u)\g(dy,du)
\end{equation}
(see (\ref{A3-T}) and (\ref{M22})).
$ \ \Box$

\medskip

REMARK. Note that the fact that   the relations similar to (\ref{convergence-to-W-dis-3}),   (\ref{convergence-to-W-dis-4}),  (\ref{convergence-to-W-dis-5}) are valid is,
to a certain extent, a common knowledge in the MDP community,  see, e.g., relevant results in \cite{Jean-H-2}, \cite{Jean-H-3} and
  Theorem 3.6 in   \cite{Vega}.

\medskip

Corollary \ref{e-three-equalities} can be strengthened. Namely, the following proposition is valid.

\begin{Proposition}\label{Propo-separation-consec}
The relations (\ref{convergence-to-W-dis}), (\ref{convergence-to-W-dis-1}) and
(\ref{convergence-to-W-dis-2}) are valid if and only if (\ref{convergence-to-W-dis-3}), (\ref{convergence-to-W-dis-4}) and (\ref{convergence-to-W-dis-5}) (respectively) are valid for any continuous $k(y,u)$.
\end{Proposition}

{\bf Proof.} The proof is given is Section \ref{Sec-Appendix}. Note here only that it is  based on the separation theorem  (see, e.g., \cite{Rudin}, p. 59).  $\ \Box$

As can be seen, the set $W$  can be obtained by formally  taking $\e = 0 $ in (\ref{w-disc}). Thus, the IDLP problem (\ref{D1}) can be considered to be \lq\lq perturbed" with respect to the \lq\lq reduced" IDLP problem  (\ref{M22}). Let us  verify that
\begin{equation}\label{e-SingP-1-1}
\ \limsup_{\e\rightarrow 0}W(\e,y_0)\subset W .
\end{equation}
In fact, let $\g_i\in W(\e_i,y_0), \ i=1,2,... $, $\e_i\to 0$ and  $\g_i\to\g$. Since the functions  under the integrals in the expression (\ref{w-disc}) for $W(\e,y_0)$
 are continuous and, therefore, bounded (due to the assumptions made), we can substitute $\e_i $ for $\e$ and pass to the limit with $i\to\infty$ in this expression. We will, thus, verify that  $\g\in W$, which  establish the validity of the inclusion (\ref{e-SingP-1-1}). The latter implies that
\begin{equation}\label{e-SingP-1}
 \liminf_{\e\rightarrow 0}k^* (\e,y_0)\geq k^* .
\end{equation}
(Note that     (\ref{e-SingP-1})  also follows directly from the relations (\ref{convergence-to-W-dis-3}) and (\ref{convergence-to-W-dis-4}).)
Inequality (\ref{e-SingP-1}) can be strict. That is, the optimal value of the IDLP problem (\ref{D1}) can be discontinuous at $\e=0$. Problems characterized by  such a discontinuity  are called   {\it singularly perturbed} (see  \cite{Kostia-0}, \cite{Kostia-1}, \cite{Kostia-2}, \cite{GR-1} and \cite{Per-Gai}). In line with Corollary \ref{e-three-equalities}, the strict inequality in (\ref{e-SingP-1}) may occur if the  Ces\`aro and Abel limits of the optimal values are dependent on initial conditions (the so called \lq\lq non-ergodic" case). In fact, as demonstrated by the  example below, these limits may exist, be equal to each other, and be strictly greater than $k^*$.

{\bf Example 1.} Let the dynamics be one-dimensional and be described by the  equation (compare with (\ref{A1}))
\begin{equation}\label{example-1-0}
y(t+1)= y(t)u(t)s(t) \ \ \ \forall \ t=0,1,...\ ,
\end{equation}
where
$
\ Y=[-1,1]
$ and  $U(y)=\{-1, 1\}$ (that is, the control can be either equal to $1$ or to $-1$). Assume that $s(t)$ takes only two values: $s(t)=1$ with probability $\frac{3}{4}$ and $s(t)=-1$ with probability $\frac{1}{4} $. Consider problem (\ref{A112-1}) with
$\
k(y,u)=y
$. It can be readily understood, that, in this example,  the  plan  $\pi^* = \{\pi_t^*(y), \ t=0,1,...\} $, where, for any $ t=0,1,.... $,
\begin{equation}\label{e-ex-1-1}
\pi_t^*(y)= +1\ \ \ {\rm for}\ \  y\in [-1,0] \ \ \ {\rm and }\ \ \ \pi_t^*(y)= -1 \ \ \ {\rm for}\ \  y\in (0,1],
 \end{equation}
is optimal in both problem (\ref{A112-1}) and problem (\ref{A112-2}) (as this is the plan that maximizes the probability for the state variable to be negative). The optimal values of problems (\ref{A112-1}) and  (\ref{A112-2}) can be evaluated to be as follows
\begin{equation}\label{e-ex-1-2}
v_T(y_0)=-\frac{1}{2}|y_0| +\frac{1}{T}\left(y_0 + \frac{1}{2}|y_0|\right)  \ \ \ \forall\ y_0\in Y,
 \end{equation}
 \begin{equation}\label{e-ex-1-3}
 h_{\e}(y_0) = -\frac{1}{2}|y_0| +\e \left(y_0 + \frac{1}{2}|y_0|\right)  \ \ \ \forall\ y_0\in Y.
 \end{equation}
(By a direct substitution, one can  verify that $v_T(y_0)$  and $  V_{\e}(y_0)\BYDEF \e^{-1} h_{\e}(y_0) $, defined in accordance with (\ref{e-ex-1-2}) and (\ref{e-ex-1-3}), satisfy the dynamic programming equations (\ref{e-App-1}) and (\ref{e-App-2}), respectively.) From (\ref{e-ex-1-2}) and (\ref{e-ex-1-3}) it follows that
\begin{equation}\label{e-ex-1-3-1}
\lim_{\e\rightarrow 0} k^* (\e,y_0) = \lim_{\e\rightarrow 0} h_{\e}(y_0) = \lim_{T\rightarrow \infty} v_T(y_0) = -\frac{1}{2}|y_0|
\end{equation}
and
$$
k^*=\lim_{\e\rightarrow 0}\min_{y_0'\in Y}h_{\e}(y_0')=\lim_{T\rightarrow \infty}\min_{y_0'\in Y}v_{T}(y_0')=-\frac{1}{2}.
$$
Thus, inequality (\ref{e-SingP-1}) is strict in this case if $|y_0|\neq 1 $.

Since the reduced IDLP problem (\ref{M22}) does not involve the dependence on the initial condition $y_0$, it is not surprising that its optimal value is not equal to (or may not even provide a good estimate for)
the  Ces\`aro and Abel limits of the optimal values in case the latter depend on $y_0$. In the next section, we will introduce an IDLP problem that allows one to capture such a dependence.

Note that it is well known (and also  readily verifiable) that the set of test functions used in the definitions of the sets $W(\e,y_0) $ and $W$ can be expanded, with the
latter  being representable in the form:
\begin{equation}\label{w-disc-2}
\begin{aligned}
W(\e,y_0)=&\bigg\{\g\in \P(G)\bigg|\, \int_G \big((1-\e)\big(E\big[\ph(f(y,u,s))\big]-\ph(y)\big)
 +\e(\ph(y_0)-\ph(y))\big)\gm =0\ \\&\ \ \ \ \ \ \ \ \ \ \ \ \quad\quad\quad\quad\quad\quad\quad\quad\quad\quad\quad \  \forall\ \ph\in\mathcal{B}(Y)\bigg\},
\end{aligned}
\end{equation}
\begin{equation}\label{w-disc-2-1}
W=\bigg\{\g\in \P(G)\bigg|\, \int_G \big(E\big[\ph(f(y,u,s))\big]-\ph(y)\big)\gm =0 \ \ \forall\ \ph\in \mathcal{B}(Y)\bigg\},
\end{equation}
where $\mathcal{B}(Y) $ stands for the space of  bounded Borel functions on $Y$.
In addition (and it is readily verifiable too)
the sets  $W(\e,y_0) $ and $W$  allow  the following representations:
\begin{equation}\label{w-disc-1}
\begin{aligned}
W(\e,y_0)=&\bigg\{\g\in \P(G)\bigg|\ \gamma_1(Q)=  (1-\e)  \int_{G}P(Q|y,u)\gamma(dy,du) +\e 1_{Q}(y_0)\ \ \forall {\rm  \ Borel} \ Q\subset Y \bigg\},
\end{aligned}
\end{equation}
\begin{equation}\label{M17-1}
\begin{aligned}
W=\bigg\{\g\in \P(G)\bigg|\, \gamma_1(Q) = \int_{G}P(Q|y,u)\gamma(dy,du)
\ \ \forall {\rm  \ Borel} \ Q\subset Y\bigg\},
\end{aligned}
\end{equation}
where
 $\gamma_1$ is the marginal of $\gamma$, that is,
\begin{equation}\label{e-continuous-Borel-2}
\gamma_1(Q)= \int_{G}1_Q(y) \gamma (dy,du) ,
\end{equation}
and  $P(dy|y,u)  $ is the transition law associated with system (\ref{A1}), that is,
\begin{equation}\label{e-continuous-Borel-1}
P(Q|y,u) = E[1_Q(f(y,u,s))] \ \  \ \forall \ (y,u)\in G.
\end{equation}

REMARK.
The validity of the representations (\ref{w-disc-1}) and (\ref{M17-1}) make the results established by Theorem \ref{Th-discounting-1} and   Corollary \ref{e-three-equalities} consistent
with  well known results in Markov control processes theory;
see \cite{Altman1999,Borkar1991,Borkar2002,Jean-H-2,Jean-H-3,Gonzalez-Hernandes,Piun-1997,PZ-2011} and references therein. Many of the latter are obtained under assumptions that are lighter than the assumptions we are using in this paper.  Note that some of our assumptions can be relaxed too. For example, the assumption about compactness of the state space $Y$ can be replaced by the assumption about the tightness of the set of occupational measures that make the results of Theorem \ref{Th-discounting-1} valid.
However, to make the presentation more expository, we stick to using simpler (albeit more restrictive) assumptions.

\section{Augmented IDLP problem: Upper and lower bounds for the Ces\`aro and Abel limits of the  optimal values}\label{Sec-Augmented-estimates}

If an LP problem is singular perturbed (SP), then the (independent of $\e$) LP problem that defines the limit of the optimal value  as $\e\rightarrow 0 $
can be constructed by augmenting the reduced problem with some  additional constraints and decision variables. Such an augmentation was established
to be effective for the SP  LP problems in finite dimensions (see \cite{Kostia-1}, \cite{Kostia-2}, \cite{Per-Gai} as well as \cite{HK-1}, \cite{HK-2}). For the SP IDLP problem
(\ref{D1}), by formally following the construction of \cite{Kostia-1} and \cite{Kostia-2}, one can arrive at
the IDLP problem
\begin{equation}\label{BB1}
\inf_{(\g,\xi)\in \O(y_0)} \int_G k(y,u)\gm \BYDEF k^*(y_0),
\end{equation}
where
\begin{equation}\label{eq-Omega}
\begin{aligned}
&\O(y_0)\BYDEF\{(\g,\xi)\in \P(G)\times \Mp(G)|\, \g\in W,\,\\
&\int_{G}(\ph(y_0)-\ph(y))\gm+\int_G\big(E\big[\ph(f(y,u,s))\big]-\ph(y)\big)\xm=0
\quad \hbox{for all } \ph\in C(Y)\},
\end{aligned}
\end{equation}
where $\mathcal{M}_+ (G) $ stands for the space of nonnegative finite measures defined on Borel subsets of $G$.
This problem is obtained by augmenting  the reduced problem (\ref{M22}) with additional constraints involving a new \lq\lq decision variable" $\xi $. Note that
the presence of the  additional constraints imply that $\ \O(y_0)\subset W \ \ \forall\ y_0\in Y$,  and, therefore, the optimal values of problems (\ref{M22}) and
(\ref{BB1})  are related by the inequality
$$
k^*\leq k^*(y_0) \ \ \ \forall \ y_0\in Y.
$$

 The problem dual to the augmented IDLP problem (\ref{BB1}) can be written in the form (see Section \ref{Sec-LP2})
 \begin{equation}\label{BB8}
\sup_{(\mu,\psi,\eta)\in \D(y_0)} \mu \BYDEF d^*(y_0),
\end{equation}
where $\D(y_0)$ is the set of triplets $(\mu,\psi(\cdot),\eta(\cdot))\in \reals\times C(Y)\times C(Y)$ that for all $(y,u)\in G$ satisfy the inequalities
\begin{equation}\label{BB7}
\begin{aligned}
&k(y,u)+(\psi(y_0)-\psi(y))+E[\eta(f(y,u,s))]-\eta(y)-\mu\ge 0,\\
&E[\psi(f(y,u,s))]-\psi(y)\ge 0.
\end{aligned}
\end{equation}
Note that the optimal value of problem (\ref{BB8})  can be equivalently represented as
\begin{equation}\label{CC9}
d^*(y_0)=\sup_{(\psi,\eta)\in C(Y)\times C(Y)}\min_{(y,u)\in G}\{k(y,u)+(\psi(y_0)-\psi(y))+E[\eta(f(y,u,s))]-\eta(y)\},
\end{equation}
where $\psi$  satisfies the second inequality in \eqref{BB7}.

The following proposition establishes the validity of the {\it weak duality} inequality.


\begin{Proposition}\label{Prop-weak-duality-Omega}
The optimal values of (\ref{BB1}) and (\ref{BB8}) are related by the inequality
\begin{equation}\label{eq-weak-duality-Omega-Borel}
d^*(y_0)\leq  k^*(y_0).
\end{equation}
\end{Proposition}
{\bf Proof.}
Take any $(\g,\xi)\in \O(y_0)$ and $(\mu,\psi,\eta)\in \D( y_0)$. Integrating the first inequality in \eqref{BB7} with respect to $\g$ and taking into account that $\g\in W$, we conclude that
$$
\int_G k(y,u)\gm+\int_G(\psi(y_0)-\psi(y))\gm\ge \mu.
$$
Since $(\g,\xi)\in \O(y_0)$, from the second inequality in \eqref{BB7} it follows that
$$
\int_G(\psi(y_0)-\psi(y))\gamma(dy,du)=-\int_G \big(E[\psi(f(y,u,s))]-\psi(y)\big)\xi(dy,du)\le 0.
$$
Therefore,
$$
\int_G k(y,u)\gm  \ge \mu .
$$
Taking first $inf $ over all $(\g,\xi)\in \O(y_0)$ in the left-hand-side and then $sup$ over all  $(\mu,\psi,\eta)\in \D( y_0)$ in the right-hand-side, one establishes the validity of (\ref{eq-weak-duality-Omega-Borel}).

{\bf Example 1 (continuation).} The augmented IDLP problem (\ref{BB1}) takes the form
\begin{equation}\label{BB1-ex-1}
\inf_{(\g,\xi)\in \O(y_0)} \int_G y\gm = k^*(y_0),
\end{equation}
where  $\O(y_0)$ is the set of pairs $(\g,\xi)\in \P(G)\times \Mp(G)$ that satisfy the equations
\begin{equation}\label{BB1-ex-2}
\int_{G}\left(\frac{3}{4}\ph(yu) + \frac{1}{4}\ph(-yu) -\ph(y)\right)\gm = 0 \ \  \ \forall \ \ph\in C([-1,1]),
\end{equation}
\begin{equation}\label{BB1-ex-3}
\int_{G}(\ph(y_0)-\ph(y))\gm = -\int_{G}\left(\frac{3}{4}\ph(yu) + \frac{1}{4}\ph(-yu) -\ph(y)\right)\xm \ \ \ \forall \ \ph\in C([-1,1]),
\end{equation}
and where $G=Y\times U = [-1,1]\times \{-1,1\} $ in this case. The corresponding dual problem (see (\ref{CC9})) is
\begin{equation}\label{BB1-ex-4}
\sup_{(\psi,\eta)\in C([-1,1])\times C([-1,1])}\min_{(y,u)\in G}\left\{y+(\psi(y_0)-\psi(y))+\frac{3}{4}\eta(yu))+\frac{1}{4}\eta(-yu) -\eta(y)\right\}=d^*(y_0),
\end{equation}
where the  $\psi$ functions are assumed to satisfy the inequality
\begin{equation}\label{BB1-ex-5}
\frac{3}{4}\psi(yu)+\frac{1}{4}\psi(-yu) -\psi(y)\geq 0 \ \ \ \forall \ y\in [-1,1], \ \ \forall \ u\in \{-1, 1\}.
\end{equation}
If a function $\ph$ is even, then $\ \frac{3}{4}\ph(yu) + \frac{1}{4}\ph(-yu) -\ph(y)\equiv 0$ (since $u$ is either equal to $1$ or to $-1$). Therefore, (\ref{BB1-ex-2}) is satisfied for
all $\g\in \P(G) $,  while (\ref{BB1-ex-3}) is converted to $\ \int_{G}(\ph(y_0)-\ph(y))\gm = 0$ in this case. The latter equality implies that
$\ \int_{G} |y_0|^l\gm = \int_{G} |y|^l\gm$ for any $ l=1,2,...$, which, in turn, implies that
$$
\g(Y_{y_0}) =1, \ \ \ {\rm where}\ \ \ Y_{y_0}\BYDEF \{y: |y|=|y_0|\}.
$$
Thus, the constraints  (\ref{BB1-ex-3}) ensure that the  occupational measures $\g$ generated by the state-control trajectories satisfy the property $\g(Y\setminus Y_{y_0}) = 0 $. This is consistent with the system's dynamics (see (\ref{example-1-0})), according to which the  only states attended by the state trajectories are
$y_0$ and $-y_0$.

Let
\begin{equation*}\label{BB1-ex-6}
\bar{\g}(dy,du)\BYDEF \left(\frac{3}{4}\delta_{-|y_0|}(dy) + \frac{1}{4}\delta_{|y_0|}(dy)    \right)\delta_{\kappa(y)}(du), \ \ \ \ \  \bar{\xi}(dy,du) \BYDEF \delta_{y_0}(dy)\delta_{\kappa(y)}(du),
\end{equation*}
where $\delta_{a} $ stands for the Dirac measure concentrated at $a$, and where $\kappa(y) $ is equal to $1$ for $y\in [-1,0] $ and equal to $-1$ for
$y\in (0,1] $ (that is, for an arbitrary function $q(u)$ on $U$, $\ \int_{U}q(u)\delta_{\kappa(y)}(du) = q(1)\ \forall \ y\in [-1,0] $  and  $\ \int_{U}q(u)\delta_{\kappa(y)}(du) = q(-1) \ \forall \ y\in (0,1] $).

Via a direct substitution into (\ref{BB1-ex-2}) and (\ref{BB1-ex-3}), it can be verified that $(\bar{\g}, \bar{\xi})\in \O(y_0)$ (note that it is sufficient to verify
the validity of (\ref{BB1-ex-2}) and (\ref{BB1-ex-3}) only for the even and odd test functions $\ph(\cdot)$). Therefore,
\begin{equation}\label{BB1-ex-7}
k^*(y_0)\leq \int_G y\bar{\g}(dy,du)= -\frac{1}{2}|y_0|.
\end{equation}
On the other hand, it can also be verified that the pair of functions $(\bar{\psi}(y), \bar{\eta}(y)) $,
\begin{equation}\label{BB1-ex-8}
\bar{\psi}(y)\BYDEF -\frac{1}{2}|y|,\ \ \ \ \ \ \ \bar{\eta}(y)\BYDEF \left(y + \frac{1}{2}|y|\right),
\end{equation}
satisfy the relationships
\begin{equation*}\label{BB1-ex-9}
\min_{y\in [-1,1]}\min_{u\in \{-1,1\}} \left\{y+(\bar{\psi}(y_0)-\bar{\psi}(y))+\frac{3}{4}\bar{\eta}(yu))+\frac{1}{4}\bar{\eta}(-yu) -\bar{\eta}(y)\right\}= -\frac{1}{2}|y_0|,
\end{equation*}
$$
\frac{3}{4}\bar{\psi}(yu))+\frac{1}{4}\bar{\psi}(-yu) -\bar{\psi}(y)= 0 \ \ \ \forall \ y\in [-1,1], \ \ \forall \ u\in \{-1, 1\}.
$$
Therefore (compare the latter with (\ref{BB1-ex-4}) and (\ref{BB1-ex-5})),
$$
-\frac{1}{2}|y_0|\leq d^*(y_0).
$$
This inequality, along with (\ref{eq-weak-duality-Omega-Borel}) and (\ref{BB1-ex-7}), allows one to conclude that the optimal value of the IDLP problem (\ref{BB1-ex-1})
and the optimal value of the dual problem (\ref{BB1-ex-4}) are equal (that is, the strong duality equality is valid) and also that $(\bar{\g}, \bar{\xi})$ is an optimal solution of the former and $(\bar{\psi}(y), \bar{\eta}(y)) $ is an optimal solution of the latter. Note that the common optimal value of  problems (\ref{BB1-ex-1}) and (\ref{BB1-ex-4}) coincides with the Ces\`aro and  Abel  limits
(\ref{e-ex-1-3-1}).

Theorem \ref{Th-lower-upper-limit-main} and Corollary \ref{Cor-Strong-dual} stated below establish that the optimal values of the augmented IDLP problem and its dual  give  upper and (respectively) lower bounds
for the Ces\`aro and  Abel  limits of the optimal values, the existence and the equality of the latter being ensured if there is no duality gap.
To state these results, consider the following problem
\begin{equation}\label{BB8-B}
\sup_{(\mu,\psi,\eta)\in \hat{\D}(y_0)} \mu \BYDEF \hat d^*(y_0),
\end{equation}
where (in contrast to (\ref{BB8})) the {\it sup} is over  the set $\hat{\D}(y_0)$ consisting  of the triplets $(\mu,\psi(\cdot),\eta(\cdot))\in \reals\times \mathcal{B}(Y)\times \mathcal{B}(Y)$ that for all $(y,u)\in G$ satisfy the inequalities
\begin{equation}\label{BB7-B}
\begin{aligned}
&k(y,u)+(\psi(y_0)-\psi(y))+E[\eta(f(y,u,s))]-\eta(y)-\mu\ge 0,\\
&E[\psi(f(y,u,s))]-\psi(y)\ge 0.
\end{aligned}
\end{equation}
Problem (\ref{BB8-B}) is a \lq\lq relaxed" version of the dual problem (\ref{BB8}) that corresponds to the presentation of the feasible set $\O(y_0)$ in the form (\ref{eq-Omega-Borel}). Note that the optimal value of problem (\ref{BB8-B})  can be equivalently represented as
\begin{equation}\label{CC9-B}
\hat d^*(y_0)=\sup_{(\psi,\eta)\in \mathcal{B}(Y)\times \mathcal{B}(Y)}\inf_{(y,u)\in G}\{k(y,u)+(\psi(y_0)-\psi(y))+E[\eta(f(y,u,s))]-\eta(y)\},
\end{equation}
where $\psi$  satisfies the second inequality in \eqref{BB7-B} (compare with (\ref{CC9})). Note also that
\begin{equation}\label{CC9-B-1}
d^*(y_0)\leq \hat d^*(y_0)\leq k^*(y_0),
\end{equation}
the second inequality in (\ref{CC9-B-1})  being established similarly to Proposition \ref{Prop-weak-duality-Omega} (see Remark at the end of this section).

\begin{Theorem}\label{Th-lower-upper-limit-main}
The lower and upper  Ces\`aro/Abel  limits of the optimal value functions in problems \eqref{A112-1} and \eqref{A112-2} satisfy the inequalities:
\begin{equation}\label{BB3-relaxed-genearl}
\begin{aligned}
& \hat d^*(y_0)\le \liminf_{T\to \infty} v_T(y_0)\le    \limsup_{T\to \infty} v_T(y_0)\le k^*(y_0)\ \ \ \forall \ y_0\in Y, \\
& \hat d^*(y_0)\le \liminf_{\e\rightarrow 0} h_{\e}(y_0)\le    \limsup_{\e\rightarrow 0} h_{\e}(y_0)\le k^*(y_0) \ \ \ \forall \ y_0\in Y,\\
\end{aligned}
\end{equation}
where $k^*(y_0)$ is the optimal value of the augmented IDLP problem (\ref{BB1}) and $\hat d^*(y_0)$ is the optimal value of its dual (\ref{BB8-B}).
\end{Theorem}
{\bf Proof.} The proof  of the theorem is given in Section \ref{Sec-Proof-Ab-Ch-estimates}. $ \ \Box$

\medskip

REMARK. Theorem \ref{Th-lower-upper-limit-main} is one of the main results of the paper. The estimates from below were known (see \cite{Gonzalez-Hernandes}) and, in fact, they are relatively easy to be verified (see the proof
of Proposition \ref{Prop-est-from-below}). The proof of the estimates from above is much more involved. It is based on the use of some dynamic programming and LP related results. Note that the proof we present in this paper is similar to the proof of Theorem 3.1 in \cite{BGS-2019}, where the estimates from above were obtained in a purely deterministic setting. However, in contrast to the aforementioned theorem, we do not assume that the \lq\lq before limit" optimal value functions are continuous, this requiring a more fine analysis for establishing the desired results.

\begin{Corollary}\label{Cor-Strong-dual}
Let, for a given $y_0\in Y $, the strong duality equality be valid:
\begin{equation}\label{eq-strong-duality-1}
k^*(y_0)= \hat d^*(y_0).
\end{equation}
 Then the  Ces\`aro and  Abel limits of the optimal values exist and are equal:
\begin{equation}\label{eq-strong-duality-2}
\lim_{T\rightarrow\infty} v_{T}(y_0) =\lim_{\e\rightarrow 0} h_{\e}(y_0) = k^*(y_0)= \hat d^*(y_0).
\end{equation}
\end{Corollary}
Note that, if
\begin{equation}\label{eq-strong-duality-1-1}
k^*(y_0)=  d^*(y_0)
\end{equation}
(as in Example 1), then, by (\ref{CC9-B-1}), $ d^*(y_0) = \hat d^*(y_0) = k^*(y_0)$. That is, (\ref{eq-strong-duality-1}) is valid, with (\ref{eq-strong-duality-2}) taking the form
$$
\lim_{T\rightarrow\infty} v_{T}(y_0) =\lim_{\e\rightarrow 0} h_{\e}(y_0) = k^*(y_0)=  d^*(y_0).
$$
REMARK. Note that the set $\O(y_0)$ allows also the  representation
\begin{equation}\label{eq-Omega-Borel}
\begin{aligned}
&\O(y_0)=\{(\g,\xi)\in \P(G)\times \Mp(G)|\, \g\in W,\,\\
&\int_{G}(\ph(y_0)-\ph(y))\gm+\int_G\big(E\big[\ph(f(y,u,s))\big]-\ph(y)\big)\xm=0
\quad \hbox{for all } \ph\in \mathcal{B}(Y)\}
\end{aligned}
\end{equation}
as well as the representation
\begin{equation}\label{eq-Omega-Borel-1}
\begin{aligned}
&\O(y_0)=\{(\g,\xi)\in \P(G)\times \Mp(G)|\, \ \gamma_1(Q) -\int_{G}P(Q|y,u)\gamma(dy,du) = 0
\ \ \forall {\rm  \ Borel} \ Q\subset Y ,\,\\
&\xi_1 (Q)-\int_G P(Q|y,u)\xm  + \gamma_1(Q) = 1_Q(y_0)
\ \ \forall {\rm  \ Borel} \ Q\subset Y \},
\end{aligned}
\end{equation}
where $\gamma_1$ and $\xi_1$ are marginals of $\gamma$ and $\xi$ (see comments at the end of Section \ref{Sec-Time-discounting}). Problem  (\ref{BB1}), in which the feasible set is defined in accordance with (\ref{eq-Omega-Borel-1}), has
been considered in  \cite{Gonzalez-Hernandes}, where it was shown that the corresponding dual problem has the form (\ref{BB8-B}) (thus, the second inequality in (\ref{CC9-B-1}) is just a version of the weak duality inequality).
Also in  \cite{Gonzalez-Hernandes},  sufficient conditions for   the equality (\ref{eq-strong-duality-1}) (the strong duality) to be valid  have been studied. Note
that the strong duality may not be true in the general case. An example, in which the Ces\`aro and Abel limits of the optimal values are not equal to each other, and, therefore, by Theorem \ref{Th-lower-upper-limit-main}, there is a duality gap, is given   in \cite{Ilya-arch}.

As stated in Corollary \ref{Cor-Strong-dual}, both the limit
\begin{equation}\label{chesaro-lim-1-1}
	\lim_{T\rightarrow\infty} v_{T}(y_0)\BYDEF v(y_0)
	\end{equation}
and the limit
\begin{equation}\label{abel-lim-3-1}
	\lim_{\e \rightarrow 0} h_\e(y_0)\BYDEF h(y_0)
	\end{equation}
exist and are equal to $\hat d^*(y_0)$ if (\ref{eq-strong-duality-1}) is true, or they exist and are equal to $ d^*(y_0)$ if (\ref{eq-strong-duality-1-1}) is valid. These  statements are complimented by the following theorem.

\begin{Theorem}\label{ThN1}
	{\rm (a)} Let limit  (\ref{chesaro-lim-1-1}) exist  for any $y_0\in Y$ (that is, $v_{T}(\cdot)$ converges to a function $v(\cdot)$ point-wisely  on $Y$), and let the limit function
  $v(\cdot)$ be continuous. Then
	\begin{equation}\label{CC10}
	v(y_0)= d^*(y_0)
	\end{equation}
for any $y_0\in Y $.
	{\rm (b)} Let limit  (\ref{abel-lim-3-1}) exist  for any $y_0\in Y$ (that is, $ h_\e(\cdot)$ converges to a function $h(\cdot)$ point-wisely on $Y$), and let the limit function
  $h(\cdot)$ be continuous. Then
	\begin{equation}\label{abel-lim-4-1}
	h(y_0)= d^*(y_0)
	\end{equation}
for any $y_0\in Y $.
\end{Theorem}
{\bf Proof.} The proof of the theorem is given in Section \ref{Sec-Proof-Ab-Ch-estimates}.  $\ \Box$

\medskip

REMARK.    Theorem \ref{ThN1} is a stochastic analog of Theorem 4.2 in \cite{BGS-2019}, where a similar statement was established in a deterministic setting.
In proving  the theorem, we use the dynamic programming and LP related results.

\begin{Corollary}\label{Cor-d-hat-d}
If the conditions of Theorem \ref{ThN1} {\rm (a)}   and/or the conditions of Theorem \ref{ThN1} {\rm (b)} are satisfied, then
\begin{equation}\label{lim-1-1-1}
d^*(y_0) = \hat d^*(y_0)\ \ \ \forall \ y_0\in Y.
\end{equation}
\end{Corollary}
{\bf Proof.} If the conditions of {\rm (a)} are satisfied, then, by  (\ref{BB3-relaxed-genearl}),
\begin{equation}\label{chesaro-lim-1-1-1}
	\hat d^*(y_0)\leq v(y_0) \ \ \ \forall \ y_0\in Y.
	\end{equation}
This, along with (\ref{CC9-B-1}) and (\ref{CC10}), implies (\ref{lim-1-1-1}). If the conditions of {\rm (b)} are satisfied, then
\begin{equation}\label{abel-lim-1-1-1}
	\hat d^*(y_0)\leq h(y_0) \ \ \ \forall \ y_0\in Y.
	\end{equation}
Similarly to (\ref{chesaro-lim-1-1-1}), this also  implies (\ref{lim-1-1-1}).  $\ \Box$

\medskip

Let us conclude this section with an example, in which  the limit  (\ref{chesaro-lim-1-1}) exists  for any $y_0\in Y$ (that is, $ v_T(\cdot)$ converges to  $  v(\cdot)$ pointwisely) but the limit
function $v(\cdot) $ is discontinuous.

\medskip

{\bf Example 2.} Let the dynamics be one-dimensional and be described by the  equation
\begin{equation}\label{e-ex-2-1}
y(t+1)= u(t)s(t) \ \ \ \forall \ t=0,1,...\ ,
\end{equation}
with
$
\ Y=[-1,1]$ and with $\ U(y)\subset \reals$ defined as follows
\begin{equation}\label{e-ex-2-2}
U(y)=[-1, y] \ \ {\rm for} \ \ y\in [-1,0); \ \ \ \ \ U(0)=[-1,1]; \ \ \ \ \ U(y)=[y,1] \ \ {\rm for} \ \ y\in (0,1].
\end{equation}
 Assume that $s(t)$ takes  two values: $s(t)=1$ with probability $\frac{1}{2}$ and $s(t)=\frac{1}{4}$ with probability $\frac{1}{2} $. Consider problem (\ref{A112-1}) with
$\
k(y,u)=y
$.
 It can be readily seen, that, in this example, the control plan  $\pi^* = \{\pi_t^*(y), \ t=0,1,...\} $, where, for any $ t=0,1,.... $,
\begin{equation}\label{e-ex-2-extra-1}
\pi_t^*(y)= -1 \ \ \ {\rm for}\ \  y\in [-1,0] \ \ \ {\rm and }\ \ \ \pi_t^*(y)= y \ \ \ {\rm for}\ \  y\in (0,1],
 \end{equation}
is optimal in  problem   (\ref{A112-1}) (and in problem (\ref{A112-2})), the corresponding optimal state-control trajectory being  as follows:
  \begin{equation}\label{e-extra-extra-3}
\begin{aligned}
&u(t)= - 1 \ \ \forall  \ t\geq 0 \ \  {\rm and}  \ \ y(t)= - s(t-1) \ \ \forall  \ t\geq 1  \ \  {\rm if}  \ \ y_0\in [-1,0];\\
&u(t)= y(t) \ \ \forall  \ t\geq 0 \ \  {\rm and}  \ \ y(t)=  s(t-1)\cdots s(0)y_0 \ \ \forall  \ t\geq 1  \ \  {\rm if}  \ \ y_0\in (0,1].
\end{aligned}
\end{equation}
 The optimal value function of the problem  (\ref{A112-1}) can be verified to be as follows:
\begin{equation}\label{eq-union-6-0}
v_T(y_0)= - \frac{5}{8} + \frac{1}{T}\Big(y_0 +\frac{5}{8} \Big) \ \  \ {\rm for}  \ \ \ y_0\in [-1,0] \ \ \ \ \  {\rm and } \ \ \ \ \
v_T(y_0)=\frac{1}{T}\Big( \frac{8}{3}y_0\Big)
\Big( 1 -\Big(\frac{5}{8} \Big)^T\Big)
\ \  {\rm for}  \ \ y_0\in (0,1],
\end{equation}
with
\begin{equation}\label{eq-union-6}
v(y_0)=\lim_{T\rightarrow\infty}v_T(y_0)= - \frac{5}{8} \ \  \ {\rm for}  \ \ \ y_0\in [-1,0] \ \ \ \ \  {\rm and } \ \ \ \ \
v(y_0)=\lim_{T\rightarrow\infty}v_T(y_0)= 0 \ \  {\rm for}  \ \ y_0\in (0,1].
\end{equation}
The augmented IDLP problem (\ref{BB1}) takes in this case the form
\begin{equation}\label{BB1-ex-2-5}
\inf_{(\g,\xi)\in \O(y_0)} \int_G y\gm = k^*(y_0),
\end{equation}
where  $\O(y_0)$ is the set of pairs $(\g,\xi)\in \P(G)\times \Mp(G)$ that satisfy the equations
\begin{equation}\label{BB1-ex-2-6}
\int_{G}\left(\frac{1}{2}\ph(u) + \frac{1}{2}\ph\big(\frac{u}{4}\big) -\ph(y)\right)\gm = 0 \ \  \ \forall \ \ph\in C([-1,1]),
\end{equation}
\begin{equation}\label{BB1-ex-2-7}
\int_{G}(\ph(y)-\ph(y_0))\gm = \int_{G}\left(\frac{1}{2}\ph(u) + \frac{1}{2}\ph\big(\frac{u}{4}\big) -\ph(y)\right)\xm \ \ \ \forall \ \ph\in C([-1,1]),
\end{equation}
and where $G = G_1\cup G_2\cup G_3$, with
$
\ G_1\BYDEF \{(y,u)\ | \ u\in [-1, y), \ y\in [-1,0)\},\ \   G_2\BYDEF \{(0,u)\ | \ u\in [-1, +1]\}$ and $\ G_3\BYDEF \{(y,u)\ | \ u\in [y,1], \ y\in (0,1]\}  $.
The optimal value of the relaxed dual problem is presentable in the form (see (\ref{CC9-B})):
\begin{equation}\label{CC9-B-ex-2}
\sup_{(\psi,\eta)\in \mathcal{B}(Y)\times \mathcal{B}(Y)}\inf_{(y,u)\in G}\{y+(\psi(y_0)-\psi(y))+\frac{1}{2}\eta(u ) + \frac{1}{2}\eta\big(\frac{u}{4}\big)-\eta(y)\}= \hat d^*(y_0),
\end{equation}
where the functions $\psi(\cdot)$ are assumed to satisfy the inequality
\begin{equation}\label{e-ex-2-4}
\frac{1}{2}\psi(u ) + \frac{1}{2}\psi\big(\frac{u}{4}\big)  - \bar\psi(y)\geq 0 \ \ \ \forall (y,u)\in G
\end{equation}
Define  $\bar \psi(\cdot)$ and  $ \bar\eta(\cdot) $ by the equations:
\begin{equation}\label{e-ex-2-3}
\begin{aligned}
&\bar \psi(y)\BYDEF - \frac{5}{8} \ \  \ {\rm for}  \ \ \ y\in [-1,0] \ \ \ \ \  {\rm and } \ \ \ \ \
\bar\psi(y)\BYDEF 0 \ \  {\rm for}  \ \ y\in (0,1];\\
&\bar \eta(y)\BYDEF y+ \frac{5}{8} \ \  \ {\rm for}  \ \ \ y\in [-1,0] \ \ \ \ \  {\rm and } \ \ \ \ \
\bar \eta(y)\BYDEF \frac{8}{3}y \ \  {\rm for}  \ \ y\in (0,1].
\end{aligned}
\end{equation}
Via the direct substitution  that the function $\bar \psi(\cdot)$ satisfies (\ref{e-ex-2-4})
and that
\begin{equation}\label{e-ex-2-5}
\begin{aligned}
& \inf_{(y,u)\in G}\{y+(\bar\psi(y_0)-\bar\psi(y))+\frac{1}{2}\bar \eta(u ) + \frac{1}{2}\bar \eta\big(\frac{u}{4}\big)-\bar\eta(y)\}= - \frac{5}{8} \ \  \ {\rm for}  \ \ \ y_0\in [-1,0],\\
&\inf_{(y,u)\in G}\{y+(\bar\psi(y_0)-\bar\psi(y))+\frac{1}{2}\bar\eta(u ) + \frac{1}{2}\bar\eta\big(\frac{u}{4}\big)-\bar\eta(y)\}= 0 \ \  \ {\rm for}  \ \ \ y_0\in (0,1].
\end{aligned}
\end{equation}
Thus (see (\ref{eq-union-6}) and (\ref{CC9-B-ex-2})), $\ \hat d^*(y_0)\geq v(y_0) $. Therefore, by (\ref{BB3-relaxed-genearl}),
\begin{equation}\label{e-ex-2-6}
\hat d^*(y_0)=v(y_0) \ \ \ \forall\ y_0\in [-1,1].
\end{equation}
For $ \ y_0\in [-1,0]$, let
\begin{equation}\label{e-ex-2-7}
\bar{\g}(dy,du)\BYDEF \left(\frac{1}{2}\delta_{(-1,-1)}(dy,du) + \frac{1}{2}\delta_{(-\frac{1}{4},-1)}(dy,du) \right),
\ \ \ \ \  \bar{\xi}(dy,du) \BYDEF \delta_{(y_0,-1)}(dy,du),
\end{equation}
where $\delta_{(\alpha,\beta)}\in \mathcal{P}(G) $ is the Dirac measure concentrated at a point $(\alpha,\beta)\in G $. It can be readily
verified  that $(\bar{\g}, \bar{\xi})\in \Omega (y_0) $ (with $\Omega (y_0)$ being defined by (\ref{BB1-ex-2-6}), (\ref{BB1-ex-2-7})) and that
$$
k^*(y_0)\leq \int_G y\bar{\g}(dy,du)= - \frac{5}{8}.
$$
Due to (\ref{BB3-relaxed-genearl}), (\ref{eq-union-6}) and  (\ref{e-ex-2-6}), the latter implies that
\begin{equation}\label{e-ex-2-8}
k^*(y_0)= v(y_0) = \hat d^*(y_0)  \ \ \ \forall\ y_0\in [-1,0].
\end{equation}
For $ \ y_0\in (0,1]$, take
\begin{equation*}\label{e-ex-2-9}
\bar{\g}(dy,du)\BYDEF  \delta_{(0,0)}(dy,du),
\ \ \ \ \  \bar{\xi}(dy,du) \BYDEF \delta_{(y_0,0)}(dy,du).
\end{equation*}
Again, it is easy to verify that
$(\bar{\g}, \bar{\xi})\in \Omega (y_0) $ and that
$$
k^*(y_0)\leq \int_G y\bar{\g}(dy,du)=0.
$$
Therefore, similarly to (\ref{e-ex-2-8}), we may conclude that the strong duality equality is valid:
\begin{equation}\label{e-ex-2-8-1}
k^*(y_0)= v(y_0) = \hat d^*(y_0)  \ \ \ \forall\ y_0\in (0,1].
\end{equation}
 In addition, we also may  conclude that $ (\bar\gamma, \bar \xi)$ defined in (\ref{e-ex-2-7}) is an optimal solution of
the IDLP problem (\ref{BB1-ex-2-5}) and that $(\bar \psi(\cdot), \bar\eta(\cdot))$ defined in (\ref{e-ex-2-3}) is an optimal  solution of the  dual problem (\ref{CC9-B-ex-2}).

\section{Another representation for the dual optimal value; Periodic regime generating controls}\label{Sec-Strong-Duality}
Define the set $\ D(y_0)\subset \mathcal{P}(G)  $ by the equation
\begin{equation}\label{eq-D-y-0}
\begin{aligned}
&D(y_0)\BYDEF\Big\{\g \in \M (G) \ | \  \exists \  \xi\in \Mp(G) \ \ {\rm such \ that}\ \ \int_{G}(\ph(y)-\ph(y_0))\gm \\
&=\int_G\big(E\big[\ph(f(y,u,s))\big]-\ph(y)\big)\xm
\quad \forall\ \ph\in C(Y)\Big\},
\end{aligned}
\end{equation}
where $\mathcal{M}(G) $ stands for the space of finite signed measures defined on Borel subsets of $G$. Obviously,
$$
W\cap D(y_0) = \{\gamma \ | \ (\gamma , \xi)\in \Omega (y_0)\},
$$
and the  problem (\ref{BB1}) can be rewritten as follows:
\begin{equation}\label{BB1-D-y-0-1}
 \inf \Big\{\int_G k(y,u)\gm \ | \ \gamma\in W\cap D(y_0)\Big\}=k^*(y_0).
\end{equation}
Along with the problem  (\ref{BB1-D-y-0-1}), let us consider the problem
\begin{equation}\label{BB1-D-y-0-2}
\min \Big\{\int_G k(y,u)\gm \ | \ \gamma\in W\cap D_1(y_0)\Big\} \BYDEF k^{**}(y_0),
\end{equation}
where 
\begin{equation}\label{BB1-D-y-0-2-1}
\begin{aligned}
&\ \ \ \ \ \ \ \ \ \ \ \ \ \ \ \ \ \ \ \ \  D_1(y_0)  = \{ \gamma\in \M(G)\  |\ \exists \xi_l\in \M_+(G), \ l=1,2,...,\ \  {\rm such\ that} \\
&\int_{G}(\ph(y)-\ph(y_0))\gm = \lim_{l\to\infty}\int_G\big(E\big[\ph(f(y,u,s))\big]-\ph(y)\big)\xi_l(dy,du)
\quad \forall\ \ph\in C(Y)\}.
\end{aligned}
\end{equation}
Note that, as can be readily verified, $D_1(y_0)$ is a convex and weak$^*$ closed set, and $cl (D(y_0))\subset D_1(y_0)$ (with $cl(\cdot)$ standing for the weak$^*$ closure).  
\begin{Proposition}\label{Lemma-conv-to cl-D-y-0}
The following inclusion is valid
\begin{equation}\label{BB1-D-y-0-2-2}
\limsup_{T\rightarrow\infty}\Gamma_T(y_0)\subset W\cap D_1(y_0).
\end{equation}
\end{Proposition}
{\bf Proof.}
Due to (\ref{convergence-to-W-dis-2}), $\limsup_{T\rightarrow\infty}\Gamma_T(y_0)\subset W$. Hence, we only need to prove that
\begin{equation}\label{BB1-D-y-0-2-3}
\limsup_{T\rightarrow\infty}\Gamma_T(y_0)\subset  D_1(y_0).
\end{equation}
Take an arbitrary $\gamma\in \limsup_{T\rightarrow\infty}\Gamma_T(y_0) $. By  definition of $\limsup$, it means that there exist a sequence   $T_l\to\infty $ and a sequence
$\gamma_l\in \Gamma_{T_l}(y_0)$ such that $\g_l\to \g   $. The fact that $\gamma_l\in \Gamma_{T_l}(y_0)$ implies that it is generated by some
 control plans $\pi_l$ with the corresponding state-control trajectory $(y^{\pi_l, y_0 }(\cdot), u^{\pi_l, y_0 }(\cdot))$. That is, $\ \gamma_l=\gamma^{\pi_l , y_0, T_l} $, and
 (see (\ref{G88}))
 $$
 \int_{G} q(y,u) \g_l(dy,du)={1\o T_l}E \left[\sum_{t=0}^{T_l-1} q(y^{\pi_l , y_0 }(t),u^{\pi_l , y_0 }(t))\right]
 $$
 for any continuous $q(y,u)$. Using this equality with $q(y,u)=\ph(y)-\ph(y_0) $, we obtain
 $$
\int_{G}(\ph (y)-\ph (y_0))\gamma_l(dy,du)= \frac{1}{T_l}\sum_{t=0}^{T_l-1}E[\ph(y^{\pi_l, y_0}(t))-\ph(y_0)]
$$
\vspace{-0.2cm}
$$
=\frac{1}{T_l}\sum_{t=0}^{T_l-1}\left(\sum_{t'=0}^{t-1}E[\ph(y^{\pi_l, y_0}(t'+1)) - \ph(y^{\pi_l, y_0}(t'))] \right)
$$
\vspace{-0.2cm}
$$
=\frac{1}{T_l}E\left[\sum_{t=T_0}^{T_l -1}\left(\sum_{t'=0}^{t-1}\big(E[\ph(f(y^{\pi_l, y_0}(t'), u^{\pi_l, y_0}(t'), s(t')))|y^{\pi_l, y_0}(t')] - \ph(y^{\pi_l, y_0}(t'))\big) \right)\right].
$$
Since
$$
E[\ph(f(y^{\pi_l, y_0}(t'), u^{\pi_l, y_0}(t'), s(t')))|y^{\pi_l, y_0}(t')] = \bar \ph(y^{\pi_l, y_0}(t'), u^{\pi_l, y_0}(t'))  \ \ a.s. ,
$$
with
\begin{equation}\label{e-phi-bar-1-1}
\bar \ph (y,u) \BYDEF E[\ph(f(y,u,s))] ,
\end{equation}
we obtain
\begin{equation}\label{e-phi-bar-1-2}
\int_{G}(\ph (y)-\ph (y_0))\gamma_l(dy,du) = \frac{1}{T_l}E\left[\sum_{t=0}^{T_l -1}\left(\sum_{t'=0}^{t-1}\big( \bar \ph(y^{\pi_l, y_0}(t'), u^{\pi_l, y_0}(t')) - \ph(y^{\pi_l, y_0}(t'))\big) \right)\right].
\end{equation}
 Define $\zeta_l\in C^*(G)$ by the equation
$$
\langle \zeta_l ,q\rangle ={1\o T_l}E\left[\sum_{t=0}^{T_l -1}\sum_{t'=0}^{t-1}q(y^{\pi_l, y_0}(t'), u^{\pi_l, y_0}(t'))\right]
 \ \ \ \ \forall \ q\in C(G).
$$
Due to Riesz representation theorem, there exists $\xi_l \in \Mp(G)$ such that
$$
\langle \zeta_l,q\rangle =\int_G q(y,u)\xi_l (dy,du)  \ \ \ \ \forall \  q(\cdot,\cdot)\in C(G).
$$
Therefore, (\ref{e-phi-bar-1-2}) can be rewritten as
$$
\int_{G}(\ph (y)-\ph (y_0))\gamma_l(dy,du)= \langle \zeta_l ,\bar \ph(y,u) - \ph  (y)\rangle = \int_{G}(\bar \ph(y,u)-\ph(y))\,\xi_l (dy,du).
$$
Passing to the limit with $l\to\infty$ in this equality and having in mind that $\gamma_l\to \g $, we obtain (see also (\ref{e-phi-bar-1-1}))
$$
\int_{G}(\ph (y)-\ph (y_0))\gamma (dy,du) = \lim_{l\to\infty}\int_G\big(E\big[\ph(f(y,u,s))\big]-\ph(y)\big)\xi_l(dy,du)
\quad \hbox{for all } \ph\in C(Y)
$$
Hence (see (\ref{BB1-D-y-0-2-1})), $\g\in D_1(y_0)$. Since $\g$ is an arbitrary element of  $ \limsup_{T\rightarrow\infty}\Gamma_T(y_0) $, the validity of (\ref{BB1-D-y-0-2-3}) is

\smallskip 

established. $ \ \Box$

From Proposition \ref{Lemma-conv-to cl-D-y-0} it follows that the optimal value of (\ref{BB1-D-y-0-2}) gives a lower bound for
$\ \liminf_{T\rightarrow\infty} v_{T}(y_0)$. This, however, does not provide an improvement for an estimate from below in (\ref{BB3-relaxed-genearl}) since, as established by the
proposition below, the optimal value of (\ref{BB1-D-y-0-2}) is equal to the optimal value of the dual problem  (\ref{BB8}).

\begin{Proposition}\label{Prop-cl-D-y-0-dual}
The optimal value  of the problem (\ref{BB1-D-y-0-2}) is equal to the optimal value of the problem (\ref{BB8}):
\begin{equation}\label{BB1-D-y-0-3}
k^{**}(y_0)= d^*(y_0).
\end{equation}
\end{Proposition}

{\bf Proof.} The proof is given in Section \ref{Sec-LP2}. $ \ \Box$

\medskip

REMARK. The proof of Proposition \ref{Prop-cl-D-y-0-dual} is based on the fact that the subvalue of an IDLP problem is equal to the optimal value of its dual   (see, e.g., Theorem 3 in \cite{And-1}).

\begin{Corollary}\label{Cor-conv-feasible}
Assume that, for any continuous $k(y,u)$, the point-wise limit (\ref{chesaro-lim-1-1}) exists and the limit optimal value function $v(\cdot)$ is continuous.
Then
\begin{equation}\label{convergence-to-W-dis-2-y-0}
\lim_{T\rightarrow \infty}\rho_H (\bar{\rm co} \G_{T}(y_0),W\cap D_1(y_0) )=0.
\end{equation}
\end{Corollary}
{\bf Proof.}
Note that
$$
v_T(y_0)= \inf_{\g\in  \G_{T}(y_0) }\int_{G}k(y,u)\g(dy,du)= \min_{\g\in \bar{\rm co} \G_{T}(y_0) }\int_{G}k(y,u)\g(dy,du).
$$
Therefore, from Theorem \ref{ThN1} (a) and from Proposition \ref{Prop-cl-D-y-0-dual} it follows that, for any continuous $k(y,u)$,
$$
\lim_{T\to\infty}\min_{\g\in \bar{\rm co} \G_{T}(y_0) }\int_{G}k(y,u)\g(dy,du) = \min_{\g\in W\cap D_1(y_0) }\int_{G}k(y,u)\g(dy,du).
$$
Continuing from this point in the same way as in the proof of  Proposition \ref{Propo-separation-consec} (that is, using Blaschke's selection theorem \cite{Klein}) and
the separation theorem \cite[p. 59]{Rudin})), one can establish the validity of (\ref{convergence-to-W-dis-2-y-0}).
$\ \Box$

\begin{Corollary}\label{Cor-the-only-conclusion}
The strong duality equality (\ref{eq-strong-duality-1-1}) is valid if
\begin{equation}\label{eq-D-y-0-7}
cl (W\cap  D(y_0)) = W\cap  D_1(y_0).
\end{equation}
\end{Corollary}
{\bf Proof.} The proof follows from  Proposition \ref{Prop-cl-D-y-0-dual} and from (\ref{BB1-D-y-0-1}), (\ref{BB1-D-y-0-2}).
$ \ \Box$

\medskip

REMARK. We do not provide sufficient conditions for the validity of (\ref{eq-D-y-0-7}) in the present paper leaving investigating of this matter for the future research.

\medskip

DEFINITION. For a given $y_0\in Y$, we will say that a control plan $\pi$ is periodic regime generating (PRG) if there exist  integer $T_0\geq 0$ and $\mathcal{T}>0$ such that,
for any continuous $q(y,u)$,
\begin{equation}\label{e-opt-4}
	E[q(y^{\pi, y_0}(t),u^{\pi, y_0}(t)) ]= E[q(y^{\pi, y_0}(t+\mathcal{T}),u^{\pi, y_0}(t+\mathcal{T})) ]\ \ \  \forall \ t\geq T_0.
	\end{equation}

Consider the problem
\begin{equation}\label{A112-1-per}
 \inf_{\pi\in \Pi_{per}(y_0)}\left\{\lim_{T\rightarrow\infty}{1\o T}E\left[\sum_{t=0}^{T-1} k(y^{\pi, y_0}(t),u^{\pi, y_0}(t))\right]\right\}\BYDEF V_{per}(y_0),
\end{equation}
where $\Pi_{per}(y_0)\subset \Pi$ stands for  the set of all PRG control plans.
Note that, as can be readily understood,
$$
 V_{per}(y_0)\geq \liminf_{T\to \infty} v_T(y_0)\ \ \ \forall \ y_0\in Y.
$$
The following statement is valid.
\begin{Proposition}\label{Prop-PGR}
If for some $y_0\in Y$,
\begin{equation}\label{A112-1-per-ineq-1}
  V_{per}(y_0) =  \liminf_{T\to \infty} v_T(y_0),
\end{equation}
then the limit (\ref{chesaro-lim-1-1}) exists, and
\begin{equation}\label{A112-1-per-ineq}
 v(y_0)=k^*(y_0).
\end{equation}
\end{Proposition}
{\bf Proof.} Due to the assumed validity of (\ref{A112-1-per-ineq-1}) and due to the upper bound in (\ref{BB3-relaxed-genearl}),
$$
 V_{per}(y_0)\leq k^*(y_0).
$$
Hence, if one shows that
\begin{equation}\label{e-V-per-k}
V_{per}(y_0)\geq k^*(y_0),
\end{equation}
it would imply
$$
k^*(y_0) =  V_{per}(y_0) =  \liminf_{T\to \infty} v_T(y_0) \leq \limsup_{T\to \infty} v_T(y_0)\leq k^*(y_0),
$$
which, in turn, would imply
the existence of the limit (\ref{chesaro-lim-1-1}) and the validity of the equality (\ref{A112-1-per-ineq}). That is, the proposition will be proved if
we show that (\ref{e-V-per-k}) is true. Thus, we will be aiming at proving (\ref{e-V-per-k}).

Note  that, as mentioned above,  the IDLP problem (\ref{BB1}) can be  rewritten in the form (\ref{BB1-D-y-0-1}).
Note also that problem (\ref{A112-1-per}) can be equivalently rewritten in the form
\begin{equation}\label{A3-T-per}
\min_{\g\in \G_{per}(y_0)} \int_{G} k(y,u)\g(dy,du)=V_{per}(y_0),
\end{equation}
where $\G_{per}(y_0) $ is the set of occupational measures generated by the PRG control plans: $\gamma\in \G_{per}(y_0)  $ if and only if
\begin{equation}\label{e-proof-ineq-01}
\g(Q)={1\o \mathcal{T}} E\left[\sum_{t=T_0}^{T_0+\mathcal{T}-1} 1_Q(y^{\pi , y_0 }(t),u^{\pi , y_0 }(t))\right]\ \ \ \ \forall \ \ {\rm Borel} \ Q\subset G
\end{equation}
for some
$\pi\in \Pi_{per}(y_0)$, with $(y^{\pi , y_0 }(\cdot),u^{\pi , y_0 }(\cdot))$ being the state-control trajectory generated by $\pi$  and with $T_0 $,  $\mathcal{T} $
being as in (\ref{e-opt-4}).
Due to (\ref{BB1-D-y-0-1}) and (\ref{A3-T-per}), the validity of (\ref{e-V-per-k}) will be established if we show that
\begin{equation}\label{e-proof-ineq}
 \G_{per}\subset W\cap D(y_0).
\end{equation}
It can be readily understood that
$$
\G_{per}(y_0)\subset \limsup_{T\rightarrow\infty}\Gamma_T(y_0)\subset W
$$
(the latter inclusion follows from (\ref{convergence-to-W-dis-2})). Consequently (see (\ref{eq-D-y-0})), to prove (\ref{e-proof-ineq}), it is sufficient to prove that, for any $\gamma\in \G_{per}(y_0)$,
there exists $\xi\in M_+(G)$ such that
\begin{equation}\label{e-proof-ineq-1}
 \int_G(\ph(y)-\ph(y_0))\gm=\int_G(E[\ph(f(y,u,s))]-\ph(y))\xm\ \ \ \ \forall \ \ph\in C(Y) .
\end{equation}
To prove that this is the case, observe that, if $\gamma\in \G_{per}(y_0)$, then, by (\ref{e-proof-ineq-01}),
$$
\int_{G}q(y,u)\gamma(dy,du)= {1\o \mathcal{T}} E\left[\sum_{t=T_0}^{T_0+\mathcal{T}-1} q(y^{\pi , y_0 }(t),u^{\pi , y_0 }(t))\right]
$$
for any continuous $q(y,u)$. Therefore, for any continuous $\ph(y)$,
$$
\int_{G}(\ph (y)-\ph (y_0))\gamma(dy,du)= \frac{1}{\mathcal{T}}\sum_{t=T_0}^{T_0+\mathcal{T}-1}E[\ph(y^{\pi, y_0}(t))-\ph(y_0)]
$$
\vspace{-0.4cm}
$$
=\frac{1}{\mathcal{T}}\sum_{t=T_0}^{T_0+\mathcal{T}-1}E[\ph (y^{\pi, y_0}(t))-\ph(y^{\pi, y_0}(T_0))] + E[\ph(y^{\pi, y_0}(T_0)) -\ph(y_0)]
$$
\vspace{-0.2cm}
$$
=\frac{1}{\mathcal{T}}\sum_{t=T_0}^{T_0+\mathcal{T}-1}\left(\sum_{t'=T_0}^{t-1}E[\ph(y^{\pi, y_0}(t'+1)) - \ph(y^{\pi, y_0}(t'))] \right) +
\sum_{t'=0}^{T_0-1}E[\left(\ph(y^{\pi, y_0}(t'+1)) - \ph(y^{\pi, y_0}(t')) \right)]
$$
\vspace{-0.2cm}
$$
=\frac{1}{\mathcal{T}}E\left[\sum_{t=T_0}^{T_0+\mathcal{T}-1}\left(\sum_{t'=T_0}^{t-1}\big(E[\ph(f(y^{\pi, y_0}(t'), u^{\pi, y_0}(t'), s(t')))|y^{\pi, y_0}(t')] - \ph(y^{\pi, y_0}(t'))\big) \right)\right]
$$
\vspace{-0.2cm}
$$
+ E\left[\sum_{t'=0}^{T_0-1}\big(E[\ph(f(y^{\pi, y_0}(t'), u^{\pi, y_0}(t'), s(t')))|y^{\pi, y_0}(t')] - \ph(y^{\pi, y_0}(t'))\big) \right].
$$
Since
$$
E[\ph(f(y^{\pi, y_0}(t'), u^{\pi, y_0}(t'), s(t')))|y^{\pi, y_0}(t')] = \bar \ph(y^{\pi, y_0}(t'), u^{\pi, y_0}(t'))  \ \ a.s. ,
$$
with
\begin{equation}\label{e-phi-bar}
\bar \ph (y,u) \BYDEF E[\ph(f(y,u,s))] ,
\end{equation}
we obtain
$$
\int_{G}(\ph (y)-\ph (y_0))\gamma(dy,du) = \frac{1}{\mathcal{T}}E\left[\sum_{t=T_0}^{T_0+\mathcal{T}-1}\left(\sum_{t'=T_0}^{t-1}\big( \bar \ph(y^{\pi, y_0}(t'), u^{\pi, y_0}(t')) - \ph(y^{\pi, y_0}(t'))\big) \right)\right]
$$
\vspace{-0.2cm}
\begin{equation}\label{e-phi-bar-1}
 +  E\left[\sum_{t'=0}^{T_0-1}\big(\bar \ph(y^{\pi, y_0}(t'), u^{\pi, y_0}(t')) - \ph(y^{\pi, y_0}(t'))\big) \right].
\end{equation}
Define $\zeta\in C^*(G)$ by the equation
$$
\langle \zeta ,q\rangle ={1\o \mathcal{T}}E\left[\sum_{t=T_0}^{T_0+\mathcal{T}-1}\sum_{t'=T_0}^{t-1}q(y^{\pi, y_0}(t'), u^{\pi, y_0}(t'))\right]
+ E\left[\sum_{t'=0}^{T_0-1}q(y^{\pi, y_0}(t'), u^{\pi, y_0}(t'))\right] \ \ \ \ \forall \ q\in C(G).
$$
Due to Riesz representation theorem, there exists $\xi \in \Mp(G)$ such that
$$
\langle \zeta,q\rangle =\int_G q(y,u)\xi (dy,du)  \ \ \ \ \forall \  q\in C(G).
$$
Therefore, (\ref{e-phi-bar-1}) can be rewritten as
$$
\int_{G}(\ph (y)-\ph (y_0))\gamma(dy,du)= \langle \zeta ,\bar \ph(y,u) - \ph  (y)\rangle = \int_{G}(\bar \ph(y,u)-\ph(y))\,\xi (dy,du).
$$

In accordance with the definition of $\bar \ph(y,u) $ (see (\ref{e-phi-bar})), the latter is equivalent to (\ref{e-proof-ineq-1}). This completes the 

\bigskip 
\bigskip

proof of the proposition. $\ \Box$

\begin{Corollary}\label{Cor-strong-dual-1}
If the conditions of Theorem \ref{ThN1} {\rm (a)} are satisfied and if the equality (\ref{A112-1-per-ineq-1}) is valid, then the strong duality equality (\ref{eq-strong-duality-1-1}) is true.
\end{Corollary}

{\bf Example 1 (continuation).}
Let us  show that the optimal control plan $\pi^*$ defined in (\ref{e-ex-1-1}) is PRG.  In fact, the use of $\pi^*$ leads to that, for any $t\geq 0$,
 $$
y^{\pi^*, y_0}(t+1)=\Bigg\{
\begin{array}{rrrl}
-y^{\pi^*, y_0}(t) \ & \ \ \ \ \ {\rm with\ probability}&  \ \ \ \ \ \ \ \ \ \ \ \ 3/4,
\\
y^{\pi^*, y_0}(t) \ & \ \ \ \  \ {\rm with\ probability}&  \ \ \ \ \ \ \ \ \ \ \ \  \ 1/4
\end{array}
$$
if $y^{\pi^*, y_0}(t)>0$, and
$$
y^{\pi^*, y_0}(t+1)=\Bigg\{
\begin{array}{rrrl}
-y^{\pi^*, y_0}(t) \ & \ \ \ \ \ {\rm with\ probability}&  \ \ \ \ \ \ \ \ \ \ \ \ 1/4,
\\
y^{\pi^*, y_0}(t) \ & \ \ \ \  \ {\rm with\ probability}&  \ \ \ \ \ \ \ \ \ \ \ \  \ 3/4
\end{array}
$$
if $y^{\pi^*, y_0}(t)<0$. Since $|y^{\pi^*, y_0}(t)|=|y_0|$ for any $t\geq 0$, we can come to the conclusion that
$$
y^{\pi^*, y_0}(t+1)=\Bigg\{
\begin{array}{rrrl}
-|y_0| \ & \ \ \ \ \ {\rm with\ probability}&  \ \ \ \ \ \ \ \ \ \ \ \ 3/4,
\\
|y_0| \ & \ \ \ \  \ {\rm with\ probability}&  \ \ \ \ \ \ \ \ \ \ \ \  \ 1/4
\end{array}
$$
for any $t\geq 0$. Thus,
\begin{equation*}\label{e-opt-4-11}
	E[q(y^{\pi^*, y_0}(t),u^{\pi^*, y_0}(t)) ]= \frac{3}{4}q(-|y_0|,1) + \frac{1}{4}q(|y_0|,-1)\ \ \  \forall \ t\geq 1,
	\end{equation*}
for any continuous $q(y,u)$.  Consequently, (\ref{e-opt-4}) is satisfied with $T_0=1$ and $\mathcal{T}=1$.

\medskip

{\bf Example 2 (continuation).}  The optimal control plan defined by (\ref{e-ex-2-extra-1}) is PRG for $y_0\in [-1,0]$ since in this case (see (\ref{e-extra-extra-3}))
$$
E[q(y(t),u(t))] =\frac{1}{2}q(-1,-1) + \frac{1}{2}q(-1,-\frac{1}{4}) \ \ \ \ \ \ \forall \ t\geq 1 .
$$
However, this control plan is not PRG for $y_0\in (0,1]$.

\medskip

\section{Sufficient and necessary conditions for long-run average optimality}\label{Sec-opt-cond}

DEFINITION. A pair $(\bar \psi(\cdot), \bar \eta(\cdot))\in \mathcal{B}(Y)\times \mathcal{B}(Y) $ will be called  an {\it optimal solution} of
(\ref{BB8-B}) if it satisfies the inequalities (compare with (\ref{BB7-B}))
\begin{equation}\label{BB7-3-1}
\begin{aligned}
&k(y,u)+(\bar\psi(y_0)-\bar\psi(y))+E[\bar\eta(f(y,u,s))]-\bar\eta(y)\ge \hat d^*(y_0),\\
&E[\bar\psi(f(y,u,s))]-\bar\psi(y)\ge 0
\end{aligned}
\end{equation}
for all $(y,u)\in G$.

\begin{Proposition}\label{PN2}	
{\rm (a)} If, for a given $y_0\in Y$,  the limit (\ref{chesaro-lim-1-1}) exists and
\begin{equation}\label{abel-lim-1-11}
	v(y_0) = \hat d^*(y_0),
	\end{equation}
then
 a pair  $ (\bar \psi(\cdot), \bar \eta(\cdot)) $ is an optimal solution of (\ref{BB8-B}) if and only if $\ \bar \psi(\cdot)$ satisfies the second inequality in (\ref{BB7-B}) and
	\begin{equation}\label{min}
	\begin{aligned}
	\inf_{(y,u)\in G}\{k(y,u)-\bp(y)+E[\bar\eta(f(y,u,s))]-\bar\eta(y)\}=v(y_0)-\bp(y_0).
	\end{aligned}
	\end{equation}
	{\rm (b)} Let the limit (\ref{chesaro-lim-1-1}) exist  and (\ref{abel-lim-1-11})  be valid for any $y_0\in Y$. If
$\ \bar \eta(\cdot)\in \mathcal{B}(Y)$   is such that
\begin{equation}\label{min-V}
	\begin{aligned}
	\inf_{(y,u)\in G}\{k(y,u)- v(y)+E[\bar\eta(f(y,u,s))]-\bar\eta(y)\}=0,
	\end{aligned}
	\end{equation}
 then the pair
$ ( \bar \psi(\cdot), \bar \eta(\cdot)) $, where $\bar \psi(\cdot) = v(\cdot) $
		is an optimal solution of problem  (\ref{BB8-B}).
	\end{Proposition}
{\bf Proof.} By (\ref{CC9-B}), the first inequality in (\ref{BB7-3-1}) is equivalent to the equality
\begin{equation}\label{min-0}
\inf_{(y,u)\in G}\{k(y,u)+\bp(y_0)-\bp(y)+E[\bar\eta(f(y,u,s))]-\bar\eta(y)\}=\hat d^*(y_0).
\end{equation}
Also,
(\ref{min-0}) is equivalent to  (\ref{min}) (due to  (\ref{abel-lim-1-11})). Thus {\rm (a)} is proved.

If $\bar \eta(\cdot) $ is such that (\ref{min-V}) is satisfied, then the pair $ ( \bar \psi(\cdot), \bar \eta(\cdot)) $, where $\bar \psi(\cdot) = v(\cdot) $, satisfies (\ref{min}). Therefore, due to {\rm (a)} and  due to the fact  that $v(\cdot)$ satisfies the second inequality in (\ref{BB7-B}) (see (\ref{V51-1}) below), this pair is an optimal solution of  (\ref{BB8-B}). This proves
{\rm (b)}. $\ \Box$

\medskip

Consider the  optimal control problem
\begin{equation}\label{e-LRAOC-1}
	\inf_{\pi \in \Pi} \limsup_{T\to \infty} {1\over T}E\left[\sum_{t=0}^{T-1} k(y^{\pi, y_0}(t),u^{\pi, y_0}(t))\right]\BYDEF V(y_0) .
	\end{equation}
As can be readily seen,
\begin{equation}\label{e-known}
V(y_0)\geq \liminf_{T\rightarrow\infty} v_{T}(y_0),
\end{equation}
and, therefore,
\begin{equation}\label{e-known-1}
V(y_0)\geq \hat d^*(y_0).
\end{equation}
 The following proposition gives sufficient conditions for a control plan $\pi$ to be optimal in (\ref{e-LRAOC-1}) and for the equality
 \begin{equation}\label{e-opt-3}
V(y_0)=\hat d^*(y_0)
\end{equation}
to be valid.

 \begin{Proposition}\label{Prop-optim-cond-suf}
	 Let  an optimal solution $(\bar\psi(\cdot),\bar\eta(\cdot))$  of  \eqref{BB8-B} exist.
	For a control plan $\pi \in \Pi$ with the associated state-control trajectory $(y^{\pi, y_0}(\cdot), u^{\pi, y_0}(\cdot))$  to be optimal in  problem (\ref{e-LRAOC-1}) and for the equality (\ref{e-opt-3}) to be valid,
it is sufficient that there exists an  integer $T_0>0$ such that
\begin{equation}\label{e-opt-1}
\begin{aligned}
&k(y^{\pi, y_0}(t),u^{\pi, y_0}(t))+(\bar\psi(y_0)-\bar\psi(y^{\pi, y_0}(t)))+E\big[\bar\eta(f(y^{\pi, y_0}(t),u^{\pi, y_0}(t),s(t)))|y^{\pi, y_0}(t)\big]\\ &\ \ \ \ \ \ \ \ \ \ \ \ \ \ \ \ \ \ \ \ \ \ \ \ \ \ \ \ \ \ \ \ \ \ \ \ \ \ \ \ \ \ \ \ \ \ \ \ \ \ \ \ \ \ \ \ \ \ \ \ \ \
\    -\bar\eta(y^{\pi, y_0}(t)) =  \hat d^*(y_0) \ \ a.s. \ \ \forall \ t \geq T_0 ,
\end{aligned}
\end{equation}
 and
	\begin{equation}\label{e-opt-2}
	E[\bar\psi(y^{\pi, y_0}(t)) ]= \bar\psi(y_0)\ \ \ \forall \ t\geq T_0.
	\end{equation}
	\end{Proposition}
{\bf Proof.}
Taking the mathematical expectation of the sum (from $t=T_0$ to $T-1$, with $T\geq T_0+1$) of the equations in (\ref{e-opt-1}) and taking into account
(\ref{e-opt-2}), one can obtain (similarly to (\ref{e-long})):
\begin{equation*}\label{e-long-1}
	\begin{split}
	\frac{1}{T}E\left[\sum_{t=T_0}^{T-1}k(y^{\pi, y_0}(t),u^{\pi, y_0}(t))\right] + \frac{1}{T}E\left[\sum_{t=T_0}^{T-1}E[\bar\eta\big(f(y^{\pi, y_0}(t), u^{\pi, y_0}(t), s(t))\big) | y^{\pi, y_0}(t)] - \bar\eta(y^{\pi, y_0}(t))\right] &\\
	=\frac{1}{T}E\left[\sum_{t=T_0}^{T-1}k(y^{\pi, y_0}(t),u^{\pi, y_0}(t))\right] + \frac{1}{T}E\left[\sum_{t=T_0}^{T-1}E[\bar\eta(y^{\pi, y_0}(t+1)) | y^{\pi, y_0}(t)] - \bar\eta(y^{\pi, y_0}(t))\right]\\
	=\frac{1}{T}E\left[\sum_{t=T_0}^{T-1}k(y^{\pi, y_0}(t),u^{\pi, y_0}(t))\right] + \frac{1}{T}E\left[\sum_{t=T_0}^{T-1}\bar\eta(y^{\pi, y_0}(t+1)) - \bar\eta(y^{\pi, y_0}(t))\right]\\
	=\frac{1}{T}E\left[\sum_{t=T_0}^{T-1}k(y^{\pi, y_0}(t),u^{\pi, y_0}(t))\right] + \frac{1}{T}E\left[\bar\eta(y^{\pi, y_0}(T)) - \bar\eta(y^{\pi, y_0}(T_0))\right]\\
= \frac{T-T_0}{T}\hat d^*(y_0).
	\end{split}
	\end{equation*}
Therefore,
\begin{equation*}\label{e-long-2}
	\begin{split}
\frac{1}{T}E\left[\sum_{t=0}^{T-1}k(y^{\pi, y_0}(t),u^{\pi, y_0}(t))\right]= \frac{T-T_0}{T}\hat d^*(y_0) + \frac{1}{T}E\left[\sum_{t=0}^{T_0-1}k(y^{\pi, y_0}(t),u^{\pi, y_0}(t))\right]&\\
-\frac{1}{T}E\left[\bar\eta(y^{\pi, y_0}(T)) - \bar\eta(y^{\pi, y_0}(T_0))\right],
\end{split}
	\end{equation*}
and, consequently,
$$
\lim_{T\rightarrow\infty}\frac{1}{T}E\left[\sum_{t=0}^{T-1}k(y^{\pi, y_0}(t),u^{\pi, y_0}(t))\right] = \hat d^*(y_0).
$$
Thus, by (\ref{e-known-1}), $\pi$ is optimal in (\ref{e-LRAOC-1}), and (\ref{e-opt-3}) is valid. $\ \Box$

\medskip

Let us now establish that the fulfillment  of (\ref{e-opt-1}), (\ref{e-opt-2})  is also a necessary
condition for a PRG control plan $\pi$ to be optimal and for (\ref{e-opt-3}) to be valid.

\begin{Proposition}\label{Prop-nec-opt-cond}
Let  an optimal solution $(\bar\psi(\cdot),\bar\eta(\cdot))$  of  \eqref{BB8-B} exist.
 If a PRG control plan $\pi$ (that is, a control plan such that (\ref{e-opt-4}) is satisfied) is optimal in (\ref{e-LRAOC-1})  and if (\ref{e-opt-3}) is valid,
 then (\ref{e-opt-1}) and (\ref{e-opt-2}) are satisfied for any $t\geq T_0$, where $T_0$ is as in
(\ref{e-opt-4}).
\end{Proposition}
{\bf Proof.}
Due to optimality of the PRG control plan $\pi$ and due to (\ref{e-opt-4}), (\ref{e-opt-3}),
\begin{equation}\label{e-opt-5}
\hat d^*(y_0)=\lim_{T\rightarrow\infty}\frac{1}{T}E\left[\sum_{t'=0}^{T-1}	k(y^{\pi, y_0}(t'),u^{\pi, y_0}(t')) \right]= \frac{1}{\mathcal{T}}E\left[\sum_{t'=t}^{t+\mathcal{T}-1}	k(y^{\pi, y_0}(t'),u^{\pi, y_0}(t')) \right]\ \ \ \ \  \forall \ t\geq T_0.
	\end{equation}
By (\ref{BB7-3-1}), the following inequalities hold true:
\begin{equation}\label{BB7-3-a.s.}
\begin{aligned}
&k(y^{\pi, y_0}(t'),u^{\pi, y_0}(t'))+(\bar\psi(y_0)-\bar\psi(y^{\pi, y_0}(t')))+E\big[\bar\eta(f(y^{\pi, y_0}(t'),u^{\pi, y_0}(t'),s(t')))|y^{\pi, y_0}(t')\big]\\ &\ \ \ \ \ \ \ \ \ \ \ \ \ \ \ \ \ \ \ \ \ \ \ \ \ \ \ \ \ \ \ \ \ \ \ \ \ \ \ \ \ \ \ \ \ \ \ \ \ \ \ \ \ \ \ \ \ \ \ \ \ \ \ \ \ \ \ \ \
\ \ \ \ \ \ \ \ \ \ \ \ \ \ \ \   -\bar\eta(y^{\pi, y_0}(t')) \ge  \hat d^*(y_0) \ \ \ a.s. \ \ \forall\ t'= 0,1,...
\end{aligned}
\end{equation}
\begin{equation}\label{BB7-3-a.s.1}
E[\bar\psi(f(y^{\pi, y_0}(t'),u^{\pi, y_0}(t'),s(t')))|y^{\pi, y_0}(t')]-\bar\psi(y^{\pi, y_0}(t'))\ge 0  \ \ \ a.s. \ \ \forall\ t'= 0,1,... \ .
\end{equation}
 Take an arbitrary $t\geq T_0$ and take the mathematical expectation of the sum of the inequalities (\ref{BB7-3-a.s.}) from
$t'=t$ to $t'=t+\mathcal{T}-1$. Using (\ref{e-opt-5}), one  obtains
\begin{equation}\label{BB7-3-a.s.2}
\sum_{t'=t}^{t+\mathcal{T}-1}E[\bar\psi(y_0)-\bar\psi(y^{\pi, y_0}(t'))] + E\left[\sum_{t'=t}^{t+\mathcal{T}-1} \big(E\big[\bar\eta(f(y^{\pi, y_0}(t'),u^{\pi, y_0}(t'),s(t')))|y^{\pi, y_0}(t')\big]-\bar\eta(y^{\pi, y_0}(t'))\big)  \right]\geq 0.
\end{equation}
By (\ref{e-opt-4}),
\begin{equation}\label{BB7-3-a.s.5}
\begin{aligned}
 E\left[\sum_{t'=t}^{t+\mathcal{T}-1} \big(E\big[\bar\eta(f(y^{\pi, y_0}(t'),u^{\pi, y_0}(t'),s(t')))|y^{\pi, y_0}(t')\big]-\bar\eta(y^{\pi, y_0}(t'))\big)  \right]
 \\
 = E\left[\sum_{t'=t}^{t+\mathcal{T}-1} \big(\bar\eta(y^{\pi, y_0}(t'+1))-\bar\eta(y^{\pi, y_0}(t'))\big)  \right] = E[\bar\eta(y^{\pi, y_0}(t+\mathcal{T}))]&
 -E[\bar\eta(y^{\pi, y_0}(t))]=0.
\end{aligned}
\end{equation}
Hence, from (\ref{BB7-3-a.s.2}) it follows that
\begin{equation}\label{BB7-3-a.s.3}
\sum_{t'=t}^{t+\mathcal{T}-1}E[\bar\psi(y_0)-\bar\psi(y^{\pi, y_0}(t'))] \geq 0.
\end{equation}
From (\ref{BB7-3-a.s.1}), on the other hand, it follows that
$$
E[\bar\psi(y^{\pi, y_0}(t'+1))- \bar\psi(y^{\pi, y_0}(t'))] \geq 0\ \ \forall\ t'=0,1,... \ \ \Rightarrow \ \ E[\bar\psi(y^{\pi, y_0}(t'))- \bar\psi(y_0)] \geq 0
\ \ \forall\ t'=0,1,...\ .
$$
The latter and (\ref{BB7-3-a.s.3}) imply that
\begin{equation}\label{BB7-3-a.s.4}
E[\bar\psi(y_0)-\bar\psi(y^{\pi, y_0}(t'))] = 0 \ \ \ \ \forall\ t'= t, t+1,...,t+\mathcal{T}-1,
\end{equation}
which, in turn, implies (\ref{e-opt-2}) (since $t$ is an arbitrary integer that is greater or equal than $T_0$).

Let us now show that the inequality (\ref{BB7-3-a.s.}) is  satisfied a.s. as the equality for any $t'=t\geq T_0$. Assume it is not the case. Then
\begin{equation*}
\begin{aligned}
&E\big[k(y^{\pi, y_0}(t),u^{\pi, y_0}(t))+(\bar\psi(y_0)-\bar\psi(y^{\pi, y_0}(t)))+E\big[\bar\eta(f(y^{\pi, y_0}(t),u^{\pi, y_0}(t),s(t)))|y^{\pi, y_0}(t)\big]\\ &\ \ \ \ \ \ \ \ \ \ \ \ \ \ \ \ \ \ \ \ \ \ \ \ \ \ \ \ \ \ \ \ \ \ \ \ \ \ \ \ \ \ \ \ \ \ \ \ \ \ \ \ \ \ \ \ \ \ \ \ \ \ \ \ \ \ \ \ \
\ \ \ \ \ \ \ \ \ \ \ \ \ \ \ \   -\bar\eta(y^{\pi, y_0}(t))\big] > \hat d^*(y_0),
\end{aligned}
\end{equation*}
and, consequently,
\begin{equation*}\label{BB7-3-a.s.6}
\begin{aligned}
 \frac{1}{\mathcal{T}}E\left[\sum_{t'=t}^{t+\mathcal{T}-1}k(y^{\pi, y_0}(t'),u^{\pi, y_0}(t')) \right]
+ \frac{1}{\mathcal{T}}E\left[\sum_{t'=t}^{t+\mathcal{T}-1}\big(\bar\psi(y_0)-\bar\psi(y^{\pi, y_0}(t'))\big)\right]
 \\
 + \frac{1}{\mathcal{T}}E\left[\sum_{t'=t}^{t+\mathcal{T}-1} \big(E\big[\bar\eta(f(y^{\pi, y_0}(t'),u^{\pi, y_0}(t'),s(t')))|y^{\pi, y_0}(t')\big]-\bar\eta(y^{\pi, y_0}(t'))\big)  \right] > \hat d^*(y_0).
\end{aligned}
\end{equation*}
By virtue of (\ref{BB7-3-a.s.5}) and (\ref{BB7-3-a.s.4}), this leads to
$$
\frac{1}{\mathcal{T}}E\left[\sum_{t'=t}^{t+\mathcal{T}-1}k(y^{\pi, y_0}(t'),u^{\pi, y_0}(t')) \right]> \hat d^*(y_0),
$$
which contradicts (\ref{e-opt-5}). This contradiction proves the required statement. $\ \Box$

\medskip

REMARK. In accordance with Propositions \ref{Prop-optim-cond-suf} and \ref{Prop-nec-opt-cond}, for a control plan $\pi$ to be optimal it is sufficient and also necessary (if it is PRG) that
the corresponding state-control trajectory $(y^{\pi, y_0}(t),u^{\pi, y_0}(t))$  satisfies the equalities (\ref{e-opt-1}), (\ref{e-opt-2}) for $t\geq T_0$ ($T_0$ being some positive integer). Note that, by (\ref{min}), the equality (\ref{e-opt-1}) is equivalent to
$$
(y^{\pi, y_0}(t),u^{\pi, y_0}(t))={\rm argmin}_{(y,u)\in G}\{k(y,u)-\bar\psi(y)+ \bar{\bar \eta}(y,u)-\bar\eta(y)\} \ \ a.s.\ \ \forall\ t\geq T_0 ,
$$
where $\bar{\bar \eta}(y,u)\BYDEF E[\bar\eta (f(y,u,s))] $.
This leads to
$$
u^{\pi, y_0}(t)={\rm argmin}_{u\in U(y^{\pi, y_0}(t))}\{k(y^{\pi, y_0}(t),u)+ \bar{\bar \eta}(y^{\pi, y_0}(t),u)\}  \ \ a.s.\ \ \forall\ t\geq T_0,
$$
which, in turn, implies that the  feedback control
\begin{equation}\label{e-feedback}
\begin{aligned}
u^f(y)\BYDEF{\rm argmin}_{u\in U(y)}\{k(y ,u)+ \bar{\bar \eta}(y,u)\} = {\rm argmin}_{u\in U(y)}\{k(y ,u)+E[\bar\eta (f(y,u,s))]\}
\end{aligned}
\end{equation}
is optimal for $T\geq T_0$ provided that the solution of system  (\ref{A1}) obtained with the use of this control satisfies (\ref{e-opt-2}).

\medskip

{\bf Example 1 (continuation).} The pair of functions $(\bar \psi(y), \bar \eta(y)) $ defined in (\ref{BB1-ex-8}) is
an optimal solution of the dual problem  (\ref{BB1-ex-4}) in this case, and   it can be readily verified that
$$
E[\bar\eta (f(y,u,s))] = E[\bar\eta (yus)]= \frac{3}{4} \bar\eta (yu) + \frac{1}{4} \bar\eta (-yu)= \frac{1}{2} yu +  \frac{1}{2}|y|.
$$
Therefore, the feedback control  (\ref{e-feedback}) takes  the form
\begin{equation*}\label{e-feedback-ex-1}
u^f(y)= {\rm argmin}_{u\in \{-1,1\}}\{y+\frac{1}{2} yu +  \frac{1}{2}|y|\}= \Bigg\{
\begin{array}{rrrl}
 1 \ & \ \ \ \ \ {\rm if}&  \ \ \ \ \ \ \ \ \ \ \ \ y\in [-1,0)
\\
-1 \ & \ \ \ \  \ {\rm if}&  \ \ \ \ \ \ \ \ \ \ \ \  \ y\in (0,1].
\end{array}%
\end{equation*}
 That is, by using formula (\ref{e-feedback}), we obtain the optimal
 control plan $\pi^*$  (see (\ref{e-ex-1-1})). Since, due to (\ref{example-1-0}),  $|y(t)|=|y_0|$ for any $t=0,1,...$, the equality   (\ref{e-opt-2})  is satisfied automatically in this example (see (\ref{BB1-ex-8})).

\medskip

{\bf Example 2 (continuation).} The pair of functions $(\bar \psi(y), \bar \eta(y)) $  defined in  (\ref{e-ex-2-3}) is an optimal solution of the dual problem  and,
having in mind that $E[\ph (us)]=\frac{1}{2}\ph (u)+ \frac{1}{2}\ph (\frac{u}{4})  $,  we obtain (see (\ref{e-ex-2-3}))
\begin{equation*}\label{e-extra-extra-1}
\begin{aligned}
&E[\bar\eta(us)]= \frac{1}{2}\Big(u+\frac{5}{8}\Big) + \frac{1}{2}\Big(\frac{u}{4}+\frac{5}{8}\Big)=\frac{5}{8}u + \frac{5}{8} \ \ \ \ {\rm for }\ \ \ \ u\in [-1,0],\\
&E[\bar\eta(us)]= \frac{1}{2} \Big(\frac{8}{3}u\Big) +  \frac{1}{2} \Big(\frac{8}{12}u\Big)= \frac{5}{3}u\ \ \ \ {\rm for }\ \ \ \ u\in (0,1].
\end{aligned}
\end{equation*}
Hence, in accordance with (\ref{e-feedback}),
\begin{equation*}\label{e-extra-extra-2}
\begin{aligned}
&u^f(y)= {\rm argmin}_{u\in [-1,y]}\Big\{y +\frac{5}{8}u + \frac{5}{8}\Big\}= -1\ \ \ \ {\rm if }\ \ \ \ y\in [-1,0),\\
& u^f(y)= {\rm argmin}_{u\in [-1,1]}\Big\{y +\frac{5}{8}u + \frac{5}{8}\Big\}= -1\ \ \ \ {\rm if }\ \ \ \ y=0,                                                                                          \\
&u^f(y)= {\rm argmin}_{u\in [y,1]}\Big\{y+\frac{5}{3}u\Big\}= y\ \ \ \ {\rm if }\ \ \ \ y\in (0,1].
\end{aligned}
\end{equation*}
 That is,   (\ref{e-feedback}) defines the optimal control plan $\pi^*$ in this example too (see (\ref{e-ex-2-extra-1})). Note  that the  state trajectory
 obtained with the use of this control plan satisfies (\ref{e-opt-2}) (see  (\ref{e-extra-extra-3})  and (\ref{e-ex-2-3})).

 \medskip

In Example 1,  the pair consisting of the first  term and the second term multiplied by $T$  in the expressions for the optimal value $v_T(\cdot)$ (see (\ref{e-ex-1-2})) coincides with 
the optimal solution (\ref{BB1-ex-8}) of the dual problem (\ref{BB1-ex-4}). Similarly, in Example 2, the pair consisting of the first  term and the limit
as $T\to\infty$ of the second term (multiplied by $T$) in the expressions  for the optimal value $v_T(\cdot)$  (see (\ref{eq-union-6-0})) coincides with 
the optimal solution (\ref{e-ex-2-3}) of the dual problem (\ref{CC9-B-ex-2}).

We conclude this section with a statement that gives sufficient conditions for this to be true in the general case. (Note that all these conditions are satisfied in both Example 1 and Example 2.)

\begin{Proposition}\label{Prop-dual-sol}
Assume that the following conditions are satisfied:

(i) The optimal value functions $v_T(y_0)$ is presented in the form
\begin{equation}\label{e-last-1}
v_T(y_0)= v(y_0) + \frac{1}{T}\eta_T(y_0) \ \ \ \forall \ y_0\in G;
\end{equation}
(ii) The equality (\ref{abel-lim-1-11}) is valid for any $y_0\in Y$ and
\begin{equation}\label{e-last-2}
E[v(f(y_0,u,s))] = v(y_0)  \ \ \ \forall \ u\in U(y_0),\ \ \forall \ y_0\in G;
\end{equation}
(iii) The function $\eta_T(\cdot) $ converges in the uniform metric to a function $\bar \eta(\cdot)$; that is
\begin{equation}\label{e-last-3}
\lim_{T\rightarrow\infty}\sup_{y\in Y }|\eta_T(y)-\bar \eta (y)| = 0 .
\end{equation}
Then the pair $(\bar \psi(\cdot),  \bar \eta (\cdot))$, where $\bar \psi(\cdot) = v(\cdot) $, is an optimal solution of the dual problem
 (\ref{BB8-B}).
\end{Proposition}
{\bf Proof.} Firstly, note that from (\ref{e-App-1}) it follows that
\begin{equation}\label{e-last-4}
\inf_{(y,u)\in G}\{-T v_T(y)+ k(y,u) + (T-1)E[v_{T-1}(f(y,u,s))]\}=0 , \ \ \ T=1,2,... \ .
\end{equation}
By substituting (\ref{e-last-1}) into (\ref{e-last-4}), we obtain
$$
\inf_{(y,u)\in G}\{-T v(y) - \eta_T(y) + k(y,u) + (T-1)E[v(f(y,u,s))] + E[\eta_{T-1}(f(y,u,s))]\}=0 , \ \ \ T=1,2,... \ ,
$$
which, after taking into account (\ref{e-last-2}), leads to
$$
\inf_{(y,u)\in G}\{k(y,u) -v(y)   + E[\eta_{T-1}(f(y,u,s))]- \eta_T(y)\}=0 , \ \ \ T=1,2,... \ .
$$
Passing to the limit with $T\rightarrow\infty$ in the equality above implies the validity of (\ref{min-V}) (due to (\ref{e-last-3})). The statement
follows now from Proposition \ref{PN2}{\rm (b)}. $\ \Box$

\section{Proofs of Theorems \ref{Th-lower-upper-limit-main} and  \ref{ThN1}}\label{Sec-Proof-Ab-Ch-estimates}

The proof of Theorem \ref{Th-lower-upper-limit-main} follows from Propositions \ref{Prop-est-from-below} and \ref{Prop-est-from-above} that are stated and proved below.

\begin{Proposition}\label{Prop-est-from-below} The following estimates are valid:
\begin{equation}\label{est-from-below}
\begin{aligned}
\liminf_{T\to \infty} v_T(y_0)\geq \hat d^*(y_0) \ \ \forall \ y_0\in Y, \\
\liminf_{\e\rightarrow 0} h_{\e}(y_0)\geq \hat d^*(y_0) \ \ \ \forall \ y_0\in Y.\\
\end{aligned}
\end{equation}
\end{Proposition}
{\bf Proof.} Let us prove the first inequality in (\ref{est-from-below}). Assume it is not true, that is, $\hat d^*(y_0) > \liminf_{T \rightarrow \infty}v_T(y_0) $. Then, there exists $\beta > 0$ and a pair of functions $(\psi(\cdot), \eta(\cdot)) \in \mathcal{B}(Y) \times \mathcal{B}(Y)$, with $\psi(\cdot)$ satisfying the second inequality in (\ref{BB7-B}), such that
	\begin{equation*}
	k(y,u) + (\psi(y_0) - \psi(y)) + E[\eta(f(y,u,s))] - \eta(y) \ge \liminf_{T' \rightarrow \infty}v_{T'}(y_0) + \beta,
	\end{equation*}
	for all $(y,u) \in G$.
From the latter, it follows that for any control plan $\pi \in \Pi$ with the associated state-control trajectory $(y^{\pi, y_0}(\cdot), u^{\pi, y_0}(\cdot))$,
	\begin{equation}\label{e-LB-1}
	\begin{split}
	E\big[k(y^{\pi, y_0}(t), u^{\pi, y_0}(t)) +(\psi(y_0) - \psi(y^{\pi, y_0}(t))) + \bar{\eta}(y^{\pi, y_0}(t), u^{\pi, y_0}(t)) - \eta(y^{\pi, y_0}(t))\big]& \\
	\ge \liminf_{T' \rightarrow \infty}v_{T'}(y_0) + \beta \ \ \forall \ t = 0,1,\ldots&
	\end{split}
	\end{equation}
	where
\begin{equation}\label{e-bar-eta}
\bar{\eta}(y,u) \BYDEF E[\eta(f(y,u,s))].
\end{equation}
 Due to the fact that $\psi(\cdot)$ satisfies the second inequality in (\ref{BB7-B}),
\begin{equation*}
	\begin{split}
	E\big[ \psi(y^{\pi, y_0}(t+1)) - \psi(y^{\pi, y_0}(t))\big]  = E\big[ E[\psi(y^{\pi, y_0}(t+1))|y^{\pi, y_0}(t) ] - \psi(y^{\pi, y_0}(t))\big]& \\
	= E\big[ E\big[\psi\big(f(y^{\pi, y_0}(t), u^{\pi, y_0}(t), s(t))\big)|y^{\pi, y_0}(t) \big] - \psi(y^{\pi, y_0}(t))\big]\geq 0 \ \ \forall \ t = 0,1,\ldots&\ .
	\end{split}
	\end{equation*}
Consequently,
\begin{equation}\label{e-LB-2}
E\big[ \psi(y^{\pi, y_0}(t))- \psi(y_0)\big]\geq 0  \ \ \forall \ t = 0,1,\ldots \ ,
\end{equation}
and from (\ref{e-LB-1}) it follows that
\begin{equation*}
	\begin{split}
	E\big[k(y^{\pi, y_0}(t), u^{\pi, y_0}(t))  + \bar{\eta}(y^{\pi, y_0}(t), u^{\pi, y_0}(t)) - \eta(y^{\pi, y_0}(t))\big]& \\
	\ge \liminf_{T' \rightarrow \infty}v_{T'}(y_0) + \beta \ \ \forall \ t = 0,1,\ldots&\ ,
	\end{split}
	\end{equation*}
which, in accordance with (\ref{e-bar-eta}), is equivalent to
\begin{equation*}
	\begin{split}
	E[k(y^{\pi, y_0}(t), u^{\pi, y_0}(t))] + E\big[E[\eta\big(f(y^{\pi, y_0}(t), u^{\pi, y_0}(t), s(t))\big) | y^{\pi, y_0}(t)] - \eta(y^{\pi, y_0}(t))\big] & \\
	\ge \liminf_{T' \rightarrow \infty}v_{T'}(y_0) + \beta \ \ \forall \ t = 0,1,\ldots&\ .
	\end{split}
	\end{equation*}
From these inequalities it follows that
	\begin{equation}\label{e-long}
	\begin{split}
	\frac{1}{T}E\left[\sum_{t=0}^{T-1}k(y^{\pi, y_0}(t),u^{\pi, y_0}(t))\right] + \frac{1}{T}E\left[\sum_{t=0}^{T-1}E[\eta\big(f(y^{\pi, y_0}(t), u^{\pi, y_0}(t), s(t))\big) | y^{\pi, y_0}(t)] - \eta(y^{\pi, y_0}(t))\right] &\\
	=\frac{1}{T}E\left[\sum_{t=0}^{T-1}k(y^{\pi, y_0}(t),u^{\pi, y_0}(t))\right] + \frac{1}{T}E\left[\sum_{t=0}^{T-1}E[\eta(y^{\pi, y_0}(t+1)) | y^{\pi, y_0}(t)] - \eta(y^{\pi, y_0}(t))\right]\\
	=\frac{1}{T}E\left[\sum_{t=0}^{T-1}k(y^{\pi, y_0}(t),u^{\pi, y_0}(t))\right] + \frac{1}{T}E\left[\sum_{t=0}^{T-1}\eta(y^{\pi, y_0}(t+1)) - \eta(y^{\pi, y_0}(t))\right]\\
	=\frac{1}{T}E\left[\sum_{t=0}^{T-1}k(y^{\pi, y_0}(t),u^{\pi, y_0}(t))\right] + \frac{1}{T}E\left[\eta(y^{\pi, y_0}(T)) - \eta(y^{\pi, y_0}(0))\right]\\
	\ge \liminf_{T' \rightarrow \infty}v_{T'}(y_0) + \beta. &
	\end{split}
	\end{equation}
Hence,
\begin{equation*}
	\frac{1}{T}E\left[\sum_{t=0}^{T-1}k(y^{\pi, y_0}(t),u^{\pi, y_0}(t))\right] + \frac{1}{T}\left(\sup_{y\in Y}\eta(y) - \eta(y_0)\right) \ge \liminf_{T' \rightarrow \infty}v_{T'}(y_0) + \beta.
	\end{equation*}
Since the above inequality holds for any control plan $\pi \in \Pi$, we may conclude that
	\begin{equation*}
	v_T(y_0) + \frac{1}{T}\left(\sup_{y\in Y}\eta(y) - \eta(y_0)\right) \ge \liminf_{T' \rightarrow \infty}v_{T'}(y_0) + \beta.
	\end{equation*}
	By taking $\liminf_{T \rightarrow \infty}$ in the left-hand side of this expression, one obtains a contradiction.
	Thus, the first inequality in (\ref{est-from-below}) is proved.

Let us now prove the second inequality in (\ref{est-from-below}). Assume it is not true. Then there exists $\beta > 0$ and a pair of functions $(\psi(\cdot), \eta(\cdot)) \in \mathcal{B}(Y) \times \mathcal{B}(Y)$, with $\psi(\cdot)$ satisfying the second inequality of (\ref{BB7-B}) such that
	\begin{equation*}
	k(y,u) + (\psi(y_0) - \psi(y)) + E[\eta(f(y,u,s))] - \eta(y) \ge \liminf_{\e' \rightarrow 0}  h_{\e'}(y_0) + \beta.
	\end{equation*}
	It follows that, for any control plan $\pi \in \Pi$ with the associated state-control trajectory $(y^{\pi, y_0}(\cdot), u^{\pi, y_0}(\cdot))$,
\begin{equation*}
	\begin{split}
	E\big[k(y^{\pi, y_0}(t), u^{\pi, y_0}(t)) +(\psi(y_0) - \psi(y^{\pi, y_0}(t))) + \bar{\eta}(y^{\pi, y_0}(t), u^{\pi, y_0}(t)) - \eta(y^{\pi, y_0}(t))\big]& \\
	\ge \liminf_{\e' \rightarrow 0} h_{\e'}(y_0) + \beta \ \ \forall \ t = 0,1,\ldots\ ,&
	\end{split}
	\end{equation*}
	where $\bar{\eta}(y,u)$ is as in (\ref{e-bar-eta}). By (\ref{e-LB-2}), the latter implies that
\begin{equation*}
	\begin{split}
	E\big[k(y^{\pi, y_0}(t), u^{\pi, y_0}(t))  + \bar{\eta}(y^{\pi, y_0}(t), u^{\pi, y_0}(t)) - \eta(y^{\pi, y_0}(t))\big]& \\
	\ge \liminf_{\e' \rightarrow 0} h_{\e'}(y_0) + \beta \ \ \forall \ t = 0,1,\ldots\ ,&
	\end{split}
	\end{equation*}
which, in view of (\ref{e-bar-eta}), can be rewritten as follows
	\begin{equation*}
	E[k(y^{\pi, y_0}(t),u^{\pi, y_0}(t))] + E\big[E[\eta\big(f(y^{\pi, y_0}(t), u^{\pi, y_0}(t), s(t))\big) | y^{\pi, y_0}(t)\big] - \eta(y^{\pi, y_0}(t))] \ge \liminf_{\e' \rightarrow 0}  h_{\e'}(y_0) + \beta.
	\end{equation*}
Hence, for any $\e\in (0,1) $,
\begin{equation}
	\label{eq:TFPP}
	\begin{split}
	\e E\left[\sum_{t=0}^{\infty}(1-\e)^tk(y^{\pi, y_0}(t),u^{\pi, y_0}(t))\right] + \e\sum_{t=0}^{\infty}(1-\e)^tE\left[E[\eta(y^{\pi, y_0}(t), u^{\pi, y_0}(t), s(t)) | y^{\pi, y_0}(t)] - \eta(y^{\pi, y_0}(t))\right] &\\
	=\e E\left[\sum_{t=0}^{\infty}(1-\e)^tk(y^{\pi, y_0}(t),u^{\pi, y_0}(t))\right] + \e\sum_{t=0}^{\infty}(1-\e)^tE\left[E[\eta(y^{\pi, y_0}(t+1))| y^{\pi, y_0}(t)]-\eta(y^{\pi, y_0}(t))\right] &\\
	=\e E\left[\sum_{t=0}^{\infty}(1-\e)^tk(y^{\pi, y_0}(t),u^{\pi, y_0}(t))\right] + \e\sum_{t=0}^{\infty}(1-\e)^tE\left[\eta(y^{\pi, y_0}(t+1))-\eta(y^{\pi, y_0}(t))\right] &\\
	\ge \liminf_{\e' \rightarrow 0}  h_{\e'}(y_0) + \beta.	 &
	\end{split}
	\end{equation}
One may verify that
	\begin{equation*}
	\label{eq:V_1}
	\begin{split}
	&\sum_{t=0}^{k}(1-\e)^tE\left[\eta(y^{\pi, y_0}(t+1))-\eta(y^{\pi, y_0}(t))\right] \\
	&= E\left[-\eta(y_0)+(1-\e)^k\eta(y^{\pi, y_0}(k+1)) + \e\sum_{t=1}^{k}(1-\e)^{t-1}\eta(y^{\pi, y_0}(t))\right]
	\end{split}
	\end{equation*}
	for any $k=1,2,\ldots$. This implies the validity of the following inequalities:
	\begin{equation}
	\label{eq:V_2}
	\begin{split}
	&\e\left|\sum_{t=0}^{k}(1-\e)^tE[\eta(y^{\pi, y_0}(t+1))-\eta(y^{\pi, y_0}(t))] \right|\\
	&\le \e E\left[|\eta(y_0)| + (1-\e)^k|\eta(y^{\pi, y_0}(k+1))| + \e\sum_{t=1}^{k}(1-\e)^{t-1}|\eta(y^{\pi, y_0}(t))|\right]\le 3\e M_\eta  ,
	\end{split}
	\end{equation}
where $M_\eta \BYDEF \sup_{y\in Y}|\eta(y)| $. From (\ref{eq:TFPP}) and (\ref{eq:V_2}) it follows that
$$
\e E\left[\sum_{t=0}^{\infty}(1-\e)^tk(y^{\pi, y_0}(t),u^{\pi, y_0}(t))\right] + 3\e M_\eta\geq \liminf_{\e' \rightarrow 0}  h_{\e'}(y_0) + \beta.
$$
Since this inequality is valid for  any control plan $\pi \in \Pi$, we may conclude that
\begin{equation*}
	 h_{\e}(y_0) + 3\e M_\eta \ge \liminf_{\e' \rightarrow 0}  h_{\e'}(y_0)+ \beta.
	\end{equation*}
	By taking $\liminf_{\e \rightarrow 0}$ in the left-hand-side of the expression above, we obtain a contradiction. This proves the second inequality in (\ref{est-from-below}).
$\ \Box$

The proof of the estimates from above in (\ref{BB3-relaxed-genearl}) is based on the following lemma.

\begin{Lemma}\label{Lem-upper-bound}
For any natural $T$,
\begin{equation}\label{e-before-lim-1}
\int_{G} v_T(y)\,\g(dy,du)\le   \int_{G} k(y,u)\,\g(dy,du)\ \ \ \forall \ \gamma\in W.
\end{equation}
Also, for any $\e\in (0,1) $,
\begin{equation}\label{e-before-lim-2}
\int_{G} h_{\e}(y)\,\g(dy,du)\le   \int_{G} k(y,u)\,\g(dy,du) \ \ \ \forall \ \gamma\in W.
\end{equation}
\end{Lemma}
{\bf Proof.}
Take an arbitrary $\g \in W$. It is well known that the measure $\g$ can
be represented as follows:
$$
\g(dy,du) = \pi(du|y) \g_1(dy),
$$
where $\g_1$ is the marginal of $\g$ on $Y$.
Note that the stochastic kernel $\pi(du|y)$ can be associated with a randomized stationary policy,
under which the induced state-control
 process has
the invariant probability measure $\g$, with
 $\g_1$
being the corresponding invariant probability  measure  of the state process  under the policy $\pi$. That is,
$$
\g_1(Q) = \int_{Y} P_\pi(Q|y) \g_1(dy),
\quad \forall \ \mbox{Borel} \ Q \subset Y,
$$
where
$$
P_\pi(Q|y) = \int_{\hat{U}} P(Q|y,u) \pi(du|y) ,
$$
$P(Q|y,u)$ being defined in (\ref{e-continuous-Borel-1}).

Under our assumptions about the model,
$v_T(y)$, which originally was defined as the optimal value function for \lq\lq non-randomized" controls, will remain such in case  the randomized controls are allowed (see, e.g., Sections~3.2,~3.3 in \cite{Jean-H-2} or Theorem 2(iii) in \cite{Feinbergatal2012}). Therefore,
$\ T \int_{G} v_T(y)\,\g(dy,du)= T \int_{Y} v_T(y)\,\g_1(dy)$  can be interpreted as the optimal $T$-stage cost for
the initial distribution $\g_1$, and
\begin{equation}
\label{e-before-lim-1-T}
T \int_{G} v_T(y)\,\g(dy,du) \le T \int_{G} k(y,u)\,\g(dy,du)
\end{equation}
since the right-hand-side in the expression above  represents
the total cost after $T$ stages provided that the process starts from the initial distribution
$\g_1$ and it is controlled by the stationary policy $\pi$. By dividing (\ref{e-before-lim-1-T}) by $T$, we obtain (\ref{e-before-lim-1}).

Let us now prove  (\ref{e-before-lim-2}). The optimal value function in  the infinite  horizon problem with time discounting, $h_{\e}(y)$, being defined  for the class of non-randomized controls, remains such in case the randomized controls 
are allowed (see, e.g., Section 4.2 in \cite{Jean-H-2} or Theorem 2(v) in \cite{Feinbergatal2012}). Therefore,
$\  \e^{-1}\int_{G} h_{\e}(y)\,\g(dy,du) =   \e^{-1}\int_{Y}  h_{\e}(y)\,\g_1(dy)$  is the optimal infinite horizon cost
for the initial distribution $\g_1$. Since the expression $ \ \e^{-1}\int_{G} k(y,u)\,\g(dy,du)$
  represents
the total infinite horizon discounted cost in case the process starts from the initial distribution
$\g_1$ and  the stationary control policy $\pi$ is used, we may write down
$$
\e^{-1}\int_{G} h_{\e}(y)\,\g(dy,du)\le   \e^{-1} \int_{G} k(y,u)\,\g(dy,du) .
$$
Multiplying the latter by $\e$, we obtain (\ref{e-before-lim-2}).
$\ \Box$

\medskip

REMARK. The proof above was proposed by an anonymous reviewer. It is much shorter than the original authors' proof, which was similar to the proof of Lemma 3.2 in \cite{BGS-2019}. Note that, the latter, in contrast to the former, was not based on    results    that establish the optimality of non-randomized controls.

\begin{Proposition}\label{Prop-est-from-above}  The following estimates are valid:
\begin{equation}\label{est-from-above}
\begin{aligned}
\limsup_{T\to \infty} v_T(y_0)\leq  k^*(y_0) \ \ \forall \ y_0\in Y, \\
\limsup_{\e\rightarrow 0} h_{\e}(y_0)\leq  k^*(y_0) \ \ \ \forall \ y_0\in Y.\\
\end{aligned}
\end{equation}
\end{Proposition}
{\bf Proof.} Let us prove the first inequality in (\ref{est-from-above}). Due to the fact that the function $v_T(\cdot)$ is lower  semicontinuous, there
exists a sequence of continuous functions  $v_{l,T}(\cdot), \ l=1,2,...,$ such that
\begin{equation}\label{est-from-above-1}
v_{l,T}(y)\leq v_{l+1,T}(y)\leq v_{T}(y)\ \ \forall \ y\in Y,\ \ \forall\ l=1,2,..., \ \  {\rm and } \ \  \lim_{l\rightarrow\infty}v_{l,T}(y) = v_{T}(y)\ \ \forall \ y\in Y
\end{equation}
and, for any $l=1,2,... $,
\begin{equation}\label{est-from-above-2}
\max_{y\in Y}|v_{l,T}(y)|\leq \sup_{y\in Y}|v_{T}(y)| \leq \max_{(y,u)\in Y\times U}|k(y,u)|\BYDEF M  .
\end{equation}
(see, e.g., Theorem A6.6 in \cite{Ash}). Let
\begin{equation}\label{e-kappa}
\k_{l,T}(y)\BYDEF v_{T}(y) - v_{l,T}(y).
\end{equation}
Note that, by  (\ref{est-from-above-1}) and (\ref{est-from-above-2}),
\begin{equation}\label{e-kappa-1}
0\leq \k_{l+1,T}(y)\leq \k_{l,T}(y)\leq 2 M\ \ \forall \ y\in Y,\ l=1,2,..., \ \  {\rm and } \ \  \lim_{l\rightarrow\infty}\k_{l,T}(y) =0 \ \ \forall \ y\in Y.
\end{equation}
Consider the following IDLP problem
\begin{equation*}\label{BB21-1}
\sup_{(\psi,\eta)\in Q_l(T)} \psi(y_0)\BYDEF  d_l^*(T,y_0),
\end{equation*}
where
$Q_l(T)$ is the set of pairs $(\psi(\cdot),\eta(\cdot))\in C(Y)\times C(Y)$ that satisfy the inequalities
\begin{equation}\label{BB25-1}
\begin{aligned}
&k(y,u)-\psi(y)+E[\eta(f(y,u,s))]-\eta(y)\ge 0,\\
&E[\psi(f(y,u,s))]-\psi(y)\ge -\frac{2M}{T} - E[\k_{l,T}(f(y,u,s))]\ \ \ \ \forall \ (y,u)\in G.
\end{aligned}
\end{equation}
Let us show that, for an arbitrary small  $\beta > 0 $ and for every $l=1,2,...$, there exists a function $\eta_{l, T,\beta}(\cdot)\in C(Y) $ such that
\begin{equation}\label{BB21-1-1}
\left(\psi_{l,T,\beta}(\cdot), \eta_{l,T,\beta}(\cdot)\right)\in Q_l(T), \ \ \ {\rm where} \ \ \psi_{l,T,\beta}(\cdot)\BYDEF v_{l,T}(\cdot)-\beta.
\end{equation}
Note that, if the inclusion above is established,  it would imply that
\begin{equation}\label{BB24-1}
v_{l,T}(y_0)-\beta\le  d^*_l(T, y_0).
\end{equation}
Let us  first verify that there exists $\eta_{l,T,\beta}(\cdot)\in C(Y) $ such that the pair $(\psi_{l,T,\beta}(\cdot), \eta_{l,T,\beta}(\cdot))$ satisfies the first inequality in (\ref{BB25-1}). To this end, note
that the inequality (\ref{e-before-lim-1}) is equivalent to the inequality
\begin{equation*}
\int_{G}(k(y,u)-v_{T}(y))\,\gm\ge 0\quad \hbox{for all }\g\in W,
\end{equation*}
which, due to the fact that  $v_{l,T}(y) \le v_{T}(y)$ (see (\ref{est-from-above-1})), implies that
\begin{equation*}
\int_{G}(k(y,u)-v_{l,T}(y))\,\gm\ge 0\quad \hbox{for all }\g\in W .
\end{equation*}
The latter, in turn, is equivalent to
\begin{equation}\label{CC15-1}
\min_{\g\in W}\int_{G}(k(y,u)-v_{l,T}(y))\,\gm\ge 0.
\end{equation}
The problem on the left hand side of (\ref{CC15-1}), i.e.,
\begin{equation}\label{CC8-1}
\min_{\g\in W}\int_{G}(k(y,u)-v_{l,T}(y))\,\gm,
\end{equation}
is an IDLP problem, which is similar to (\ref{M22}) (with $k(y,u)-v_{l,T}(y)$ instead of $k(y,u)$ in the objective function).
The problem dual to (\ref{CC8-1}) is of the form (compare with (\ref{M21}))
\begin{equation}\label{BB12-1}
\sup_{\eta\in C(Y)}\inf_{(y,u)\in G}\{k(y,u)-v_{l,T}(y)+E[\eta(f(y,u,s))]-\eta(y)\}.
\end{equation}
By Proposition \ref{SD-1}, the optimal values of \eqref{CC8-1} and \eqref{BB12-1} are equal. Therefore, \eqref{CC15-1} is equivalent to
\begin{equation}\label{BB15-1}
\sup_{\eta\in C(Y)}\inf_{(y,u)\in G}\{k(y,u)-v_{l,T}(y)+E[\eta(f(y,u,s))]-\eta(y)\}\ge 0.
\end{equation}
From \eqref{BB15-1} it follows that, for any $\beta>0$, there exists a function $\eta_{l, T,\beta}(\cdot)\in C(Y)$ such that
\begin{equation*}\label{e-eps-feasib-1}
k(y,u)-v_{l,T}(y)+E[\eta_{l, T,\beta}(f(y,u,s))]-\eta_{l, T,\beta}(y)\ge -\beta\quad \hbox{for all }(y,u)\in G.
\end{equation*}
The latter implies that  the pair $(\psi_{l, T,\beta}(\cdot), \eta_{l, T,\beta}(\cdot))$, where $\psi_{l, T,\beta}(\cdot):= v_{l,T}(\cdot)-\beta $, satisfies the first inequality in (\ref{BB25-1}).
Let us now  verify that the function $\psi_{l, T,\beta}(\cdot)= v_{l,T}(\cdot)-\beta$ satisfies the second inequality in (\ref{BB25-1}).
In accordance with the dynamic programming principle (see (\ref{e-App-1})),  for any $T\geq 1$,
\begin{equation}\label{e-before-lim-1-1}
Tv_{T}(y)\le k(y,u)+(T-1)E[v_{T-1}(f(y,u,s))] \ \ \ \forall\ (y,u)\in G.
\end{equation}
Also, as can be readily seen,
\begin{equation}\label{e-before-lim-1-1-1}
(T-1)E[v_{T-1}(f(y,u,s))]\le T E[v_{T}(f(y,u,s))] + M \ \ \ \forall\ (y,u)\in G,
\end{equation}
where $M$ is as in (\ref{est-from-above-2}).
By (\ref{e-before-lim-1-1}) and (\ref{e-before-lim-1-1-1}),
$$
Tv_{T}(y)\le k(y,u) +  T E[v_{T}(f(y,u,s))] + M \leq  T E[v_{T}(f(y,u,s))] +  2M.
$$
Consequently (having in mind (\ref{e-kappa})),
$$
v_{l,T}(y) \leq v_{T}(y)\leq E[v_{l,T}(f(y,u,s))] +  \frac{2M}{T}+ E[\kappa_{l, T}(f(y,u,s))]  \ \ \
$$
$$
 \Rightarrow \ \ \ \psi_{l, T,\beta}(y) \leq E[\psi_{l, T,\beta}(f(y,u,s))] +  \frac{2M}{T} + E[\kappa_{l, T}(f(y,u,s))].
$$
Thus, $\psi_{l, T,\beta}(\cdot)= v_{l,T}(\cdot)-\beta$ satisfies the second inequality in (\ref{BB25-1}). Hence,
(\ref{BB21-1-1}) is valid and, consequently, (\ref{BB24-1}) is valid too. Moreover, the latter  implies that
\begin{equation}\label{BB24-1-1}
v_{l,T}(y_0)\le  d^*_l(T, y_0)
\end{equation}
since $\beta >0 $ in (\ref{BB24-1}) is arbitrary small.

By Lemma \ref{L-equivalince},
\begin{equation}\label{e-inequality-lemma-1}
 d^*_l(T, y_0)\leq k^*_l(T, y_0),
\end{equation}
where
\begin{equation*}\label{BB1-T}
 k^*_l(T, y_0)\BYDEF\inf_{(\g,\xi)\in \O(y_0)}\left\{ \int_G k(y,u)\gm  + \frac{2M}{T}\int_G\xi(dy,du) + \int_GE[\k_{l,T}(f(y,u,s))]\xi(dy,du)\right\}.
\end{equation*}
(Note that, to adjust the notations used above and the notations used in Lemma \ref{L-equivalince}, one should write $ d^*_l(T, y_0)$ and $ k^*_l(T, y_0)$ as $d^*_{\theta_{T,l}}(y_0) $ and $ k^*_{\theta_{T,l}}(y_0)$, where $\theta_{T,l}(y,u)= \frac{2M}{T} + E[\k_{l,T}(f(y,u,s))] $.)

From (\ref{BB24-1-1}) and (\ref{e-inequality-lemma-1}) it follows that
\begin{equation}\label{BB1-T-1}
v_{l,T}(y_0)\le k^*_l(T, y_0).
\end{equation}
As can be readily seen, $k^*_l(T, y_0) $ is monotone decreasing in $l$ (due to (\ref{e-kappa-1})). Therefore, there exists a limit $\lim_{l\rightarrow\infty}k^*_l(T, y_0) $. Let us show that
\begin{equation}\label{BB1-T-2-1}
\lim_{l\rightarrow\infty}k^*_l(T, y_0) = k^*(T, y_0),
\end{equation}
where
\begin{equation*}\label{BB1-T-3}
 k^*(T, y_0)\BYDEF\inf_{(\g,\xi)\in \O(y_0)}\left\{ \int_G k(y,u)\gm  + \frac{2M}{T}\int_G\xi(dy,du) \right\}.
\end{equation*}
Firstly, note that $\ \lim_{l\rightarrow\infty}k^*_l(T, y_0)\geq k^*(T, y_0) $ (since $\ k^*_l(T, y_0)\geq k^*(T, y_0) $ for any $l=1,2,...$). To show the validity of the opposite inequality, take an arbitrary small $\beta > 0 $ and choose $(\gamma', \xi')\in \Omega(y_0) $ such that
$$
 \int_G k(y,u)\gamma'(dy,du)  + \frac{2M}{T}\int_G\xi'(dy,du)\leq k^*(T, y_0) + \beta.
$$
Then
$$
k^*_l(y_0, T)\leq \int_{Y\times U} k(y,u)\gamma'(dy,du)+\frac{2M}{T}\int_{Y\times U} \xi'(dy,du)+ \int_GE[\k_{l,T}(f(y,u,s))]\xi'(dy,du)
$$
$$
\leq  k^*(y_0, T) + \beta + \int_GE[\k_{l,T}(f(y,u,s))]\xi'(dy,du).
$$
Since $\ \lim_{l\rightarrow \infty}\int_GE[\k_{l,T}(f(y,u,s))]\xi'(dy,du) = 0 $ (by (\ref{e-kappa-1}) and the Monotone Convergence Theorem; see, e.g., Theorem 1.6.2 in \cite{Ash}), it follows that
$$
 \lim_{l\rightarrow \infty}k^*_l(y_0, T)\leq k^*(y_0, T) + \beta \ \ \ \Rightarrow \ \ \ \lim_{l\rightarrow \infty}k^*_l(y_0, T)\leq k^*(y_0, T)
$$
(the latter being due to the fact that $\beta$ is arbitrary small).
Thus, (\ref{BB1-T-2-1}) is valid, and (along with   (\ref{est-from-above-1}) and (\ref{BB1-T-1}))  it implies
 that
\begin{equation}\label{BB1-T-1-1}
v_T(y_0) \le k^*(T, y_0).
\end{equation}
The function $k^*(T, y_0)$ is monotone decreasing in $T$ and $k^*(T, y_0)\geq k^*(y_0) $ for any $T$. Therefore, there exists a limit
$\ \lim_{T\rightarrow\infty}k^*(T, y_0) \geq k^*(y_0) $. In fact, arguing as above, one can establish that
$$
\lim_{T\rightarrow\infty}k^*(T, y_0) = k^*(y_0).
$$
This and (\ref{BB1-T-1-1}) establish the validity of the first inequality  in (\ref{est-from-above}).

Let us now prove the second inequality in (\ref{est-from-above}) (this proof being very similar to that of the first one).
Since that the function $h_{\e}(\cdot)$ is lower  semicontinuous, there
exists a sequence of continuous functions  $h_{l,\e}(\cdot), \ l=1,2,...,$ such that
\begin{equation}\label{est-from-above-1-h}
h_{l,\e}(y)\leq h_{l+1,\e}(y)\leq h_{\e}(y)\ \ \forall \ y\in Y,\ \ \forall\ l=1,2,..., \ \  {\rm and } \ \  \lim_{l\rightarrow \infty}h_{l,\e}(y) = h_{\e}(y)\ \ \forall \ y\in Y
\end{equation}
and, for any $l=1,2,... $,
\begin{equation}\label{est-from-above-2-h}
\max_{y\in Y}|h_{l,\e}(y)|\leq \sup_{y\in Y}|h_{\e}(y)| \leq \max_{y\in Y}|k(y,u)|= M  .
\end{equation}
(see Theorem A6.6 in \cite{Ash}). Let
\begin{equation*}\label{e-kappa_a}
\k_{l,\e}(y)\BYDEF h_{\e}(y) - h_{l,\e}(y).
\end{equation*}
Note that, by  (\ref{est-from-above-1-h}) and (\ref{est-from-above-2-h}),
\begin{equation*}\label{e-kappa-1-h}
0\leq \k_{l+1,\e}(y)\leq \k_{l,\e}(y)\leq 2 M\ \ \forall \ y\in Y,\ l=1,2,..., \ \  {\rm and } \ \  \lim_{l\rightarrow \infty}\k_{l,\e}(y) =0 \ \ \forall \ y\in Y.
\end{equation*}
Consider the IDLP problem
\begin{equation*}\label{BB21-1-2}
\sup_{(\psi,\eta)\in Q_l(\e)} \psi(y_0)\BYDEF d^*_l(\e,y_0),
\end{equation*}
where
$Q_l(\e)$ is the set of pairs $(\psi(\cdot),\eta(\cdot))\in C(Y)\times C(Y)$ that satisfy the inequalities
\begin{equation}\label{BB25-1-2}
\begin{aligned}
&k(y,u)-\psi(y)+E[\eta(f(y,u,s))]-\eta(y)\ge 0,\\
&E[\psi(f(y,u,s))]-\psi(y)\ge -2M\e - E[\kappa_{l, \e}(f(y,u,s))]\ \ \ \ \forall \ (y,u)\in G.
\end{aligned}
\end{equation}
Let us show that, for an arbitrary small  $\beta > 0 $, there exists  a function $\eta_{l, \e,\beta}(\cdot)\in C(Y) $ such that
\begin{equation}\label{BB21-1-3}
\left(\psi_{l, \e,\beta}(\cdot), \eta_{l, \e,\beta}(\cdot)\right)\in Q_l(\e), \ \ \ {\rm where} \ \ \psi_{l, \e,\beta}(\cdot)\BYDEF h_{l, \e}(\cdot)-\beta,
\end{equation}
with the inclusion above  implying that
\begin{equation}\label{BB24-3}
h_{l, \e}(y_0)-\beta\le  d^*_l(\e, y_0).
\end{equation}
To verify (\ref{BB21-1-3}), let us first show
that there exists $\eta_{l, \e,\beta}(\cdot)\in C(Y) $ such that the pair $(\psi_{l, \e,\beta}(\cdot), \eta_{l, \e,\beta}(\cdot))$ satisfies the first inequality in (\ref{BB25-1-2}). To this end, let us  rewrite the inequality (\ref{e-before-lim-2}) in the form
\begin{equation*}
\int_{G}(k(y,u)-h_{\e}(y))\,\gm\ge 0 \ \ \ \ \forall \ \g\in W.
\end{equation*}
By (\ref{est-from-above-1-h}), the latter implies
\begin{equation*}
\int_{G}(k(y,u)-h_{l,\e}(y))\,\gm\ge 0 \ \ \ \ \forall \ \g\in W,
\end{equation*}
which is equivalent to
\begin{equation}\label{CC15-3}
\min_{\g\in W}\int_{G}(k(y,u)-h_{l,\e}(y))\,\gm\ge 0.
\end{equation}
The problem on the left hand side of (\ref{CC15-3}), i.e.,
\begin{equation}\label{CC8-3}
\min_{\g\in W}\int_{G}(k(y,u)-h_{l,\e}(y))\,\gm,
\end{equation}
is an IDLP problem, which is similar to (\ref{M22}) (with $k(y,u)-h_{l,\e}(y)$ instead of $k(y,u)$ in the objective function).
The problem dual to (\ref{CC8-3}) is of the form (compare with (\ref{M21}))
\begin{equation}\label{BB12-3}
\sup_{\eta\in C(Y)}\inf_{(y,u)\in G}\{k(y,u)-h_{l,\e}(y)+E[\eta(f(y,u,s))]-\eta(y)\}.
\end{equation}
By Proposition \ref{SD-1}, the optimal values of \eqref{CC8-3} and \eqref{BB12-3} are equal. Therefore, \eqref{CC15-3} is equivalent to
\begin{equation}\label{BB15-1-3}
\sup_{\eta\in C(Y)}\inf_{(y,u)\in G}\{k(y,u)-h_{l,\e}(y)+E[\eta(f(y,u,s))]-\eta(y)\}\ge 0.
\end{equation}
From \eqref{BB15-1-3} it follows that, for any $\beta>0$, there exists a function $\eta_{l,\e,\beta}(\cdot)\in C(Y)$ such that
\begin{equation*}\label{e-eps-feasib-1-3}
k(y,u)-h_{l,\e}(y) +E[\eta_{l,\e,\beta}(f(y,u,s))]-\eta_{l,\e,\beta}(y)\ge -\beta\ \ \ \ \forall \ (y,u)\in G.
\end{equation*}
The latter implies that  the pair $(\psi_{l,\e,\beta}(\cdot), \eta_{l,\e,\beta}(\cdot))$, where $\psi_{l,\e,\beta}(\cdot):= h_{l,\e}(\cdot)-\beta $, satisfies the first inequality in (\ref{BB25-1-2}).

To  verify that the function $\psi_{l,\e,\beta}(\cdot)= h_{l,\e}(\cdot)-\beta$ satisfies the second inequality in (\ref{BB25-1-2}), note that
by  (\ref{e-App-2}) of Proposition \ref{Prop-lsc},
\begin{equation*}\label{e-before-lim-2-1}
h_{\e}(y)\le \e k(y,u)+(1-\e) E[h_{\e}(f(y,u,s))]\ \ \ \forall \ (y,u)\in G .
\end{equation*}
The latter implies that
$$
h_{\e}(y)\le E[h_{\e}(f(y,u,s))] + \e(k(y,u) - E[h_{\e}(f(y,u,s))]) \ \ \ \forall \ (y,u)\in G.
$$
This, in turn, leads to
\begin{equation*}\label{e-before-lim-2-2}
h_{l, \e}(y) \le h_{\e}(y)\le  E[h_{l,\e}(f(y,u,s)))]+ 2M\e + E[\kappa_{l,\e}(f(y,u,s))] \ \ \ \forall \ (y,u)\in G
\end{equation*}
(since, as can be readily seen, $\ \max_{y\in Y}|h_{\e}(y)|\leq M $).
Thus, $\psi_{l,\e,\beta}(\cdot)= h_{l,\e}(\cdot)-\beta$ satisfies the second inequality in (\ref{BB25-1-2}), and, therefore,
(\ref{BB24-3}) is valid, the latter implying that
\begin{equation}\label{BB24-1-1-h}
h_{l, \e}(y_0)\le  d^*_l(\e, y_0)
\end{equation}
(since $\beta >0 $ in (\ref{BB24-3}) is arbitrary small).

By Lemma \ref{L-equivalince},
\begin{equation}\label{e-inequality-lemma-1-h}
 d^*_l(\e, y_0)\leq k^*_l(\e, y_0),
\end{equation}
where
\begin{equation*}\label{BB1-T-h}
 k^*_l(\e, y_0)\BYDEF\inf_{(\g,\xi)\in \O(y_0)}\left\{ \int_G k(y,u)\gm  + 2M\e\int_G\xi(dy,du) + \int_GE[\k_{l,\e}(f(y,u,s))]\xi(dy,du)\right\}.
\end{equation*}
(In this case, to adjust the notations used above and the ones used in Lemma \ref{L-equivalince}, one should write $ d^*_l(\e, y_0)$ and $ k^*_l(\e, y_0)$ as $d^*_{\theta_{\e,l}}(y_0) $ and $ k^*_{\theta_{\e,l}}(y_0)$, where $\theta_{\e,l}(y,u)= 2M\e + E[\k_{l,\e}(f(y,u,s))] $.)

From (\ref{BB24-1-1-h}) and (\ref{e-inequality-lemma-1-h}) it follows that
\begin{equation}\label{BB1-T-1-h}
h_{l,\e}(y_0)\le k^*_l(\e, y_0).
\end{equation}
Using the argument similar to one used above, we can show that
\begin{equation}\label{BB1-T-2-1-h}
\lim_{l\rightarrow 0}k^*_l(\e, y_0) = k^*(\e, y_0)
\end{equation}
where
\begin{equation*}\label{BB1-T-3-h}
 k^*(\e, y_0)\BYDEF\inf_{(\g,\xi)\in \O(y_0)}\left\{ \int_G k(y,u)\gm  + 2M\e\int_G\xi(dy,du) \right\}.
\end{equation*}
Subsequently, one can show that
\begin{equation}\label{BB1-T-2-1-h-1}
\lim_{\e\rightarrow 0}k^*(\e, y_0) = k^*(y_0).
\end{equation}
The validity of the second inequality in (\ref{est-from-above}) follows from (\ref{BB1-T-1-h}), (\ref{BB1-T-2-1-h}),  (\ref{BB1-T-2-1-h-1}) (and from (\ref{est-from-above-1-h})).
$\ \Box$

\medskip

{\bf Proof of Theorem \ref{ThN1}.} If the point-wise limit (\ref{chesaro-lim-1-1}) exists, then, by Theorem \ref{Th-lower-upper-limit-main},
the limit function $v(\cdot) $  satisfies the inequality
$$
v(y_0)\geq d^*(y_0) \ \ \ \forall \ y_0\in Y.
$$
Therefore, to prove the statement (a), one needs to   show that
\begin{equation}\label{important-ineq-1}
v(y_0)\leq d^*(y_0) \ \ \forall \ y_0\in Y.
\end{equation}
Similarly,  if  the point-wise limit (\ref{abel-lim-3-1}) exists, then, by Theorem \ref{Th-lower-upper-limit-main},
the limit function
$h(\cdot)$ satisfies  the inequality
$$
h(y_0)\geq d^*(y_0)\ \ \ \forall \ y_0\in Y.
$$
Therefore, to prove the statement (b), one needs to   show that
\begin{equation}\label{important-ineq-2}
h(y_0)\leq d^*(y_0)\ \ \ \forall \ y_0\in Y.
\end{equation}
We will prove only (\ref{important-ineq-1}) (the proof of (\ref{important-ineq-2}) follows exactly the same lines). Note that, from the dynamic programming principle (\ref{e-App-1}), it follows that
$$
Tv_{T}(y)\le k(y,u)+(T-1)E[v_{T-1}(f(y,u,s))] \ \ \ \forall\ (y,u)\in G
$$
  for any $T\geq 1$.
Note also that by dividing the latter by $T$ and passing to the limit as $T\to \infty$, one obtains
\begin{equation}\label{V51-1}
v(y)\le E[v(f(y,u,s))]  \ \ \forall \ (y,u)\in G.
\end{equation}
Also, by passing to the limit as $T\to \infty$ in (\ref{e-before-lim-1}), one obtains
\begin{equation}\label{V5-1}
\int_{G}v(y)\,\g(dy,du)\le \int_{G}k(y,u)\,\g(dy,du) \ \ \forall \ \g\in W.
\end{equation}
Inequality (\ref{V5-1})  can be rewritten in the form
\begin{equation*}
\int_{G}(k(y,u)-v(y))\,\gm\ge 0\quad \hbox{for all }\g\in W,
\end{equation*}
which is equivalent to that
\begin{equation}\label{CC15}
\min_{\g\in W}\int_{G}(k(y,u)-v(y))\,\gm\ge 0.
\end{equation}
The problem in the left hand side of the above inequality,
\begin{equation}\label{CC8}
\min_{\g\in W}\int_{G}(k(y,u)-v(y))\,\gm,
\end{equation}
is an IDLP problem, whose dual is
\begin{equation}\label{BB12}
\sup_{\eta\in C(Y)}\inf_{(y,u)\in G}\{k(y,u)-v(y)+E[\eta(f(y,u,s))]-\eta(y)\}.
\end{equation}
By equation (\ref{e-SE-2}) of Proposition \ref{SD-1} (considered with $k(y,u)-v(y)$ instead of $k(y,u)$), the optimal values of \eqref{CC8} and \eqref{BB12} are equal. Therefore, \eqref{CC15} is equivalent to
\begin{equation}\label{BB15}
\sup_{\eta\in C(Y)}\inf_{(y,u)\in G}\{k(y,u)-v(y)+E[\eta(f(y,u,s))]-\eta(y)\}\ge 0.
\end{equation}
From \eqref{BB15} it follows that, for any $\beta>0$, there exists a function $\eta_{\beta}(\cdot)\in C(Y)$ such that
\begin{equation}\label{e-eps-feasib}
k(y,u)-v(y)+E[\eta_{\beta}(f(y,u,s))]-\eta_{\beta}(y)\ge -\beta\quad \hbox{for all }(y,u)\in G.
\end{equation}
Consider now the problem
\begin{equation}\label{BB21}
\sup_{(\psi,\eta)\in Q} \psi(y_0)= d^*(y_0),
\end{equation}
where
$Q$ is the set of pairs $(\psi,\eta)\in C(Y)\times C(Y)$ that satisfy the inequalities
\begin{equation}\label{BB25}
\begin{aligned}
&k(y,u)-\psi(y)+E[\eta(f(y,u,s))]-\eta(y)\ge 0,\\
&E[\psi(f(y,u,s))]-\psi(y)\ge 0\quad \hbox{for all }(y,u)\in G.
\end{aligned}
\end{equation}
Note that the optimal value of problem (\ref{BB21}) is the same as that of (\ref{BB8}) (see (\ref{e-Pert-1}) taken with $\theta = 0$). Due to (\ref{V51-1}) and (\ref{e-eps-feasib}),  the pair $(\psi_{\beta}(\cdot), \eta_{\beta}(\cdot)) $, where  $\psi_{\beta}(\cdot):=v(\cdot)-\beta$, satisfies the inequalities \eqref{BB25}. Consequently,
$$
d^*(y_0)\ge v(y_0) - \beta \ \ \ \forall \ y_0 \in Y.
$$
This proves (\ref{important-ineq-1}), since $\beta > 0$ is arbitrarily small. Thus, {\rm (a)} is proved.
 $\ \Box$

\section{Proof of Theorem \ref{Th-discounting-1}}\label{Proofs-Secondary}

Let us prove (\ref{convergence-to-W-dis}). To this end, let us first demonstrate that the following inclusion holds:
	\begin{equation}\label{P11}
	\bar{\rm co} \T_\e(y_0) \subset W(\e, y_0)\ \ \ \ \forall \e\in (0,1).
	\end{equation}
	Let $\g\in \T_\e(y_0) $. That is, $\g$ is the discounted occupational measure generated  by a control plan $\pi\in \Pi$: $\gamma = \g_d^{\pi,y_0,\e}$. The state-control trajectory $(y^{\pi, y_0}(\cdot), u^{\pi,  y_0}(\cdot))$ obtained with the use of this plan  satisfies the equality
\begin{equation*}
	E\left[\sum_{t=0}^{\infty}(1-\e)^t \psi(y^{\pi,  y_0}(t))  \right] = \psi(y_0) + E\left[(1-\e)\sum_{t=0}^{\infty}(1-\e)^t\psi(f(y^{\pi,  y_0}(t),u^{\pi,  y_0}(t),s(t)))\right]
	\end{equation*}
for any $\psi \in C(Y)$. Multiplying both sides of this equality by $\e$ and using  (\ref{G8}), we obtain
	\begin{equation*}
	\begin{split}
	\int_G \psi(y) \g_d^{\pi,y_0,\e}(dy,du) =& \int_G \e\psi(y_0)\g_d^{\pi,y_0,\e}(dy,du)\\
	&+ \int_G (1-\e) E\big[\psi(f(y,u,s))\big]\g_d^{\pi,y_0,\e}(dy,du).
	\end{split}
	\end{equation*}
	Rearranging the above, we have
	\begin{equation*}
	\int_G \big((1-\e) (E\big[\psi(f(y,u,s))\big]-\psi(y)) + \e(\psi(y_0)-\psi(y)) \big)\g_d^{\pi,y_0,\e}(dy,du) = 0,
	\end{equation*}
	implying $\g_d^{\pi,y_0,\e} \in W(\e, y_0)$. Hence, $\ \T_\e(y_0)\subset  W(\e, y_0) $. The validity of (\ref{P11})  follows from the fact that
$W(\e, y_0)$ is convex and compact.

 To prove the converse inclusion, it is sufficient to prove that the inequality
\begin{equation}\label{e-any-k}
	\min_{\g \in \bar{\rm co} \T_\e(y_0)}\int_G k(y,u)\g(dy,du) \leq \min_{\g \in W(\e, y_0)} \int_G k(y,u)\g(dy,du) \ \ \ \forall \ \e\in (0,1)
	\end{equation}
is valid for an arbitrary continuous function $k(y,u)  $. (The fact that the validity of  (\ref{e-any-k})  for any $k(y,u)  $ implies the inclusion converse to (\ref{P11})
follows from the  separation theorem;  see the proof of Proposition \ref{Propo-separation-consec} in Section  \ref{Sec-Appendix}).
 Due to (\ref{A3}) and (\ref{auxiliary-1}), the latter is equivalent to
\begin{equation}\label{e-convers-ineq}
h_{\e}(y_0)\leq k^*(\e,y_0).
\end{equation}
Let $LS$ stand for the set of bounded  lower semicontinuous functions on $G$, and, for any $\psi \in LS$, let
	\begin{equation}
	\label{eq:DGA_1}
	\mu(\psi, \e, y_0) \BYDEF \inf_{(y,u) \in G} \big\{k(y,u)+(1-\e) \big(E\big[\psi(f(y,u,s))\big]-\psi(y)\big) + \e\big(\psi(y_0)-\psi(y)\big) \big\}.
	\end{equation}
Obviously (compare with (\ref{e-dis-dual})),
\begin{equation}\label{e-dis-dual-LS}
\mu_{LS}^*( \e, y_0)\BYDEF\sup_{\psi\in LS}\mu(\psi, \e, y_0)\geq \mu^*( \e, y_0).
\end{equation}
By (\ref{e-App-2}) (see Proposition \ref{Prop-lsc}), we have
	\begin{equation*}
	\min_{u \in U(y)}\left\{k(y,u)+(1-\e) E [V_\e(f(y,u,s))] - V_\e(y)  \right\} = 0 \quad \forall y \in Y,
	\end{equation*}
	which implies
	\begin{equation*}
	\min_{(y,u) \in G}\left\{k(y,u)+(1-\e) E [V_\e(f(y,u,s))] - V_\e(y)  \right\} = 0.
	\end{equation*}
	Hence,
	\begin{equation}
	\label{eq:V_a_solves_mu_a}
	\begin{split}
	h_\e(y_0) = \e V_\e(y_0) = \min_{(y,u) \in G}\{k(y,u) + (1-\e) \big(E [V_\e(f(y,u,s))]-V_\e(y)\big)\\+\e\big(V_\e(y_0)-V_\e(y)\big)  \}= \mu(V_{\e}, \e, y_0)
\le \mu_{LS}^*( \e, y_0).
	\end{split}
	\end{equation}
Take  an arbitrary $\g \in W(\e, y_0)$ and an arbitrary $\psi \in LS$. Let $\{\psi_n\}_{n=1}^\infty$ be a bounded sequence of continuous functions such that $\psi_n(y) \rightarrow \psi(y)$ point-wise on $Y$ as $n \rightarrow \infty$ (such a sequence exists; see, e.g., \cite[Theorem A6.6]{Ash}).
	From Lebesgue's dominated convergence theorem (see, e.g., Theorem 1.6.9, p. 49 in \cite{Ash}) and from the definition of $W(\e, y_0)$, it follows that
	\begin{equation*}
	\begin{split}
	\mu (\psi,\e, y_0) &\le \int_G \big(k(y,u)+(1-\e) (E [\psi(f(y,u,s))]-\psi(y)) + \e(\psi(y_0)-\psi(y))\big)\g(dy,du) \\
	&= \lim_{n \rightarrow \infty} \int_G \big(k(y,u)+(1-\e) (E [\psi_n(f(y,u,s))]-\psi_n(y)) + \e(\psi_n(y_0)-\psi_n(y))\big)\g(dy,du)\\
	&= \int_G k(y,u) \g(dy,du).
	\end{split}
	\end{equation*}
	Taking $sup$ with respect to $\psi \in LS$ in the left-hand-side and then taking $min$ with respect to $\g \in W(\e, y_0)$ in the right-hand-side allow us to conclude that
	\begin{equation}\label{e-conv-2}
	\mu_{LS}^*( \e, y_0) \le k^*(\e,y_0).
	\end{equation}
The latter and (\ref{eq:V_a_solves_mu_a}) prove (\ref{e-convers-ineq}).
{Note that (\ref{e-conv-2}) along with (\ref{e-SE-1}),  (\ref{e-dis-dual-LS}) imply also that $\mu_{LS}^*( \e, y_0)=\mu^*( \e, y_0)$.}
Thus, the validity of (\ref{convergence-to-W-dis}) is established.

Let us  prove (\ref{convergence-to-W-dis-1}). It is straightforward to verify that
$$
\limsup_{\e \rightarrow 0}\left\{\bigcup_{y_0 \in Y} W(\e, y_0)\right\} \subset W \ \ \ \ \Rightarrow \ \ \ \ \limsup_{\e \rightarrow 0}\bar{\rm co}\left\{\bigcup_{y_0 \in Y} W(\e, y_0)\right\} \subset W ,
$$
the second inclusion being due to the fact that $W $ is convex and compact.
By (\ref{convergence-to-W-dis}), the latter is equivalent to
\begin{equation}
\label{e-Theta-1}
\limsup_{\e \rightarrow 0}\bar{\rm co} \T_{\e} \subset W .
\end{equation}
Let
	\begin{equation*}
	\label{eq:DGA-a-1}
	\mu(\psi) \BYDEF \inf_{(y,u) \in G} \big\{k(y,u)+ E\big[\psi(f(y,u,s))\big]-\psi(y)  \big\}
	\end{equation*}
and let
\begin{equation}\label{e-dis-dual-LS-1}
\mu_{LS}^*\BYDEF\sup_{\psi\in LS}\mu(\psi)\geq \mu^*
\end{equation}
(see (\ref{M21})).
By (\ref{e-App-2}),
\begin{equation*}
 k(y,u) + (1-\e) E [V_\e(f(y,u,s))]-V_\e(y)\ge 0 \ \ \ \forall \ (y,u) \in G
\end{equation*}
for any $\e \in (0,1)$. After rearranging, the latter leads to
\begin{equation*}
\begin{split}
k(y,u) + E [V_\e(f(y,u,s))] - V_\e(y) &\ge \e E[V_\e(f(y,u,s))]
\ge \min_{y' \in Y}\e V_\e(y')
= \min_{y' \in Y} h_\e(y')\ \ \ \forall\ (y,u) \in G.
\end{split}
\end{equation*}
Consequently, since $V_\e(\cdot)$ is lower semi-continuous (as established by Proposition \ref{Prop-lsc}),
\begin{equation}\label{e-LS-a}
\mu_{LS}^* \ge \min_{y \in Y} h_\e(y)\ \ \ \forall \ \e\in (0,1).
\end{equation}

Take an arbitrary $\g \in W$ and an arbitrary $\psi \in LS$. There exists a sequence of bounded continuous functions $\{\psi_n\}_{n=1}^\infty$ that converges to $\psi(y)$ point-wisely on $Y$ as $n \rightarrow \infty$.
	From Lebesgue's dominated convergence theorem and from the definition of $W$, it follows that
	\begin{equation*}
	\begin{split}
	\mu (\psi) &\le \int_G \big(k(y,u)+E [\psi(f(y,u,s))]-\psi(y) \big)\g(dy,du) \\
	&= \lim_{n \rightarrow \infty} \int_G \big(k(y,u)+E [\psi_n(f(y,u,s))]-\psi_n(y)\big)\g(dy,du)= \int_G k(y,u) \g(dy,du).
	\end{split}
	\end{equation*}
	Taking $sup$ with respect to $\psi \in LS$ in the left-hand-side and then taking $min$ with respect to $\g \in W$ in the right-hand-side lead to the inequality
	$$
	\mu_{LS}^* \le k^*.
$$
	This and (\ref{e-LS-a}) lead  to the inequality
  $$
  \min_{y \in Y} h_\e(y)\le k^* \ \ \ \forall \ \e\in (0,1),
  $$
which can be rewritten as follows (see (\ref{M22}) and (\ref{auxiliary-2}))
$$
\min_{\g \in \bar{\rm co} \T_\e}\int_G k(y,u)\g(dy,du)\leq \min_{\g\in W} \int_{G} k(y,u)\g(dy,du) \ \ \ \forall \ \e\in (0,1)
$$
The latter is valid for an arbitrary continuous $k(y,u) $. Therefore, by the separation theorem,
\begin{equation}\label{e-W-W-dis-un}
W\subset \bar{\rm co} \T_\e,
\end{equation}
which along with
(\ref{e-Theta-1}) prove (\ref{convergence-to-W-dis-1}).

Let us now establish the validity of (\ref{convergence-to-W-dis-2}). To this end, let us first show that
\begin{equation}
\label{e-GT_and_W_inclusion}
\limsup_{T \rightarrow \infty} \G_T \subset W.
\end{equation}
Take an arbitrary $\g\in \limsup_{T \rightarrow \infty}\G_T $. That is, there exist sequences  $y_0^i\in Y$, $\pi^i \in \Pi, \ i=1,2,...,$ and  $T_i$ (with $T_i \rightarrow \infty$ as $i \rightarrow \infty$) such that the corresponding sequence of occupational measures $\g^{\pi^i , y_0^i ,T_i} \in \G_{T_i}$ converges to $\g$:
$\lim_{i\rightarrow \infty}\g^{\pi^i , y_0^i ,T_i}= \g$. For any $\ph \in C(Y)$, let $\bar \ph(y,u)$ be defined in accordance with (\ref{e-phi-bar}), that is, $\bar \ph(y,u) =E[\ph(f(y,u,s))] $. (Note that $\bar \ph \in C(Y)$ too.)
It can be readily understood (see (\ref{G88})) that
\begin{equation}\label{e-helpful}
\begin{split}
\int_G E[\ph(f(y,u,s))]\g^{\pi^i , y_0^i ,T_i}_i(dy,du)& =  \int_G \bar \ph(y,u)\g^{\pi^i , y_0^i ,T_i}_i(dy,du)={1\o T_i}E \left[\sum_{t=0}^{T_i-1} \bar\ph\left(y^{\pi^i , y_0^i }(t),u^{\pi^i , y_0^i }(t)\right)\right]\\
 & = \frac{1}{T_i} E\left[ \sum_{t=0}^{T_i-1} E\left[\bar \ph\left(y^{\pi^i,  y_0^i}(t), u^{\pi^i,  y_0^i}(t)\right)\bigg|y^{\pi^i, y_0^i}(t) \right]\right]\\
& = \frac{1}{T_i} E\left[ \sum_{t=0}^{T_i-1} E\left[ \ph\left(f\left(y^{\pi^i,  y_0^i}(t), u^{\pi^i,  y_0^i}(t),s(t)\right)\right)\bigg|y^{\pi^i, y_0^i}(t) \right]\right].
\end{split}
\end{equation}
Therefore,
\begin{equation*}
\begin{split}
&\int_G (E[\ph(f(y,u,s))] - \ph(y))\g^{\pi^i , y_0^i ,T_i}_i(dy,du)\\ &= \frac{1}{T_i} E\left[  \sum_{t=0}^{T_i - 1} \left(E\left[\ph\left(f\left(y^{\pi^i, y_0^i}(t), u^{\pi^i,  y_0^i}(t),s(t)\right)\right) \bigg|y^{\pi^i, y_0^i}(t) \right] - \ph\left(y^{\pi^i, y_0^i}(t)\right)\right)\right]  \\
&= \frac{1}{T_i} E\left[ \sum_{t=0}^{T_i-1} \left(E\left[\ph\left(y^{\pi^i,  y_0^i}(t+1)\right)\bigg|y^{\pi^i, y_0^i}(t) \right]-\ph\left(y^{\pi^i, y_0^i}(t)\right)\right)\right] \\
&=  \frac{1}{T_i} E\left[ \sum_{t=0}^{T_i-1} \left(\ph\left(y^{\pi^i,  y_0^i}(t+1)\right)-\ph\left(y^{\pi^i, y_0^i}(t)\right)\right)\right] = \frac{1}{T_i}E\left[ \ph\left(y^{\pi^i, y_0^i}(T_i)\right)-\ph\left(y_0^i\right)\right]
.
\end{split}
\end{equation*}
Hence,
\begin{equation*}
\begin{split}
\int_G ( E[\ph(f(y,u,s))]- \ph(y))\g(dy,du) &= \lim_{i\rightarrow \infty} \int_G (E[\ph(f(y,u,s))] - \ph(y)) \g^{\pi^i , y_{0_i} ,T_i}_i(dy,du) = 0 \ \ \forall \ \ph\in C(Y),
\end{split}
\end{equation*}
and, consequently, $\g\in W$. This proves  (\ref{e-GT_and_W_inclusion}), which implies
 that
\begin{equation}\label{auxiliary-3-1}
\liminf_{T\rightarrow\infty}\min_{y_0\in Y}v_{T}(y_0) \geq k^*  .
\end{equation}
Let us now prove that
\begin{equation}\label{auxiliary-3-2}
\limsup_{T\rightarrow\infty}\min_{y_0\in Y}v_{T}(y_0) \leq k^*.
\end{equation}
Note that, if (\ref{auxiliary-3-2}) is proved, then together with (\ref{auxiliary-3-1}), it will imply the validity (\ref{convergence-to-W-dis-5}), which, in turn, will establish the validity of
(\ref{convergence-to-W-dis-2}) (see Proposition \ref{Propo-separation-consec}).

To prove (\ref{auxiliary-3-2}), take a sequence $\e_i, \ i=1,2,...,  $  with $\e_i \rightarrow 0$ as $i \rightarrow \infty$. By  (\ref{convergence-to-W-dis-4}) (that follows from (\ref{convergence-to-W-dis-1}), which has been already proved), there exist sequences of initial conditions $y_0^i$ and control plans $\pi^i$ with the corresponding state-control trajectories $(y^{\pi^i, y_0^i}(\cdot), u^{\pi^i, y_0^i}(\cdot))$, such that
\begin{equation*}
\e_iE\left[\sum_{t=0}^{\infty}(1-\e_i)^tk\left(y^{\pi^i, y_0^i}(t), u^{\pi^i,  y_0^i}(t)\right)\right] = k^* + \zeta_i , \ \ \ {\rm where} \ \ \ \lim_{i\rightarrow \infty}\zeta_i = 0.
\end{equation*}
 Let us use Lemma \ref{lem:infty_time_to_finite_time} with $\sigma = k^* + \zeta_i$, $\ \d = \sqrt{-\ln(1-\e_i)}$ and $\ g(t) = E[k(y^{\pi^i, y_0^i}(t), u^{\pi^i,  y_0^i}(t))]$.
 (Lemma  \ref{lem:infty_time_to_finite_time} and Lemma \ref{lem:finite_time_translated} used below are stated at the end of this section.)
 According to this lemma, there exist $T_i \ge c / \sqrt{-\ln(1-\e_i)}, \ i=1,2,...$ ($c$ being a positive constant) such that
\begin{equation*}
\frac{1}{T_i}E\left[\sum_{t=0}^{T_i-1}k\left(y^{\pi^i, y_0^i}(t), u^{\pi^i,  y_0^i}(t)\right)\right]  < k^* + \zeta_i + \sqrt{-\ln(1-\e_i)} + \frac{2M}{T_i}, \ \ i=1,2,...\ .
\end{equation*}
From the expression above it follows that $\liminf_{T \rightarrow \infty}\min_{y \in Y} v_T(y) \le k^*$ which, along with (\ref{auxiliary-3-1}), implies that
\begin{equation}
\label{e-liminf_vT_eq_k*}
\liminf_{T \rightarrow \infty}\min_{y \in Y} v_T(y) = k^*.
\end{equation}
In accordance with (\ref{e-liminf_vT_eq_k*}), there exist a sequence $T_i, \ i=1,2,...$ ($T_i\rightarrow\infty$ as $i\rightarrow\infty $) and sequences of initial conditions $y_0^i $ and control plans $\pi^i \in \Pi$ with the corresponding state-control trajectories $(y^{\pi^i, y_0^i }(\cdot),u^{\pi^i, y_0^i }(\cdot))$ such that
\begin{equation*}
\frac{1}{T_i}E\left[\sum_{t=0}^{T_i - 1}k\left(y^{\pi^i, y_0^i}(t), u^{\pi^i,  y_0^i}(t)\right)\right] = k^* + \zeta_i, \ \ \ {\rm where}\ \ \ \lim_{i\rightarrow \infty}\zeta_i = 0.
\end{equation*}
Let us use Lemma \ref{lem:finite_time_translated} with $\sigma = k^* + \zeta_i$, $\ \d = \frac{1}{T_i}$ and $\ g(t) = E[k(y^{\pi^i, y_0^i}(t), u^{\pi^i,  y_0^i}(t))]$.
  According to this lemma,  there exists $T_i^* \in \{0,1,\ldots, T_i-1\}$ such that
$$
\frac{1}{T}E\left[\sum_{t=0}^{T-1}k\left(y^{\pi^i, y_0^i}(T_i^*+t),u^{\pi^i, y_0^i}(T_i^*+t)\right)\right]  \le k^* + \zeta_i + \frac{1}{T_i}\ \ \ \forall \ T \in \{1,\ldots, T_i - T_i^*\},
$$
with $T_i - T_i^* \rightarrow \infty$ as $i \rightarrow \infty$. As can be readily understood,
\begin{equation*}
\begin{split}
\frac{1}{T}E\left[\sum_{t=0}^{T-1}k\left(y^{\pi^i, y_0^i}(T_i^*+t),u^{\pi^i, y_0^i}(T_i^*+t)\right)\right] & =  E\left[\frac{1}{T}E\left[\sum_{t=0}^{T-1}k\left(y^{\pi^i, y_0^i}(T_i^*+t),u^{\pi^i, y_0^i}(T_i^*+t)\right)\bigg|y^{\pi^i, y_0^i}(T_i) \right] \right]\\ & \geq E\left[v_T\left(y^{\pi^i, y_0^i}(T_i)\right) \right]\geq \min_{y \in Y} v_T(y).
\end{split}
\end{equation*}
Therefore,
\begin{equation*}
\min_{y \in Y} v_T(y) \le k^* + \zeta_i + \frac{1}{T_i}\ \ \ \forall \ T \in \{1,\ldots, T_i - T_i^*  \}.
\end{equation*}
The latter implies (\ref{auxiliary-3-2}), and, thus, the proof is completed.
 $\ \Box$

 \begin{Lemma}
	\label{lem:infty_time_to_finite_time}
	Let $\mathbb{Z}^*$ represent the set of non-negative integers and define the function $g:\mathbb{Z}^* \rightarrow \mathbb{R}$ such that $|g(t)| \le M$ for all $t$. Let $\e \in (0,1)$ and
	\begin{equation}
	\sigma \BYDEF \e\sum_{t=0}^{\infty}(1-\e)^t g(t).
	\end{equation}
	Then, for any $\d > 0$, there exists a positive integer $T \ge \left[\frac{\d}{(4M+4|\sigma|+\d)(-ln(1-\e))}\right]$ satisfying
	\begin{equation}
	\frac{1}{T}\sum_{t=0}^{T-1}g(t) < \sigma + \d + \frac{2M}{T}.
	\end{equation}
\end{Lemma}

\begin{Lemma}
	\label{lem:finite_time_translated}
	Let $\mathbb{Z}^*$ represent the set of non-negative integers and define the function $g:\mathbb{Z}^* \rightarrow \mathbb{R}$ such that $|g(t)| \le M$ for all $t \in \mathbb{Z}^*$. Let $T$ be an arbitrary positive integer and
	\begin{equation*}
	\sigma \BYDEF \frac{1}{T}\sum_{t = 0}^{T-1}g(t).
	\end{equation*}
	For any $\d > 0$, there exists $T^* \in \{0,1, \ldots, T-1 \}$ such that
	\begin{equation}
	\frac{1}{S}\sum_{t=0}^{S-1}g(T^* + t) \le \sigma + \d \quad \text{for all } S \in \{1,\ldots,T-T^*\}.
	\end{equation}
	Moreover,
	\begin{equation}
	l(T) \BYDEF T - T^* \rightarrow \infty \quad as \quad T \rightarrow \infty.
	\end{equation}
\end{Lemma}
The proofs of Lemmas \ref{lem:infty_time_to_finite_time} and \ref{lem:finite_time_translated} can be found in \cite{GPS-17}. Note that these lemmas are discrete
time versions of the continuous time results established in \cite{GruneSIAM98} and \cite{GruneJDE98}.

\section{Proof  of Proposition \ref{Propo-separation-consec}}\label{Sec-Appendix}
Let us prove (\ref{convergence-to-W-dis}) provided that (\ref{convergence-to-W-dis-3}) is valid, that is,
\begin{equation*}\label{eq-OPR-0-0}
\min_{\g \in \bar{\rm{co}}\T_\e(y_0)}\int_G k(y,u) \g(dy,du) = \min_{\g \in W(\e, y_0)}\int_G k(y,u) \g(dy,du)
\end{equation*}
for any continuous $k(y,u)$ (see (\ref{auxiliary-1})). Since the validity of (\ref{P11}) has already been established (as a part of the proof of Theorem \ref{Th-discounting-1}), we only need to prove
that
\begin{equation}\label{eq-OPR-1-0}
\bar{\rm{co}}\T_\e(y_0) \supset W(\e, y_0)  \ \ \forall \e \in (0, 1).
\end{equation}
In fact, the latter can be shown to be true if
\begin{equation}\label{eq-OPR-0-0-1}
\min_{\g \in \bar{\rm{co}}\T_\e(y_0)}\int_G k(y,u) \g(dy,du) \leq \min_{\g \in W(\e, y_0)}\int_G k(y,u) \g(dy,du)
\end{equation}
for any continuous $k(y,u)$.
Assume that for some $\e \in (0,1)$ the inclusion (\ref{eq-OPR-1-0}) is not true. That is, there exists $\gamma'\in W(\e, y_0) \setminus \bar{\rm{co}}\T_\e(y_0) $.
By the separation theorem (see, e.g., \cite[p. 59]{Rudin}), there exists $\hat{k}(y,u) \in C(G)$ such that
\begin{equation*}
\min_{\g \in \bar{\rm{co}}\T_\e(y_0)} \int_G \hat{k}(y,u) \g(dy,du) \ge \int_G \hat{k}(y,u) \gamma'(dy,du) + \beta\ge  \min_{\g \in W(\e, y_0)}\int_G k(y,u) \g(dy,du) +\beta
\end{equation*}
for some $\beta > 0$. This contradicts (\ref{eq-OPR-0-0-1}) considered with $k(y,u)=\hat{k}(y,u) $. Thus, (\ref{eq-OPR-1-0}) is proved.

Let us now show that (\ref{convergence-to-W-dis-1}) is true  if (\ref{convergence-to-W-dis-4}) is true, that is, if
\begin{equation}\label{eq-OPR-0-2}
\lim_{\e \rightarrow 0} \min_{\g \in \bar{\rm{co}}\T_\e} \int_G k(y,u)\g (dy,du) = \min_{\g \in W} \int_G k(y,u) \g(dy,du)	
\end{equation}
for any continuous $k(y,u)$ (see  (\ref{auxiliary-2})).
Note that from (\ref{e-Theta-1}) it follows that
\begin{equation*}
\lim_{\e \rightarrow 0} \max_{\g \in \bar{\rm{co}}\T_\e} \rho(\g, W) = 0.
\end{equation*}
Therefore, to prove (\ref{convergence-to-W-dis-1}), we only need to show that
\begin{equation}\label{eq-OPR-1-2}
\lim_{\e \rightarrow 0} \max_{\g \in W} \rho(\g, \bar{\rm{co}}\T_\e) = 0.
\end{equation}
Assume that (\ref{eq-OPR-1-2}) is not true. Then there exists a positive number $\alpha$ and sequences $\e_i > 0$, $\g_i \in W,$ $i = 1,2,\ldots,$ such that $\lim_{i \rightarrow \infty}\e_i = 0$ and
\begin{equation}\label{eq-OPR-2-2}
\lim_{i \rightarrow \infty}\rho(\g_i, \bar{\rm{co}}\T_{\e_i}) \ge \a.
\end{equation}
Due to the fact $W$ is compact in the metric $\rho$, and due to the fact (by Blaschke's selection theorem; see e.g., \cite{Klein}) the set of closed subsets of $P(G)$ is compact in the corresponding Hausdorff  metric $\rho_H$, the sequences $\{\g_i \}$ and $\{\bar{\rm{co}}\T_{\e_i} \}$ have partial limits. More specifically, there exists $\hat{\g} \in W$ and $\hat{\T} \subset P(G)$ such that, for some subsequence $\{i' \} \subset \{i \}$,
\begin{equation}\label{eq-OPR-3-2}
\lim_{i' \rightarrow \infty} \rho(\g_{i'},\hat{\g})=0, \quad \lim_{i' \rightarrow \infty} \rho_H(\bar{\rm{co}}\T_{\e_{i'}}, \hat{\T})=0.
\end{equation}
(Note that, being a limit in the Hausdorff  metric of a sequence of convex compact sets, the set $\hat{\T} $ is convex and compact.)
By passing to the limit in (\ref{eq-OPR-2-2}), one obtains that
\begin{equation*}
\rho(\hat{\g}, \hat{\T}) \ge \a \quad \implies \quad \hat{\g} \notin \hat{\T}.
\end{equation*}
Based on the separation theorem (see e.g. \cite[p. 59]{Rudin}), we may conclude that there exists $\hat{k}(y,u) \in C(G)$ such that
\begin{equation}\label{eq-OPR-4-2}
\min_{\g \in \hat{\T}}\int_{G} \hat{k}(y,u) \g(dy,du) \ge \int_{G} \hat{k}(y,u)\hat{\g}(dy,du) + \beta \ge \min_{\g \in W}\int_{G} \hat{k}(y,u)\g(dy,du)  + \beta
\end{equation}
for some $\beta > 0$.
From the second equality in (\ref{eq-OPR-3-2}), it follows that
\begin{equation*}
\lim_{i' \rightarrow \infty}\min_{\g \in \bar{\rm{co}}\T_{\e_{i'}}}\int_{G}\hat{k}(y,u)\g(dy,du) = \min_{\g \in \hat{\T}}\int_{G}\hat{k}(y,u)\g(dy,du).
\end{equation*}
Consequently, from (\ref{eq-OPR-4-2})  it follows that
\begin{equation*}
\min_{\g \in \bar{\rm{co}}\T_{\e_{i'}}}\int_{G}\hat{k}(y,u)\g(dy,du) \ge \min_{\g \in W}\int_{G} \hat{k}(y,u)\g(dy,du)   + \frac{\beta}{2}
\end{equation*}
for $i'$ large enough. This contradicts (\ref{eq-OPR-0-2}) considered with $k(y,u)=\hat k(y,u) $. Thus, the validity of (\ref{eq-OPR-1-2}) is established, and (\ref{convergence-to-W-dis-1}) is proved.

The fact that (\ref{convergence-to-W-dis-2}) is valid if (\ref{convergence-to-W-dis-5}) is true, that is, if
\begin{equation}\label{eq-OPR-0-1}
 \lim_{T \rightarrow \infty} \min_{\g \in \bar{\rm{co}}\G_T} \int_G k(y,u)\g (dy,du) = \min_{\g \in W} \int_G k(y,u) \g(dy,du)
\end{equation}
for any continuous $k(y,u)$, is established in a similar way.
$ \ \Box$

\section{Linear programming results and proofs}\label{Sub-Sec-Auxiliary-Res}

Let us note again that the conditions about the model introduced in Section \ref{Model-Prelim}
are assumed to be satisfied
in this section too (like everywhere  in the consideration above).

\subsection{Results and proofs referred to in Section \ref{Sec-Time-discounting}}\label{Sec-LP1}
Let $C^*(Y) $ stand for the space of continuous linear functionals  on $C(Y) $ and let (as in Section \ref{Sec-Strong-Duality}) $\mathcal{M}(G) $ stand for the space of finite signed measures defined on Borel subsets of $G$. For any $\e\in [0,1) $, define a linear operator $\mathcal{A}_{\e}(\cdot): \mathcal{M}(G)\mapsto\reals\times  C^*(Y)  $  as follows: for any $\gamma\in \M(G)$,
\begin{equation}\label{e-Duality-1-W}
\mathcal{A}_{\e}(\gamma)\BYDEF \left(\int_{G}\gamma(dy,du), \ a_{\e}^{\gamma}\right),
\end{equation}
where $a_{\e}^{\gamma}\in C^*(Y)$ is defined by the equation: $\ \forall \  \ph(\cdot)\in C(Y)$,
$$
 \ a_{\e}^{\gamma}(\ph) \BYDEF -\left\{
  \int_{G}\big((1-\e)(E\big[\ph(f(y,u,s))\big]-\ph(y))+\e(\ph(y_0)-\ph(y))\big)\gamma(dy,du)\right\}.
 $$
 In these notations, the sets $W(\e, y_0), \ W$ defined in (\ref{w-disc}) and (\ref{M17}) can be rewritten as follows
$$
 W(\e,y_0)= \big\{\gamma \in \mathcal{M}_+(G)\ : \ \mathcal{A}_{\e}(\gamma)= (1, {\bf 0})\big\},\ \ \ \ \  W=W(0,y_0)= \big\{\gamma \in \mathcal{M}_+(G)\ : \ \mathcal{A}_{0}(\gamma)= (1, {\bf 0})\big\},
 $$
where  ${\bf 0} $ stands for the zero element of $C^*(Y)$. Also,
problem (\ref{D1})  takes the forms
\begin{equation}\label{e-Duality-2-0-W}
\min_{\gamma \in W(\e, y_0)}\langle k, \gamma  \rangle\ = k^*(\e, y_0),
\end{equation}
where
 $\langle \cdot, \gamma  \rangle $ denotes the integral of the corresponding function over $\gamma$ (problem (\ref{M22}) taking the form (\ref{e-Duality-2-0-W})
 when $\e = 0$).
 Note that, for any $(\mu, \psi(\cdot))\in \reals\times C(Y) $,
$$
\langle \mathcal{A}_{\e}(\gamma), (\mu, \psi)\rangle =
\mu \int_{G}\gamma(dy,du) + a_{\e}^{\gamma}(\psi)
$$
$$
=  \int_{G}\big(\mu -  \big((1-\e)(E\big[\psi(f(y,u,s))\big]-\psi(y))+\e(\psi(y_0)-\psi(y))\big)\big)\gamma(dy,du).
$$

Define the linear operator  $\mathcal{A}_{\e}^*(\cdot): \reals\times C(Y) \mapsto C(G)\subset \mathcal{M}^*(G)$ in such a way that, for any $(\mu, \psi(\cdot))\in \reals\times C(Y) $,
\begin{equation}\label{e-Duality-3-1}
\mathcal{A}_{\e}^*(\mu, \psi)(y,u)\BYDEF \mu  - \big((1-\e)(E\big[\psi(f(y,u,s))\big]-\psi(y))+\e(\psi(y_0)-\psi(y))\big).
\end{equation}
Thus,
$$
\langle  \mathcal{A}_{\e}^*(\mu, \psi), \gamma \rangle =  \int_{G}\big(\mu -  \big((1-\e)(E\big[\psi(f(y,u,s))\big]-\psi(y))+\e(\psi(y_0)-\psi(y))\big)\big)\gamma(dy,du)
 = \langle \mathcal{A}_{\e}(\gamma), (\mu, \psi )\rangle  .
$$
That is, the operator $ \mathcal{A}_{\e}^*(\cdot) $ is the adjoint of  $ \mathcal{A}_{\e}(\cdot) $. The problem dual to (\ref{e-Duality-2-0-W})
is of the form (see \cite{And-1} and \cite{And-2})
\begin{equation}\label{e-Duality-4-0}
\sup_{(\mu, \psi(\cdot))\in \reals\times C(Y)} \mu = \mu^*(\e, y_0)
\end{equation}
$$
\ s.\ t.  \ \ \ \ \ \ \ \ \ \ \   \ \ \ \ \ \ \ \ \ \ \ \ \ \ \ \ \ \ \ \   \ \ \ \ \ \ \ \ \
$$
\begin{equation}\label{e-Duality-4-0-1}
 - \mathcal{A}_{\e}^*(\mu, \psi)(y,u) + k(y,u))\geq 0 \ \ \forall (y,u)\in G .
\end{equation}
Due to (\ref{e-Duality-3-1}), the constraint (\ref{e-Duality-4-0-1}) is equivalent to
$$
k(y,u)+(1-\e)(E\big[\psi(f(y,u,s))\big]-\psi(y))+\e\big(\psi(y_0)-\psi(y)\big)\geq \mu \ \ \forall (y,u)\in G .
$$
Thus, problem (\ref{e-Duality-4-0}) is equivalent to (\ref{e-dis-dual}) for $\e\in (0,1)$, and it is equivalent to (\ref{M21}) when $\e=0$.

\medskip

{\bf Proof of Proposition \ref{SD-1}.} By Theorem 6 in \cite{And-1}, to prove the proposition, it is sufficient to establish
that, for any $\e\in [0,1) $, the set $D_{\e}$,
$$
D_{\e}\BYDEF \big\{ \left(\mathcal{A}_{\e}(\gamma ),\ \langle k, \gamma  \rangle\ \right) \ : \ \gamma \in
{\cal M}_+(G)\big\}\subset \reals\times C^*(Y)\times \reals ,
$$
is closed in the weak$^*$ topology of $ \reals\times C^*(Y)\times \reals$. The proof of this is
similar to the proof of Theorem 12 in \cite{And-1}. It is based on the following two properties of the problem.

\noindent{\it Property A.}  The set  ${\cal M}_+(G) $ has a compact base. That is (see {\rm \cite{And-1}}),
\begin{equation}\label{Lemma-Anderson-2}
{\cal M}_+(G) = \{\lambda \gamma  : \ \gamma \in \mathcal{L}, \ \lambda\geq 0  \},
\end{equation}
where
$$
\mathcal{L}\BYDEF \{\gamma \in  {\cal M}_+(G) \ :
\int_{G} \gamma(dy,du)   = 1 \},
$$
with $\mathcal{L}$ being a weak$^*$ compact subset of ${\cal M}(G)$.

\medskip

\noindent{\it Property B.}  For $\gamma\in  {\cal M}_+(G)$, the equalities
\begin{equation}\label{Lemma-Anderson-3}
\mathcal{A}_{\e}(\gamma ) = (0, {\bf 0})
\end{equation}
can be valid only if $\gamma = 0 $, this being readily verifiable due to the fact that a part of the relationships  {\rm (\ref{Lemma-Anderson-3})} is the equality:
$$
\int_{G} \gamma(dy,du) = 0.
$$

Let us now prove that $D$ is closed. Let $\gamma_l\in {\cal M}_+(G)$ be such that
\begin{equation}\label{Lemma-Anderson-4}
\mathcal{A}_{\e}(\gamma_l ) - {\bf z}\rightarrow (0, {\bf 0}), \ \ \ \ \langle k,  \gamma_l \rangle\ - \beta\rightarrow 0 \ \ \ {\rm as} \ \ \ l\rightarrow\infty,
\end{equation}
where
${\bf z}\in  \reals\times C^*(Y)  $ and $\beta\in \reals$. By (\ref{Lemma-Anderson-2}), $\gamma_l  = \lambda_l \bar \gamma_l  $, where $\lambda_l\geq 0 $ and $\ \bar \gamma_l \in \mathcal{L} $. Due to compactness of $\mathcal{L} $, one may assume (without loss of generality) that  $\ \bar \gamma_l \rightarrow \bar \gamma
\in \mathcal{L} $.

Note that, due to Property B, the sequence $\lambda_l$ is bounded. Indeed, assuming that this is not the case and there exists a subsequence $\{l'\} $ of $\{l\} $ such that $\lambda_{l'}\rightarrow\infty $ as $l'\rightarrow\infty $ , one would obtain (via substitution of $\ \lambda_{l'} \bar \gamma_{l'}  $ into (\ref{Lemma-Anderson-4}) and passing to the limit with $l'\rightarrow\infty $) that
$$
\ \mathcal{A}_{\e}(\bar \gamma_{l'} )- \frac{1}{\lambda_{l'}}{\bf z}\rightarrow (0, {\bf 0}), \ \ \
  \langle k, \bar \gamma_{l'} \rangle\ -
\frac{1}{\lambda_{l'}}\beta\ \rightarrow
0 \ \ \ {\rm as} \ \ \ l\rightarrow\infty.
$$
$$
\Rightarrow \ \ \ \ \ \mathcal{A}_{\e}(\bar\gamma) = (0, {\bf 0}), \ \ \ \ \langle k, \bar\gamma  \rangle = 0.
$$
According to Property B, the latter implies that $\bar\gamma = 0 $. This contradicts to the fact that $\bar \gamma \in \mathcal{L} $. Thus, the sequence $\{\lambda_l\} $ is bounded and therefore one may assume (without loss of generality) that $\lambda_l\rightarrow \lambda  $ as $l\rightarrow\infty $. Consequently,  $\gamma_l \rightarrow \lambda \bar\gamma  $
as $l\rightarrow\infty $, and, by (\ref{Lemma-Anderson-4}),
$$
\mathcal{A}_{\e}(\lambda \bar\gamma ) = {\bf z}, \ \ \ \langle k, \lambda \bar\gamma\rangle = \beta \ \ \ \ \Rightarrow \ \ \ \
({\bf z},\beta)\in D.
$$
This proves that $D$ is closed. $\ \Box$

 \medskip

 REMARK. Note that the proof above is similar to the proof of  the absence of the duality gap in \cite{PZ-2011}.

\subsection{Results and proofs referred to in Sections \ref{Sec-Augmented-estimates} and \ref{Sec-Strong-Duality}}\label{Sec-LP2}
Define a linear operator $\mathcal{A}(\cdot): \mathcal{M}(G) \times  \mathcal{M}(G)\mapsto\reals\times  C^*(Y) \times  C^*(Y) $  as follows: for any $(\gamma, \xi)\in \M(G)\times \M(G) $,
\begin{equation}\label{e-Duality-1}
\mathcal{A}(\gamma, \xi):= \left(\int_{G}\gm, \ a^{(\gamma, \xi)}, \ b^{\gamma}\right),
\end{equation}
where $a^{(\gamma, \xi)}, \ b^{\gamma}\in C^*(Y)$ are defined for all $\ph(\cdot)\in C(Y)$ as follows
$$
\ a^{(\gamma, \xi)}(\ph)\BYDEF -\left\{\int_{G} (\ph(y_0)-\ph(y))\gm
+ \int_{G}(E[\ph(f(y,u,s))] - \ph(y))\xm\right\},
$$
\vspace{-0.4cm}
$$
\  b^{\gamma}(\ph)\BYDEF - \left\{\int_{G}(E[\ph(f(y,u,s))] - \ph(y))\gm\right\}.
$$
In these notations, the set $\Omega(y_0)$ defined in (\ref{eq-Omega}) takes the form
$$
\Omega(y_0)= \big\{(\gamma , \xi)\in \mathcal{M}_+(G)\times \mathcal{M}_+(G)\ : \ \mathcal{A}(\gamma, \xi)= (1, {\bf 0}, {\bf 0})\big\},
$$
and problem (\ref{BB1})  can be rewritten as follows
\begin{equation}\label{e-Duality-2-0}
\inf_{(\gamma , \xi)\in \Omega(y_0)}\langle k, \gamma  \rangle\ = k^*(y_0).
\end{equation}
For any $(\mu, \psi(\cdot), \eta(\cdot))\in \reals\times C(Y)\times C(Y) $, we have
\begin{equation*}
\begin{split}
\langle \A(\gamma, \xi), (\mu, \psi, \eta)\rangle &=
\mu \int_{G}\gm + a^{(\gamma, \xi)}(\psi)
+ b^{\gamma}(\eta)\\
&=  \int_{G}\left(\mu - (\psi(y_0)-\psi(y)) - ( E[\eta(f(y,u,s))] - \eta(y))\right)\gm\\
& \ \ \ - \int_{G} ( E[\psi(f(y,u,s))] - \psi(y))\xm.
\end{split}
\end{equation*}
Define now the linear operator $\mathcal{A}^*(\cdot): \reals\times C(Y)\times C(Y) \mapsto C(G)\times C(G)\subset \M^*(G)\times  \M^*(G) $ in such a way that, for any $(\mu, \psi(\cdot),\eta(\cdot))\in \reals\times C(Y)\times C(Y) $,
\begin{equation}\label{e-dual-100}
\mathcal{A}^*(\mu, \psi, \eta)(y,u)\BYDEF \big(\mu - (\psi(y_0)-\psi(y)) - ( E[\eta(f(y,u,s))] - \eta(y)), \ -( E[\psi(f(y,u,s))] - \psi(y))\big).
\end{equation}
Thus,
$$
\langle  \mathcal{A}^*(\mu, \psi, \eta), (\gamma, \xi) \rangle =  \int_{G}\left(\mu - (\psi(y_0)-\psi(y)) - ( E[\eta(f(y,u,s))] - \eta(y))\right)\gm
$$
$$
- \int_{G} ( E[\psi(f(y,u,s))] - \psi(y))\xm = \langle A(\gamma, \xi), (\mu, \psi, \eta) \rangle  .
$$
That is, the operator $ \mathcal{A}^*(\cdot) $ is the adjoint of  $ \mathcal{A}(\cdot) $. Thus, the problem dual to (\ref{e-Duality-2-0}) is of the form (see e.g. \cite{And-1} or \cite{And-2})
\begin{equation}\label{e-dual-101}
\sup_{(\mu, \psi(\cdot), \eta(\cdot))\in \reals\times C(Y)\times C(Y)} \mu = d^*(y_0)
\end{equation}
such that
\begin{equation}\label{e-dual-102}
 \mathcal{A}^*(\mu, \psi, \eta)(y,u) + (k(y,u), 0)\geq (0,0) \ \ \  \forall \ (y,u)\in G.
\end{equation}
In accordance with (\ref{e-dual-100}), the constraint (\ref{e-dual-102}) is equivalent to (\ref{BB7}), and problem
(\ref{e-dual-101}) is equivalent to 

\smallskip

(\ref{BB8}).

\medskip

{\bf Proof of Proposition \ref{Prop-cl-D-y-0-dual}. }
Let
$$
H\BYDEF\Big\{\Big( \A(\gamma, \xi), \int_{Y\times U}k(y,u)\gamma(dy,du) + r\Big)\ : \ (\gamma, \xi)\in \M_+(G)\times \M_+(G), \ r\geq 0\Big\}
$$
\vspace{-0.3cm}
$$
\subset \reals\times C^*(Y)  \times C^*(Y) \times \reals\  ,
$$
and let $\bar{H}$ stand for the closure of $H$ in the weak$^*$ topology of the space $\reals \times C^*(Y) \times C^*(Y) \times \reals$. That is,
$(\kappa , a , b, \theta) $ belongs to $\bar{H}$ (where $a , b\in C^*(Y) $ and $  \kappa , \theta \in \reals  $)  if and only if there exists  sequences $ (\gamma_l, \xi_l)\in \M_+(G)\times \M_+(G) $, \ $r_l\geq 0$, $l=1,2,...$,
such that
$$
a^{(\gamma_l, \xi_l)}(\ph)\to a(\ph),\ \ \ \ \ \ \ b^{\gamma_l  }(\ph)\to b(\ph)\ \ \ \ \ \forall \ \ph(\cdot)\in C(Y)
$$
and
$$
\int_{G}\gamma_l(dy,du)\to \kappa, \ \ \  \ \ \ \int_{Y\times U}k(y,u)\gamma_l(dy,du) + r_l\to \theta .
$$
Consider the problem
\begin{equation}\label{limits-non-ergodic-pert-dual-4-1}
\inf \{\theta \ | \ (1, {\bf 0}, {\bf 0}, \theta)\in \bar{H} \}\BYDEF k_{sub}^*(y_0).
\end{equation}
Its optimal value $k_{sub}^*(y_0)$ is called the subvalue of the IDLP problem (\ref{e-Duality-2-0}) (see \cite{And-1}, \cite{And-2}).
\smallskip
Let us show that the optimal value of (\ref{BB1-D-y-0-2}) is equal to the subvalue.

Firstly note that,
as can be readily seen, $\left(1, {\bf 0}, {\bf 0},  \int_{Y\times U}k(y,u)\gamma(dy,du)\right)\in \bar{H} $  if $\gamma\in W\cap D_1(y_0)$. Consequently,
$$
k_{sub}^*(y_0)\leq \min_{\gamma\in W\cap D_1(y_0)}\int_{Y\times U}k(y,u)\gamma(dy,du)=k^{**}(y_0).
$$
From the fact that $k_{sub}^*(y_0)$ is defined as  the optimal value in (\ref{limits-non-ergodic-pert-dual-4-1}) it follows that there exists  a sequence $(\gamma_l,\xi_l)\in \M_+(G)\times \M_+(G) $  such that $\A(\gamma_l,\xi_l) $ converges (in weak$^*$ topology)
to $(1, {\bf 0}, {\bf 0})$, with   $\int_{Y\times U}k(y,u)\gamma_l(dy,du) $ converging to $k_{sub}^*(y_0)$ as $l$ tends to infinity. That is (see (\ref{e-Duality-1})),
$$
\int_{G}\gamma_l(dy,du)\rightarrow 1, \ \ a_{(\gamma_l, \xi_l)}\rightarrow {\bf 0},\ \  b_{\gamma_l}\rightarrow {\bf 0},
$$
 \vspace{-0.4cm}
$$
\ \ \int_{G}k(y,u)\gamma_l(dy,du)\rightarrow k_{sub}^*(y_0).
$$
Without loss of generality, one may assume
that $\gamma_l$ converges in weak$^*$ topology
 to a  measure $\bar \gamma$ that satisfies the relationships
 $$
 \int_{G}\bar \gamma(dy,du)=1,  \ \  b_{\bar \gamma}= {\bf 0}\ \ \ \ \Rightarrow \ \ \ \  \bar \gamma\in W.
 $$
 Also,
 $
 \  a_{(\bar\gamma, \xi_l)}\rightarrow {\bf 0}
 $ and $\int_{G}k(y,u)\bar\gamma(dy,du)= k_{sub}^*(y_0) $. From the fact that $
 \  a_{(\bar\gamma, \xi_l)}\rightarrow {\bf 0}
 $ it follows that $\bar \gamma\in D_1(y_0)$ (see (\ref{BB1-D-y-0-2-1})). That is,  $\ \bar \gamma\in W\cap D_1(y_0) $. Consequently,
 $$
k^{**}(y_0)=\min_{\gamma\in W\cap D_1(y_0)}\int_{G}k(y,u)\gamma(dy,du)\leq \int_{G}k(y,u)\bar\gamma(dy,du)  =  k_{sub}^*(y_0).
$$
Thus, $\ k^{**}(y_0)= k_{sub}^*(y_0)$.
 To complete the proof, it is sufficient to note that the subvalue of an IDLP problem 
 
 \smallskip
 
 is equal to the optimal value of its dual  provided that the former is bounded
  (see, e.g., Theorem 3 in \cite{And-1}). That is, $k_{sub}^*(y_0) = d^*(y_0). $
$ \ \Box$

\medskip

Consider  the following \lq\lq perturbed" version of the IDLP problem (\ref{BB1}):
\begin{equation}\label{BB1-Per}
\inf_{(\g,\xi)\in \O(y_0)} \left\{\int_G k(y,u)\gm +  \int_G \theta(y,u) \xi(dy,du)\right\} \BYDEF k^*_{\theta}(y_0),
\end{equation}
where $\theta(y,u)\geq 0$ is a bounded Borel measurable  function on $Y\times U$. Using an argument similar to the one used above, one can verify that the problem dual
to (\ref{BB1-Per}) can be written in the form
\begin{equation}\label{BB8-Per}
\sup_{(\mu,\psi,\eta)\in \D_{\theta}(y_0)} \mu \BYDEF d^*_{\theta}(y_0),
\end{equation}
where $\D_{\theta}(y_0)$ is the set of triplets $(\mu,\psi(\cdot),\eta(\cdot))\in \reals\times C(Y)\times C(Y)$  that satisfy the inequalities
\begin{equation}\label{BB7-Per}
\begin{aligned}
&k(y,u)+(\psi(y_0)-\psi(y))+E[\eta(f(y,u,s))]-\eta(y)-\mu \ge 0,\\
&E[\psi(f(y,u,s))]-\psi(y)\ge -\theta(y,u)  \ \ \ \ \forall \ (y,u)\in G.
\end{aligned}
\end{equation}
(Note that  (\ref{BB1-Per}) and (\ref{BB8-Per}) coincide with (\ref{BB1}) and (\ref{BB8}) if $\ \theta(y,u) \equiv 0 $.)
Consider also the problem
\begin{equation}\label{BB21-1-Per}
\sup_{(\psi,\eta)\in Q_{\theta}} \psi(y_0)\BYDEF \bar d^*_{\theta}(y_0),
\end{equation}
where
$Q_{\theta}$ is the set of pairs $(\psi(\cdot),\eta(\cdot))\in C(Y)\times C(Y)$ that satisfy the inequalities
\begin{equation}\label{BB25-1-Per}
\begin{aligned}
&k(y,u)-\psi(y)+E[\eta(f(y,u,s))]-\eta(y)\ge 0,\\
&E[\psi(f(y,u,s))]-\psi(y)\ge -\theta(y,u)\ \ \ \ \forall \ (y,u)\in G.
\end{aligned}
\end{equation}

\begin{Lemma}\label{L-equivalince}
	The following relationships are valid:
	\begin{equation}\label{e-Pert-1}
	\bar d^*_{\theta}(y_0) =  d^*_{\theta}(y_0)\leq k^*_{\theta}(y_0) .
	\end{equation}
\end{Lemma}
{\bf Proof.} Let us prove, first, that
\begin{equation}\label{e-Pert-2}
\bar d^*_{\theta}(y_0) = d^*_{\theta}(y_0)  .
\end{equation}
In fact, the inequality $\bar d^*_{\theta}(y_0) \leq  d^*_{\theta}(y_0) $ is true (since, for any pair $\ (\psi(\cdot),\eta(\cdot))\in Q_{\theta} $, the triplet $\ (\mu , \psi(\cdot),\eta(\cdot))\in \D_{\theta}(y_0) $ with $\mu = \psi(y_0) $). Let us prove the opposite inequality. Let a triplet $\ (\mu' , \psi'(\cdot),\eta'(\cdot))\in \D_{\theta}(y_0) $ be such that $\mu'\geq d^*_{\theta}(y_0)-\beta $, with $\beta > 0 $ being arbitrarily small. Then the pair $\ (\tilde \psi'(\cdot), \eta'(\cdot))\in Q_{\theta} $, with $\tilde\psi'(y)= \psi'(y)- \psi'(y_0) + \mu'$. Since $ \tilde\psi'(y_0) = \mu'$, it leads to the inequality $\bar d^*_{\theta}(y_0)\geq \mu' \geq  d^*_{\theta}(y_0) - \beta $ and, consequently, to the inequality $\bar d^*_{\theta}(y_0) \geq  d^*_{\theta}(y_0) $. Thus, (\ref{e-Pert-2}) is proved.

Let us now prove the weak duality inequality
\begin{equation}\label{e-Pert-3}
d^*_{\theta}(y_0)\leq k^*_{\theta}( y_0).
\end{equation}
Take any $(\g,\xi)\in \O(y_0)$ and $(\mu,\psi,\eta)\in \D_{\theta}(y_0)$. Integrating the first inequality in \eqref{BB7-Per} with respect to $\g$ and taking into account that $\g\in W$ we conclude that
$$
\int_G k(y,u)\gm+\int_G(\psi(y_0)-\psi(y))\gm\ge \mu.
$$
Taking into account that $(\g,\xi)\in \O(y_0)$ and the second inequality in \eqref{BB7-Per}, we obtain
$$
\int_G(\psi(y_0)-\psi(y))\gamma(dy,du)=-\int_G (E[\psi(f(y,u,s))]-\psi(y))\xi(dy,du)\le \int_{G}\theta(y,u)\xi(dy,du).
$$
Therefore,
$$
\int_G k(y,u)\gm + \int_{G}\theta(y,u)\xi(dy,du) \ge \mu .
$$
This proves (\ref{e-Pert-3}). \hf

\section{Acknowledgement}
We would like to thank Yi Zhang from the University of Liverpool for useful comments on the present work.
We would also like to thank two anonymous reviewers for their constructive comments and, in particular the reviewer that suggested an alternative (much shorter than the original one) proof of Lemma \ref{Lem-upper-bound}.

\bibliography{bibliography}
\bibliographystyle{elsarticle-harv}

\end{document}